\documentclass[12pt]{amsart}
\usepackage{longtable}
\usepackage{amsfonts,amssymb,amscd,amsmath,epsfig}
\usepackage{latexsym}
\usepackage{mathrsfs}
\usepackage{tocvsec2}

\usepackage[%dvipdfm, 
hyperindex=true,pagebackref=true,bookmarks=true,
colorlinks=true,linkcolor=blue,citecolor=red]{hyperref}
\def\~{{\rm --}} 
\font\dfont=cmbx10 at 11pt   %%% FOR DEFIN

\title [Nonsymmetric Rogers-Ramanujan sums]
{Nonsymmetric Rogers-Ramanujan sums and 
thick Demazure modules}
\author[Ivan Cherednik]{Ivan Cherednik $^\dag$}
\author[Syu Kato]{Syu Kato $^\ddag$}
\thanks{$^\dag$ \today.
\ \ \ Partially supported by NSF grant
DMS--1363138}

\thanks{\noindent $^\ddag$ Supported in part by JSPS Grants-in-Aid 
for Scientific Research (B) JP26287004}

\address[I. Cherednik]{Department of Mathematics, UNC
Chapel Hill, North Carolina 27599, USA\\
chered@email.unc.edu}

\address[S. Kato]{ Department of Mathematics, Kyoto University,
Oiwake Kita-Shirakawa/ Sakyo Kyoto 606-8502, Japan\\
syuchan@math.kyoto-u.ac.jp}

 \def\bysame{{\bf --- }}
 \def\~{{\bf --}}
\newcommand{\comment}[1]{}

\renewcommand{\tilde}{\widetilde}
\renewcommand{\hat}{\widehat}

\newcommand{\malt}{\hbox{\tiny$\!{}^\maltese$}}

\newcommand{\ee}{e}
\newcommand*{\medoplus}{\mathbin{\scalebox{.75}
{\ensuremath{\bigoplus}}}}
\newcommand{\extra}[1]{} 
\renewcommand{\extra}[1]{{\color{red}\sf#1}}  %%ONLY FOR REMARKS

\renewcommand{\tilde}{\widetilde}
\renewcommand{\hat}{\widehat}

\newcommand{\Z}{{\mathbb Z}}
\newcommand{\Q}{{\mathbb Q}}
\newcommand{\N}{{\mathbb N}}
\newcommand{\C}{{\mathbb C}}
\newcommand{\R}{{\mathbb R}}

\def\HH{\mbox{${\mathcal H}$\kern-5.2pt${\mathcal H}$}}

\newtheorem{theorem}{Theorem}[section]

\newtheorem{proposition}[theorem]{Proposition}

\newtheorem{lemma}[theorem]{Lemma}
\newtheorem{corollary}[theorem]{Corollary}

\newtheorem{theorem }{Theorem}[section]
\newtheorem{maintheorem }[theorem]{Main Theorem}
\newtheorem{proposition }[theorem]{Proposition}
\newtheorem{definition }[theorem]{Definition}
\newtheorem{lemma }[theorem]{Lemma}
\newtheorem{corollary }[theorem]{Corollary}
\newtheorem{notation }[theorem]{Notation}
\newtheorem{remark }[theorem]{Remark}
\newtheorem{example }[theorem]{Example}

\newtheorem{ maintheorem }[theorem]{Main Theorem}
\newtheorem{ theorem}{Theorem}[section]
\newtheorem{ proposition}[theorem]{Proposition}
\newtheorem{ definition}[theorem]{Definition}
\newtheorem{ lemma}[theorem]{Lemma}
\newtheorem{ corollary}[theorem]{Corollary}
\newtheorem{ notation}[theorem]{Notation}
\newtheorem{ remark}[theorem]{Remark}
\newtheorem{ example}[theorem]{Example}

\newtheorem{thm}{Theorem}[section]
\newtheorem{prop}[thm]{Proposition}
\newtheorem{lem}[thm]{Lemma}
\newtheorem{cor}[thm]{Corollary}
\newtheorem{conj}[thm]{Conjecture}

\def\for{\  \hbox{ for } \ }
\def\iif{ \ \hbox{ if } \ }

\def\where{\  \hbox{ where } \ }
\def\and{\  \hbox{ and } \ }
\def\and{\  \hbox{ and } \ }

\def\equal{\stackrel{\,\mathbf{def}}{= \kern-3pt =}}

\def\la{\lambda}
\def\La{\Lambda}
\def\om{\omega}

\def\th{\theta}
\def\al{\alpha}

\def\de{\delta}

\def\si{\sigma}
\def\Si{\Sigma}
\def\Ga{\Gamma}
\def\ze{\zeta}

\def\vth{{\vartheta}}

\def\tal{\tilde{\alpha}}

\def\tga{\tilde{\gamma}}
\def\tGa{\tilde{\Gamma}}

\def\tw{\widetilde w}
\def\tW{\widetilde W}

\def\tz{\tilde z}
\def\tb{\tilde b}

\def\tR{\tilde R}

\def\hw{\widehat{w}}
\def\hW{\widehat{W}}

\def\0{\mathbf{0}}

\def\çF{\mathcal{F}}

\def\e{\mathcal{E}}
\def\v{\mathcal{V}}

\def\b{\mathcal{B}}

\def\lan{\langle}
\def\llb{(\!(}
\def\ran{\rangle}
\def\rrb{)\!)}
 \def\dim{{\hbox{\rm dim}}_{\mathbb C}\,}
\def\lng{\hbox{\rm{\tiny lng}}}
\def\sht{\hbox{\rm{\tiny sht}}}
\newcommand{\sq}{\phantom{1}\hfill$\qed$}

\def\HH{\mathfrak{H}}

\def\CC{\mathfrak{C}}
\def\LL{\mathfrak{L}}

\def\HH{\hbox{${\mathcal H}$\kern-5.2pt${\mathcal H}$}}
\def\HHH{\hbox{${\mathbb H}$\kern-4.2pt${\mathbb H}$}}

\font\smm=msbm10 at 12pt 
\def\symbol#1{\hbox{\smm #1}}
\def\lsmash{{\symbol n}}

\def\#{\sharp}

\begin{document}
\begin{abstract}
We consider expansions of products of theta-functions 
associated with arbitrary root systems in terms of nonsymmetric 
Macdonald polynomials at $t=\infty$ divided by their norms. 
The latter are identified with the graded
characters of Demazure slices, 
some canonical quotients of thick (upper) level-one Demazure 
modules, directly related to recent theory of generalized 
(nonsymmetric) global Weyl modules.
The symmetric Rogers-Ramanujan-type series 
considered by Cherednik-Feigin were expected to have
some interpretation of this kind; the nonsymmetric setting 
appeared necessary to achieve this. As an application, 
the coefficients of the nonsymmetric Rogers-Ramanujan
series provide 
formulas for the multiplicities of the expansions of tensor 
products of level-one Kac-Moody representations in terms 
of Demazure slices.
\end{abstract}

\maketitle {\em\small Key words: Rogers-Ramanujan series; 
Hecke algebras; nonsymmetric Macdonald polynomials; Kac-Moody 
algebras; Demazure modules; global Weyl modules.} 

\renewcommand{\baselinestretch}{1.2} 
{\textmd \tableofcontents} 
\renewcommand{\baselinestretch}{1.0}
%} 
%$\mathfrak{a,b,c,d,e,f,g,h,i,j,k,l,m,n,o,p,q,r,s,t,e^u,e^v,e^w}$
\vfill\eject 

\renewcommand{\natural}{\wr}

\setcounter{section}{0} \setcounter{equation}{0} 
\section{\sc Introduction}
Generalizing \cite{ChFB}, we expand the products of (standard)
theta-functions associated with arbitrary root systems in terms 
of {\em nonsymmetric\,} Macdonald polynomials at $t=\infty$ 
divided by their norms. The latter are identified with
the graded characters of {\em Demazure slices\,}, canonical
quotients of the level-one {\em thick\,} (upper) 
Demazure modules, which are directly related to 
recent theory of {\em generalized global Weyl modules\,}; see 
\cite{FKM,FMO,KL} and references there. We conjecture and 
partially prove (in the ``mixed case") that the expansions above,
the {\em nonsymmetric Rogers-Ramanujan series\,}, 
match the decomposition of tensor products of level-one 
integrable representations of Kac-Moody algebras in terms 
of Demazure slices. This provides formulas for the 
corresponding multiplicities, which are generally difficult 
to obtain representation-theoretically.

\subsection{\bf Major objectives}
Let us briefly overview the problems we approach and (partially) 
settle and the major techniques that are used. 
The identification of nonsymmetric Macdonald polynomials 
at $t=\infty$ divided by their norms with the graded
characters of Demazure slices is the key; interestingly,
the DAHA-based expansion of a (single) theta-function
is used in our proof. Applications of this fact to tensor 
products of arbitrary integrable Kac-Moody modules will require 
further efforts. We only consider the tensor products of
level-one modules. 
%We were able to 
%interpret representation-theoretically only some of the 
%expansions of the products of theta-functions we obtain in the 
%paper. 

\subsubsection{\sf Rogers-Ramanujan summations}
The expansions of products of standard theta-functions 
associated with an arbitrary (simple) root system in terms of 
$q$\~Hermite polynomials from \cite{ChFB} generalize the 
celebrated Rogers-Ramanujan summations and its various 
multi-rank extensions.  There are connections with 
\cite{An,War,GOW} and a vast literature on the Rogers-Ramanujan 
identities, the representation theory of Kac-Moody algebras, and 
related mathematical physics, including the so-called Nahm 
conjecture; see \cite{VZ, ChFB, CGZ}. Here the $q$\~Hermite 
polynomials, the symmetric Macdonald polynomials at $t=0$, were
sufficient. The first usage of the $q$\~Hermite 
polynomials in this context is actually due to Rogers. 

\subsubsection{\sf Demazure modules} 
When we take a single copy of the theta-function, the 
expansion in \cite{ChFB} was expected to be related to the 
filtrations of level-one representations of affine Lie algebra 
in terms of thick Demazure modules, as in 
\cite{Kas94}, but the representation-theoretical tools for 
making this a theorem were developed only recently.
\vskip 0.2cm 

Thick Demazure modules are infinite dimensional in general, as 
opposed to the {\em thin\,} (lower) Demazure modules that are 
always finite-dimensional; see e.g. \cite{Kum}. It is well known 
that thin Demazure modules of level one and nonsymmetric 
Macdonald polynomials at $t = 0$ are closely related through the 
{\em local Weyl modules\,} of (twisted) current algebras 
\cite{San,Ion1,CL,FL,FK}. Moreover, the  symmetric $q$\~Hermite 
polynomials divided by their norms were interpreted as the {\em 
global Weyl modules} \cite{LNSSS,FMS}, with some technical 
reservations. 
 Furthermore, the filtration of parabolic 
Verma modules in terms of the global Weyl modules was provided 
in \cite{CI}, which is important to us.
The latter filtration was generalized to any 
integrable highest weight modules of affine Lie algebra in 
\cite{KL}. Namely, the existence of the filtration of any such 
modules in terms of global Weyl modules was proven there. \vskip 
0.2cm 

This theory combined with \cite{FKM,FMO,OS} is essentially 
sufficient to connect {\em thick\,} Demazure modules with 
global and generalized global Weyl modules, at least for the 
types $ADE$. We use this approach and develop a systematic 
theory of {\em Demazure-Joseph functor\,} 
\cite{J} in the thick case. 
One of the applications is a connection between thin and thick 
Demazure modules; also, any twisted root systems can be 
considered. From the {\em nil-DAHA\,} perspective, the 
Demazure-Joseph 
functor is closely related to  the DAHA intertwiners from 
\cite{ChO}, which is one of the key advantages of the usage of
nonsymmetric Macdonald polynomials vs. the symmetric ones. 

%Lenart-Naito-Sagaki-Schilling-Shimozono 
%Fourier-Manning-Senesi. Chari-Ion
\subsubsection{\sf Nonsymmetric setting}
In order to fully employ these and other techniques, it is 
necessary to consider {\em all\,} Demazure modules, not only 
those stable under the action of the classical part of the 
affine Lie algebra (i.e. in the case of dominant weights). This 
corresponds to the passage from the symmetric to nonsymmetric 
Macdonald polynomials in DAHA theory. Here the limits $t\to 0$ 
and $t\to \infty$ result in very different families of 
nonsymmetric polynomials. The nonsymmetric $q$\~Hermite 
polynomials, called {\em $\overline{E}$\~polynomials} in this 
paper, correspond to $t\to 0$; they are significantly simpler 
than the {\em $E^\dag$\~polynomials} corresponding to $t\to 
\infty$. A direct relation is only in the direction from 
$E^\dag$ to $\overline{E}$. However $E^\dag$\~polynomials are 
dual to $\overline{E}$\~polynomials with the respect to some
natural inner product; this is important to us. 

 The goal of this paper is to 
present, examine and apply three explicit expansion formulas for 
the product of standard theta-functions in terms of  
$\overline{E}$\~ and $E^\dag$\~polynomials (divided by their 
norms). All three coincide with the corresponding expansion from 
\cite{ChFB} upon the {\em symmetrization\,}; the case when
we expand in terms of (nonsymmetric)
$\overline{E}$\~polynomials is actually close to 
the (symmetric) one from \cite{ChFB}. 

A general problem is actually
to expand products of the theta-functions 
multiplied by any $\overline{E}$\~ or $E^\dag$\~polynomials. 
Such formulas of course require nonsymmetric theory.
These $E$\~factors  must be omitted 
(unless for minuscule weights) if one wants to obtain 
Rogers-Ramanujan-type series of {\em $PSL(2,\Z)$\~modular type\,}, 
as in \cite{VZ,ChFB,CGZ};
their presence destroys the modularity. When 
they are absent,  we therefore expand {\em the same\,} products 
of theta-functions as in \cite{ChFB} and our formulas are
(non-trivial) {\em partitions\,} of those considered there. 
The main point of this paper is that (even without the
$E$\~factors) such partitions have deep representation-theoretic 
meaning. For instance, the DAHA intertwiners,
which require the nonsymmetric setting, are closely connected 
with the {\em Demazure-Joseph functors} \cite{J}.

\subsection{\bf Main results}
\subsubsection{\sf \texorpdfstring{{\mathversion{bold}$E$\~}}
{E-} polynomials} Let $R=\{\al\}\subset \R^n$ be a simple root 
system, $\{\al_i\}$ simple roots, $W$ the Weyl group, $P,Q$ the 
weight and root lattices. Let  $b^-$ denote the antidominant 
element (i.e. that from $P_-$) that is $W$\~conjugate to $b\in 
P$. The normalization of the standard form in $\R^n$ is $( 
\alpha, \alpha ) = 2$ for {\em short\,} roots; the corresponding 
affine system is {\em twisted\,}: 
$$
\tilde{R}=\{[\al,\nu_\al j]\mid 
\al\in R, \, j\in \Z,\, \nu_\al=(\al,\al)/2\}.
$$

See Section \ref{defNSMac} for the definition of the 
nonsymmetric Macdonald polynomials $E_b(X;q,t)$ $(b\in P)$. By 
${}^\ast$, we denote the standard DAHA anti-involution sending 
$X_b\mapsto X_b^{-1}, q\mapsto q^{-1}, t\mapsto t^{-1}$. These 
polynomials are in terms of pairwise commutative variables $X_c$ 
$(c\in P)$; $q,t=\{t_{|\al|}\}$ are their parameters. Let 
$$\overline{E}_b\equal E _b ( X; q, t\to 0 ), \hskip 1mm 
E^{\dag\ast}_{b}=E_b^\ast(t\!\to\! 0),
$$ 
$h_b^0$ be the limit of the norm of $E _b ( X,q,t)$ as $t\to 0$, 
a very explicit product of certain $(1-q^j)$. Note that the 
polynomials $E^{\dag\ast}_{b}$ naturally appear in our 
expansions, not the dag-polynomials $E^{\dag}_{b}\equal 
E_b^\ast(t\!\to\! \infty)$. 

We have a standard theta-function $\theta$ for $R$: $\theta(X)\ 
\equal\ \sum_{b\in P} q^{(b,b)/2}X_b$, and its special 
normalization $\hat{\theta}\equal \theta/\lan 
\theta\mu_\circ(t\to 0) \ran$. See (\ref{mugau}) for the formula 
for the constant term $\lan\theta\mu_\circ\ran$ of 
$\theta\mu_\circ$ in the limit $t\to 0$; 
$\mu_\circ=\mu/\lan\mu\ran$ for the DAHA $\mu$\~function. Under 
this normalization, $\hat{\theta}$ can be identified with the 
graded character of a level-one integrable representation $L$ of 
the {\em twisted\,} affinization $\widehat{\mathfrak g}$ of a 
simple Lie algebra $\mathfrak g$ corresponding to the root 
system $R$. 

\subsubsection{\sf Rogers-Ramanujan sums}
The following particular case of  formula (\ref{pggmix})
(the ``mixed formula" from Theorem \ref{MAINTHM}) is 
the key. 

\begin{thm}\label{fmain}
For any $c\in P, p\in \N,$ and ${\mathbf b}=\{b_k\in P\mid 1\le 
k\le p\}$, 
$$
\overline{E}_{c^{\iota}}\,\hat{\theta}^p=
\sum_{\mathbf{b}}C_{\mathfrak{b}}\,
\frac{q^{\left((c^--b^-_1)^2+(b^-_1-b^-_2)^2+\ldots+
(b^-_{p-1}-b^-_p)^2\right)\!/\,2}}{\prod_{k=1}^{p-1}\,h^0_{b_k}}
\,\frac{E^{\dag\ast}_{b_p}}{h_{b_p}^0},
$$
where $C_{\mathbf{b}}$ is some power of $q$ depending on 
$\mathbf{b}$, whose definition requires the theory of 
$E^\dag$\~polynomials. 
\end{thm}

Setting $c = 0$, we arrive at an expansion of $\hat{\theta}^p$ 
in terms of $E^{\dag\ast}$\~polynomials divided by their norms. 
With a reservation about the range of $b_k$ (which is $P$, not 
$P_-$) and the $q$\~factors $C_{\mathbb{b}}$, this expansion is 
quite similar to that from \cite{ChFB}. Actually, it can be 
reduced to the formula there (in the absence of $c$) using some 
theory of $E^\dag$\~polynomials. 

To understand this series from representation-theoretic 
perspective, we use the associated graded pieces of {\em 
thick\,} Demazure filtrations of level-one integrable 
representations of $\widehat{\mathfrak g}$, called in this paper 
the {\em Demazure slices\,}. \vskip 0.2cm 

The Rogers-Ramanujan theory is of course not only about the 
summations; the product formulas (if they exist) are very 
important. We expect interesting representation-theoretic 
applications here. Let us also mention the identities connecting 
different expansion of the (same) products of theta-function. 
The latter are described in \cite{ChD} for arbitrary $t$ in the 
symmetric setting; the passage to the non-symmetric sums is 
straightforward. This is closely related to the {\em topological 
vertex\,} and  {\em DAHA-Jones polynomials\,} of Hopf links. 
Another perspective is a generalization of this paper to the 
{\em spinor $q$\~Whittaker global function\,} from \cite{ChO}. 
{\em Affine Hall functions\,} from \cite{ChSphI} must be 
mentioned too, especially Section 2.5 there (devoted to the 
Kac-Moody limit).

\subsubsection{\sf Demazure slices}
Let $\mathfrak g$ be a simple finite-dimensional
Lie algebra corresponding to $R$, 
$\widehat{\mathfrak g}$ its twisted affinization; $L$ a level-one 
integrable representation of $\widehat{\mathfrak g}$ (sometimes
called {\em basic\,}). 

\begin{thm}%[$=$ Corollary \ref{MODULEMAIN-U}] !!! RMV % LATER
\label{fDs} The $E^\dag$\~polynomials divided by their norms 
are precisely the graded characters of Demazure slices. 
\end{thm}

Theorem \ref{fDs} is (homologically) dual to the results of 
\cite{San,Ion1}, though the $E^\dag$\~polynomials are significantly 
more involved than the $\overline{E}$\~polynomials. It is also a 
``nonsymmetric" analogue of that in  
%Chari-Ion 
\cite{CI} and is closely related to 
%Kato-Loktev 
\cite{KL}. Using Theorem \ref{fDs}, we derive the following
theorem.

\begin{thm}%[$=$ Corollary \ref{EXT-UL}] !!! RMV % LATER
\label{fcmain} Let $p$ be a positive integer. For a level one 
thin Demazure module $D_b$ associated to $b \in P$, the module 
$D_b^{\vee} \otimes L^{\otimes p}$ admits a filtration by the 
Demazure slices (as constituents).
\end{thm}

This is essentially our representation-theoretic interpretation 
of Theorem \ref{fmain}. Combining Theorem \ref{fcmain} and 
Theorem \ref{fmain}, we obtain (quite non-trivial) formulas for  
the multiplicities of the Demazure slices in $D_b^{\vee} \otimes 
L^{\otimes p}$. Here obviously the ``numerical" nil-DAHA approach 
is very reasonable for obtaining explicit 
formulas vs. direct representation-theoretic calculations;
this was already demonstrated in \cite{ChFB}.

%\vskip 0.2cm 

Let us provide some details. We first prove the symmetric 
analogue of Theorem \ref{fDs}, which generalizes \cite{KL} (the 
$ADE$ case) to the twisted affinizations. Then we follow 
\cite{FKM} to obtain that Demazure slices have graded characters 
(homologically) dual to those of {\em thin\,} 
Demazure modules. Using this, Theorem \ref{fDs} follows from 
the identification of the proper $\mathrm{Ext}$\~pairing and the 
pairing from the theory of Macdonald polynomials.
Papers \cite{ChO,Kas05,Kat16} are used here. Finally, 
the machinery of Demazure-Joseph functors gives 
Theorem \ref{fcmain}. 

\vskip 0.2cm
The end of the paper contains  conjectures concerning the 
interpretation of the remaining cases of the ``numerical" 
theta-function expansions.
Let us mention here that Section \ref{sec:exp-theta} is more 
compressed and technically involved than the rest of the paper. We
provide sufficient references, but  it is
mostly aimed at specialists in affine Lie algebras. 

%\vskip 0.2cm
We note that one of the key application in \cite{ChFB} was that 
Rogers-Ramanujan-type expansions there almost automatically 
satisfy the {\em level-rank duality\,}, which is generally 
involved representation-theoretically. It will be interesting to 
see which kind of level-rank duality the nonsymmetric theory 
presented in this paper has.  \vfil\eject 

\setcounter{equation}{0} 
\section{\sc Affine roots systems and DAHA}
Let $R=\{\al\}   \subset \R^n$ be a root system of type 
$A,B,...,F,G$ with respect to a euclidean form $(z,z')$ on $\R^n 
\ni z,z'$,\ $W$ the {\em Weyl group} generated by the 
reflections $s_\al$,\ $R_{+}$ the set of positive  roots 
corresponding to fixed simple roots $\al_1,...,\al_n,$\ $\Ga$ 
the Dynkin diagram with $\{\al_i, 1 \le i \le n\}$ as the 
vertices,\ $\rho=\frac{1}{2}\sum_{\al\in R_+} \al$,\ 
$R^\vee=\{\al^\vee =2\al/(\al,\al)\}.$ 

The root lattice and the weight lattice are: 
\begin{align}
& Q=\oplus^n_{i=1}\Z \al_i \subset P=\oplus^n_{i=1}\Z \om_i,
\notag
\end{align}
where $\{\om_i\}$ are fundamental weights: $ 
(\om_i,\al_j^\vee)=\de_{ij}$ for the simple coroots 
$\al_i^\vee.$ Replacing $\Z$ by $\Z_{\pm}=\{k\in\Z, \pm k\ge 
0\}$ we obtain $Q_\pm, P_\pm.$ Here and further see  
\cite{Bo,C101}. 

The form will be normalized by the condition  $(\al,\al)=2$ for 
the {\em short\,} roots in this paper. Thus $\nu_\al\equal 
(\al,\al)/2$ can be either $1,$ $\{1,2\},$ or $\{1,3\}.$ We set 
$\nu_R \equal  \{\nu_\al \mid \al \in R \}$. \noindent 
%We set 
%$\nu_i\ =\ \nu_{\al_i}, \
%\nu_R\ = \{\nu_{\al}, \al\in R\}$.
%We will sometimes use the
%notation $\nu_{\lng}$ for the long roots ($\nu_{\sht}=1$).
This normalization leads to the inclusions $Q\subset Q^\vee,  
P\subset P^\vee,$ where $P^\vee$ is defined to be generated by 
the fundamental coweights $\{\om_i^\vee\}$ dual to $\{\al_i\}$. 

We note that $Q^\vee=P$ for $C_n \ (n\ge 2)$, $P\subset Q^\vee$ 
for $B_{2n}$ and $P\cap Q^\vee=Q$ for $B_{2n+1}$; the index 
$[Q^\vee:P]$ is $2^{n-1}$ for any $B_n$ (in the sense of 
lattices). 

\smallskip

\subsection{\bf Affine Weyl groups}\label{AWG}
The vectors $\ \tal=[\al,\nu_\al j] \in \R^n\times \R \subset 
\R^{n+1}$ for $\al \in R, j \in \Z $ form the {\em affine root 
system} $\tR \supset R$; this is the so-called twisted case. 
\smallskip

The vectors $z\in \R^n$ are identified with $ [z,0]$. We add 
$\al_0 \equal [-\vth,1]$ to the simple roots for the {\em 
maximal short root} $\vth\in R_+$. It is also the {\em maximal 
positive coroot} because of the choice of normalization. The 
(dual) Coxeter number is then $h=(\rho,\vth)+1$. The 
corresponding set $\tR_+$ of positive roots equals $R_+\cup 
\{[\al,\nu_\al j],\ \al\in R, \ j > 0\}$. 

We complete the Dynkin diagram $\Ga$ of $R$ by $\al_0$ (by 
$-\vth$, to be more exact); this is called the {\em affine 
Dynkin diagram} $\tGa$. 

%One can obtain it from the
%completed Dynkin diagram from \cite{Bo} for 
%the {\em dual system}
%$R^\vee$ by reversing all arrows. 
%The number of laces between $\al_i$ and
%$\al_j$ in $\tGa$ will be denoted by $m_{ij}.$

The set of the indices of the images of $\al_0$ by all the 
diagram automorphisms of $\tGa$ will be denoted by $O$;\, 
$O=\{0\} \for E_8,F_4,G_2$. Let $O'\equal \{r\in O, r\neq 0\}$. 
The elements $\om_r$ for $r\in O'$ are the so-called minuscule 
weights: $(\om_r,\al^\vee)\le 1$ for $\al \in R_+$ (here  
$(\om_r,\vth)\le 1$ is sufficient). 
\smallskip

\subsubsection{\sf Extended Weyl groups}
Given $\tal=[\al,\nu_\al j]\in \tR,  \ b \in P$, the 
corresponding reflection in $\R^{n+1}$ is defined by the formula 
\begin{align}
&s_{\tal}(\tz)\ =\  \tz-(z,\al^\vee)\tal,\
\ b'(\tz)\ =\ [z,\ze-(z,b)],
\label{ondon}
\end{align}
where $\tz=[z,\ze] \in \R^{n+1}$. 
\smallskip

The {\em affine Weyl group} $\tW$ is generated by all $s_{\tal}$ 
(we write $\tW = \lan s_{\tal}, \tal\in \tR_+\ran)$. One can 
take the simple reflections $s_i=s_{\al_i}\ (0 \le i \le n)$ as 
its generators and introduce the corresponding notion of the 
length. This group is the semidirect product $W\lsmash Q'$ of 
its subgroups $W=$ $\lan s_\al, \al \in R_+\ran$ and $Q'=\{a', 
a\in Q\}$, where 
\begin{align}
& \al'=\ s_{\al}s_{[\al,\,\nu_{\al}]}=\
s_{[\al,\,- \nu_\al]}s_{\al}\for
\al\in R.
\label{ondtwo}
\end{align}

The {\em extended Weyl group} $ \hW$ generated by $W\and P'$ 
(instead of $Q'$) is isomorphic to $W\lsmash P'$: 
\begin{align}
&(wb')([z,\ze])\ =\ [w(z),\ze-(z,b)] \for w\in W, b\in P.
\label{ondthr}
\end{align}
From now on,  $b$ and $b',$ $P$ and $P'$ will be identified. 
\smallskip

\comment{ For $1\le i\le n$, we set $\om_i=\pi_i u_i$, where  
$u_i\in W$ is a unique element of {\em minimal\,} possible 
length in the set $\{u\}\subset W$ such that $u(\om_i)\in P_-$, 
equivalently, $w=w_0u_i$ is of {\em maximal\,} length such that 
$w(\om_i)\in P_+$. The elements $u_i$ are very explicit. Let 
$w^i_0$ be the longest element in the subgroup $W_0^{i}\subset 
W$ of the elements preserving $\om_i$; this subgroup is 
generated by simple reflections. One has: 
\begin{align}\label{ururstar}
u_i = w_0w^i_0 \hbox{\,\, and\,\, } (u_i)^{-1}=w^i_0 w_0=
u_{i^\iota} \for 1\le i\le n,
\end{align}
\smallskip
where $\iota$ is the standard  automorphism of $\Ga$. } 

\subsubsection{\sf The length on \texorpdfstring{$\hat{W}$}{}}

Given $b\in P_+$, let $w^b_0$ be the longest element in the 
subgroup $W_0^{b}\subset W$ of the elements preserving $b$. This 
subgroup is generated by simple reflections. We set 
\begin{align}
&u_{b} = w_0w^b_0  \in  W,\ \pi_{b} =
b( u_{b})^{-1}
\ \in \ \hW, \  u_i= u_{\om_i},\pi_i=\pi_{\om_i},
\label{xwo}
\end{align}
where $w_0$ is the longest element in $W,$ $1\le i\le n.$

The elements $\pi_r\equal\pi_{\om_r}, r \in O'$ and 
$\pi_0=\hbox{id}$ leave $\tGa$ invariant and form a group 
denoted by $\Pi$, 
 which is isomorphic to $P/Q$ by the natural
projection $\{\om_r \mapsto \pi_r\}$. As to $\{ u_r\}$, they 
preserve the set $\{-\vth,\al_i, i>0\}$. The relations 
$\pi_r(\al_0)= \al_r= ( u_r)^{-1}(-\vth)$ distinguish the 
indices $r \in O'$. Moreover, 
\begin{align}\label{hWPi}
& \hW  = \Pi \lsmash \tW, \where
  \pi_rs_i\pi_r^{-1}\!=\! s_j \Leftrightarrow 
\pi_r(\al_i)\!=\!\al_j,\ \,
 0\le j\le n.
\end{align}

Setting $w = \pi_r u \in \hW,\ \pi_r\in \Pi, u \in \tW,$ {\em 
the length} $l(w)$ is by definition the length of the reduced 
decomposition $u= $ $s_{i_l}...s_{i_2} s_{i_1} $ in terms of the 
simple reflections $s_i, 0\le i\le n.$ The number of  $s_{i}$ in 
this decomposition such that $\nu_i=\nu$ is denoted by   
$l_\nu(w).$ 

The length can be also defined as the cardinality $|\la(w)|$ of 
the {\em $\la$\~sets\,}: 
\begin{align}\label{lasetdef}
&\la(w)\equal\tR_+\cap w^{-1}(\tR_-)=\{\tal\in \tR_+,\
w(\tal)\in \tR_-\},\
w\in \hW;\\
&\hbox{here\, }:\, \la(w)=\cup_\nu\la_\nu(w),\
\la_\nu(w)\
\equal\ \{\tal\in \la(w),\nu({\tal})=\nu \}.
\label{xlambda}
\end{align}
See, e.g. \cite{Bo,Hu} and also \cite{C101}. 

\subsubsection{\sf Elements
\texorpdfstring{\mathversion{bold}$\pi_b, u_b$}{\em pi and u}} 
Extending the definition of $\pi_b,u_b$ from $b\in P_+$ to any 
$b\in P$, we have the following proposition. 

\begin{proposition} \label{PIOM}
 Given $ b\in P$, there exists a unique decomposition
$b= \pi_b  u_b,$ $ u_b \in W$ satisfying one of the following 
equivalent conditions: 

{(i) \  } $l(\pi_b)+l( u_b)\ =\ l(b)$ and $l( u_b)$ is the 
greatest possible, 

{(ii)\  } $ \la(\pi_b)\cap R\ =\ \emptyset$. 

The latter condition implies that $l(\pi_b)+l(w)\ =\ l(\pi_b w)$ 
for any $w\in W.$ Besides, the relation $ u_b(b) \equal b_-\in 
P_-=-P_+$ holds, which, in its turn, determines $ u_b$ uniquely 
if one of the following equivalent conditions is imposed: 

{(iii) } $l( u_b)$ is the smallest possible, 

{(iv)\ } if\, $\al\in \la( u_b)$  then $(\al,b)\neq 0$. 
\end{proposition}
\qed 

We will need the following explicit description of the sets 
$\la(b)$. For $\tal=[\al,\nu_\al j]\in \tR_+,$ one has: 
\begin{align}
\la(b) = \{ \tal,\  &( b, \al^\vee )>j\ge 0 \iif \al\in R_+,
\label{xlambi}\\
&( b, \al^\vee )\ge j> 0 \iif \al\in R_-\},
\notag \\
\la(\pi_b) = \{ \tal,\ \al\in R_-,\
&( b_-, \al^\vee )>j> 0
\iif  u_b^{-1}(\al)\in R_+,
\label{lambpi} \\
&( b_-, \al^\vee )\ge j > 0 \iif
u_b^{-1}(\al)\in R_- \}, \notag \\
\la(\pi_b^{-1}) = \{ \tal\in \tR_+,\  -&(b,\al^\vee)>j\ge 0 \},
\label{lapimin}\\
\la(u_b)\ =\ \{ \al\in R_+,\ \ \ \,&(b,\al^\vee)> 0 \}.
\label{laumin}
\end{align}
For instance, $l(b)=l(b_-)=-2(\rho^\vee,b_-)$ for 
$2\rho^\vee=\sum_{\al>0}\al^\vee.$ 
\medskip

There is an important interpretation of the length and the 
elements $\pi_b, u_b$ in terms of the following {\em affine 
action} of $\hW$ on $z \in \R^n$: 
\begin{align}
& (wb)\llb z \rrb \ =\ w(b+z),\ w\in W, b\in P,\notag\\
& s_{\tal}\llb z\rrb\ =\ z - ((z,\al^\vee)+j)\al,
\ \tal=[\al,\nu_\al j]\in \tR.
\label{afaction}
\end{align}
For instance, $(b w)\llb 0\rrb=b$ for any $w\in W.$ The relation 
to the above action is given in terms of the {\em affine 
pairing\,} $([z,l], z'+\mathsf{d})\equal (z,z')+l:$ 
\begin{align}
& (w([z,l]),w\llb z' \rrb+\mathsf{d}) \ =\
([z,l], z'+\mathsf{d}) \for \hw\in \hW,
\label{dform}
\end{align}
where we treat $(\,,\,+\mathsf{d})$ formally (one can add 
$\mathsf{d}$ to $\R^{n+1}$ and extend $(\,,\,)$ correspondingly; 
compare with Sect. \ref{uDm}). 

\subsubsection{\sf Affine Weyl chamber}
Introducing the {\em basic affine Weyl chamber\,} 
\begin{align}
&\CC^a\ =\ \bigcap_{i=0}^n \LL_{\al_i},\
\LL_{[\al,\nu_\al j]}=\{z\in \R^n,\
(z,\al)+j>0 \},
\notag \end{align}
we come to another interpretation of the $\la$\~sets: 
\begin{align}
&\la_\nu(w)\ =\  \{\tal\in \tilde{R}_+, \,   \CC^a
\not\subset w\llb \LL_{\tal}\rrb, \, \nu_{\al}=\nu \}.
\label{lamaff}
\end{align}
Equivalently, taking a vector $\xi\in \CC^a$, 
\begin{align}\label{lageomet}
\la(w)\,=\,\{\tal\in \tR\,\mid\,
(\tal^\vee,\xi+\mathsf{d})>0>(\tal^\vee,\xi'+\mathsf{d})\}
\end{align}
for $\xi'\in \hw^{-1}\llb \CC^a\rrb$. \comment{ Thus, we come to 
the following geometric description of the $\la$\~sets. 
\begin{proposition}\label{GEOMLA}
The $\la$\~sets are exactly those in the form (\ref{lageomet}) 
for an arbitrary $\xi\in \CC^a$ and an arbitrary vector $\xi'$ 
inside a certain (affine) Weyl chamber, i.e., provided that 
$(\xi',\al^\vee)\not\in \nu_{\al}\Z$ for all $\al\in R$. 
%Given a generic segment in $\R^n$ 
%from $\xi$ to $\xi'$, its consecutive
%intersections with the affine root hyperplanes arrange such
%set into a {\color{red} $\la$\~sequence
%\footnote{Is this defined earlier?}}.
\sq 
\end{proposition}
} Geometrically, $\Pi$ is the group of all elements of $\hW$ 
preserving $\CC^a$ with respect to the affine action. Similarly, 
the elements $\pi_b^{-1}$ for $b\in P$ are exactly those sending 
$\CC^a$ to the basic nonaffine Weyl chamber $\CC\equal\{z\in 
\R^n,$$ (z,\al_i)>0$ for $i>0\}.$ 
\medskip

\subsection{\bf Partial ordering on 
\texorpdfstring{\mathversion{bold}$P$}{\em P}} 
\subsubsection{\sf The definition}\label{DefOrder}
It is necessary in the theory of nonsymmetric polynomials; see 
\cite{O2,M4}. This ordering was also used in \cite{C2} for 
calculating the coefficients of $Y$\~operators. The definition 
is as follows: 
\begin{align}
&b \ll c, \ c\gg b \for
b, c\in P \iif c-b \in Q_+ \hbox{\, and\,}
b\neq c
\label{order}
\\ &b \prec c,\  c\succ b \iif b_-\ll c_- \hbox{\ or\  }
\{b_-=c_- \hbox{\ and\ } b\ll c\}.
\label{succ}
\end{align}
Recall that $b_-=c_- $ means that $b,c$ belong to the same 
$W$\~orbit. We  write  $\preceq, \succeq$ respectively if $b$
can coincide with $c$. 

The following sets 
\begin{align}
&\si(b)\equal \{c\in P, c\succeq b\},\
\si_*(b)\equal \{c\in P, c\succ b\}, \notag\\
&\si_-(b)\equal \si(b_-),\
\si_+(b)\equal \si_*(b_+)= \{c\in P, c_- \gg b_-\}.
\label{cones}
\end{align}
are convex. By {\dfont convex}, we mean that if $ c, d= 
c+r\al\in \si$ for $\al\in R_+, r\in \Z_+$, then 
\begin{align}
&\{c,\ c+\al,...,c+(r-1)\al,\ d\}\subset \si.
\label{convex}
\end{align}

\subsubsection{\sf Bruhat ordering etc}
We will use the standard Bruhat ordering. Given $w\in \hW$, the 
{\em standard Bruhat set} $\b(w)$ is formed by $u$ obtained by 
striking out any number of $\{s_{j}\}$ from a reduced 
decomposition of $w\in \hW$. The
notation is $u\le w$.
The set $\b(w)$ does not depend on 
the choice of the reduced decompositions. 

\begin{proposition }\label{BSTAL}
(i) Let $c=u\llb 0\rrb$,\, $b=w\llb 0\rrb$ and $u\in \b(w)$. 
\comment{ The latter means that $u$ can be obtained by deleting 
simple reflections from the product of any reduced 
decompositions of $\pi_b$ and $v$ in the decomposition $w=\pi_b 
v$ with $v \in W$ so that $\ell ( w ) = \ell ( \pi_b ) + 
\ell ( v )$. } Then $c\succeq  b$ and  $b-c$ is a linear 
combination of 
the non-affine components of the corresponding roots from 
$\la(w^{-1})$. Also, $c=b$ if and only if $u$ is obtained by 
striking out $s_{j}$ from $v\in W$ such that $w=\pi_b v$ and 
this product is reduced, i.e. $\ell(w)=\ell(\pi_b)+\ell(v)$. 

(ii) Letting  $b=s_i\llb c\rrb$ for $0\le i\le n,$ if the 
element $s_i\pi_c$ belongs to $\{\pi_a,\,a\in P\}$ then it 
equals $\pi_b$. It happens if and only if 
$(\al_i,c+\mathsf{d})\neq 0.$ More precisely, the following 
three conditions are equivalent: 
\begin{align}
& \{c\succ b\}\Leftrightarrow  \{(\al_i,c+\mathsf{d})>0\}
\Leftrightarrow \{s_i\pi_c=\pi_b,\ l(\pi_b)=l(\pi_c)+1\}.
\label{aljb}
\end{align}
The latter relation implies that 
$\la(\pi_b)=\pi_c^{-1}(\al_i)\cap\la(\pi_c)$. \qed 
\end{proposition }
%\vfil
The following lemma is Lemma 1.7 from \cite{C102}; it  extends 
part $(ii)$ to the case $(\al_i,c\!+\!\mathsf{d})\!=\!0.$ 

\begin{lemma}\label{VANISHL}
The condition $(\al_i,c+\mathsf{d})=0$ for $0\le i\le n$, 
equivalently, the condition $(\al_i,b+\mathsf{d})=0$ for 
$b=s_i\llb c\rrb$, implies that $u_c(\al_i)=\al_j$ for $i>0$ or 
$u_c(-\vth)=\al_j$ for $i=0$ for a proper index $j>0.$ Given $c$ 
and $i\ge 0$, the existence of such $\al_j$ and the equality 
$(\al_j,c_-)=0$ are equivalent to $(\al_i,c+\mathsf{d})=0$.\sq 
\end{lemma}
\comment{ 
\begin{proposition}\label{BSTIL}
(i) Assuming (\ref{aljb}), let $i>0$. Then $b=s_i(c),\ $ 
$b_-=c_-$\, and $u_b=u_c s_i.$ The set $\la(\pi_b)$ is obtained 
from $\la(\pi_c)$ by adding $[\al,(c_-,\al)]$ for 
$\al=u_c(\al_i)$, i.e., by replacing the inequality $( c_-, 
\al^\vee )>j> 0$ in (\ref{lambpi}) with $( c_-, \al^\vee )\ge j> 
0.$ Here $\al\in R_-$ and $(c_-,\al^\vee)=(c,\al_i^\vee)>0.$ 

(ii) In the case $i=0,\ $ the following holds: 
$b=\vth+s_{\vth}(c)$, the element $c$ is from $\si_+(b)$, 
$b_-=c_- - u_c(\vth)\in P_-$ and $u_b=u_c s_{\vth}$. For $\al= 
u_c(-\vth)=\al^\vee,$ the $\la$-inequality $( c_-, \al)\ge j> 0$ 
is replaced with $( c_-, \al)+1 \ge j> 0$; respectively, the 
root $[\al,(c_-,\al)+1]$ is added to $\la(\pi_c)$. Here 
$\al=\al^\vee\in R_-$ and $( c_-, \al)=-(c,\vth)\ge 0.$ 

(iii) For any $c\in P, r\in O',\ $ one has $\pi_r\pi_c=\pi_b$ 
where $b=\pi_r\llb c\rrb.$ Correspondingly, $u_b=u_c u_r,\ b= 
\om_r+u_r^{-1}(c),\ b_-=c_- +u_c w_0(\om_{r}).$ In particular, 
the latter weight always belongs to $P_-.$ 
\end{proposition}
\smallskip
} 

\subsection{\bf Polynomial representation} \label{SEC:POL}
For the
variables $X_1,\ldots, X_n$, let 
\begin{align}
& X_{\tb}\ =\ \prod_{i=1}^nX_i^{l_i} q^{k}
\iif \tb=[b,k],\ w(X_{\tb})\ =\ X_{w(\tb)}.
\label{Xdex}\\
&\hbox{where\ } b=\sum_{i=1}^n l_i \om_i\in P,\ k \in
\Q,\ w\in \hW.
\notag \end{align}
For instance, $X_0\equal X_{\al_0}=qX_\vth^{-1}$. We will set 
$(\tilde{b},\tilde{c})=(b,c)$, i.e. we ignore the affine 
extensions in this pairing. 

Note that $\pi_r^{-1}$ is $\pi_{r^*}$ and $u_r^{-1}$ is 
$u_{r^*}$ for $r^*\in O\ ,$ where  the reflection $\{\cdot\}^*$ is 
induced by an involution $\iota$ of the nonaffine Dynkin diagram 
$\Gamma.$ 

By $e,$ we will denote  the least natural number such that  
$\ee(P,P)/2=\Z.$  Thus $\ee=4 \for D_{2k},\ 
\ee=2 \for B_{2k} \and 
C_{k},$ otherwise $\ee=2|\Pi|$. 

\subsubsection{\sf 
\texorpdfstring{\mathversion{normal}$Y$}{\em Y}-operators} 

For $\tal=[\al,\nu_\al j] \in \tR,\ \,0\le i\le n$,\, we set 
$$t_{\tal} =t_{\al}=t_{\nu_\al}, \ \,  t_{\al_i} \equal t_i.$$
The {\em Demazure-Lusztig operators\,} are as follows: 
\begin{align}
&T_i\  = \  t_i^{1/2} s_i\ +\
(t_i^{1/2}-t_i^{-1/2})(X_{\al_i}-1)^{-1}(s_i-1),
\ 0\le i\le n.
\label{Demazx}
\end{align}

They obviously preserve $\Z[q][X_b, b\in P]$. We note that only 
the formula for $T_0$ involves $q$: 
\begin{align}
&t_0^{1/2} s_0\ +\
(t_0^{1/2}-t_0^{-1/2})(X_{\al_0}-1)^{-1}(s_0-1),
\hbox{\ where\ }\notag\\
&X_0=qX_\vth^{-1},\
s_0(X_b)\ =\ X_bX_{\vth}^{-(b,\vth)}
q^{(b,\vth)},\
\al_0=[-\vth,1].
\end{align}
We will also need $\pi_r$ $(r\in O')$; they act via the general 
formula $w(X_b)=X_{w(b)}$ for $w \in \hW$.

Given $w \in \tW, r\in O,\ $ the product 
\begin{align}
&T_{\pi_r w}\equal \pi_r\prod_{k=1}^l T_{i_k},\where
w =\prod_{k=1}^l s_{i_k},
l=l(\tw),
\label{Twx}
\end{align}
does not depend on the choice of the reduced decomposition 
(because $T_i$ satisfy the same ``braid'' relations as $s_i$ 
do). Moreover, 
\begin{align}
&T_{v}T_{w}\ =\ T_{vw}\  \hbox{ whenever}\
 l(vw)=l(v)+l(w) \for
v,w \in \hW. \label{TTx}
\end{align}
In particular, we arrive at the pairwise commutative elements 
\begin{align}
& Y_{b}\equal
\prod_{i=1}^nY_i^{l_i} \iif
b=\sum_{i=1}^n l_i\om_i\in P,\ 
Y_i\equal T_{\om_i},b\in P.
\label{Ybx}
\end{align}

The action of $T_i$, $\pi_r$ $(r\in O)$ and $X_b$, considered to 
the operators of multiplication by $X_b$ (see (\ref{Xdex})), 
induces a representation of DAHA, the abstract algebra generated 
by these operators. It is called the {\em polynomial 
representation\,}; the notation is 
$$
\v\equal \ \Z[q^{\pm1/(2\ee)},t_\nu^{\pm1/2}][X_b, b\in P].
$$
It is a faithful DAHA-module if $q$ is not a root of unity, so 
we can skip the definition of DAHA in this paper. 

\subsubsection{\sf The
\texorpdfstring{\mathversion{bold}$\mu$}{\em mu}-function} The 
following exponential notation for $t_i$ in terms of parameters 
(or complex numbers) $k_\nu$ will be convenient in quite a few 
formulas. For $\tal=[\al,\nu_\al j] \in \tR,\ \,0\le i\le n$, 
$k_\al\equal k_{\nu_\al},$\, we set 
\begin{align}\label{talqal}
&  t_{\tal} =t_{\al}=t_{\nu_\al}=q_\al^{k_\nu} ,\ \, 
q_{\al}=q^{\nu_\al},\ \,
\rho_k\equal (1/2)\,\!\sum_{\al>0} k_\al \al.
%\notag
\end{align}

For instance, by $X_\al(q^{\rho_k})=q^{(\rho_k,\al)}$, we mean 
$\prod_{\nu\in\nu_R}t_\nu^{(\rho_k,\al) / \nu}$, where $\al\in 
R$. This product contains only {\em integral} powers of 
$t_{\sht}$ and $t_{\lng}$ ($t_\al$ for short and long roots). 
Note that $(\rho_k,\al_i^\vee)=k_i=k_{\al_i}$ for $i>0$. 

The {\em truncated theta-function\,}, which is the key in the 
definition of the inner product of the polynomial DAHA 
representation (and in the theory of nonsymmetric Macdonald 
polynomials) is as follows: 
\begin{align}
&\mu(X;q,t) =\prod_{\al \in R_+}
\prod_{j=0}^\infty \frac{(1-X_\al q_\al^{j})
(1-X_\al^{-1}q_\al^{j+1})
}{
(1-X_\al t_\al q_\al^{j})
(1-X_\al^{-1}t_\al^{}q_\al^{j+1})}.\
\label{mu}
\end{align}
We will consider $\mu$ as a Laurent series with coefficients in  
the ring $\Q[t_\nu][[q]]$. The constant term of a Laurent series 
$f(X)$ will be denoted by $\lan  f \ran$ through the paper. One 
has: 
\begin{align}
&\lan\mu\ran\ =\ \prod_{\al \in R_+}
\prod_{i=1}^{\infty} \frac{ (1- q^{(\rho_k,\al)+i\nu_\al})^2
}{
(1-t_\al q^{(\rho_k,\al)+i\nu_\al})
(1-t_\al^{-1}q^{(\rho_k,\al)+i\nu_\al})
}.
\label{consterm}
\end{align}
                     %%%CHANGED:
Using that $q^{(z,\al)}=q_\al^{(z,\al^\vee)}$, we can set here 
$q^{(\rho_k,\al)+i\nu_\al}=q_\al^{(\rho_k,\al^\vee)+i}$. This 
formula is equivalent to the Macdonald constant term conjecture. 

Let $\mu_\circ\equal \mu/\lan \mu \ran$. The coefficients of the 
Laurent series $\mu_\circ$ are from the field of rationals 
$\Q(q,t)\equal\Q(q_\nu,t_\nu), \where \nu\in \nu_R.$ We set 
\begin{align}
\label{innerpro}
&\lan f,g\ran\ \equal\ \langle
 f\ {g}^* \mu_\circ\rangle\ =\
\langle g,f\rangle^* \for
f,g \in \Q(q^{\frac{1}{2\ee}},t^{\frac{1}{2}}_\nu)[X],\\
&\hbox{where\, } 
X_b^* =  X_{-b}, (t_\nu^u)^*= t_\nu^{-u}, (q^u)^* = 
q^{-u}\for u\in \Q.\notag
\end{align}
Note that  $\mu_\circ^*\ =\ \mu_\circ$. Also $\langle 
A(f),g\rangle =\langle f, A^{-1}(g)\rangle$, for the DAHA 
generators $A=X_i \ (1\le i\le n), \, \pi_r \ (r\in O), \,T_i \ 
(i\ge 0)$ and therefore for any $A=T_{w}, Y_b$. So $\v$ is 
formally a $\ast$\~unitary representation of DAHA. 

We will mostly need this pairing in the limit $t_\nu\to 0$. For 
later reference, let $\overline{\mu}\equal\mu(t_\nu\!\to\!0, 
\nu\in \nu_R)$. Then 
\begin{align}
&\overline{\mu}\ = \ \prod_{\al \in R_+}
\prod_{j=0}^\infty (1-X_\al q_\al^{j})
(1-X_\al^{-1}q_\al^{j+1}),
\label{mubar}\\
&\lan\overline{\mu}\ran\ =\ \prod_{i=1}^{n}
\prod_{j=1}^{\infty} \frac{1}
{1-q_i^j}\,, \where q_i=q^{\nu_i},\ 
\nu_i=\nu_{\al_i}=\frac{(\al_i,\al_i)}{2}.
\label{constermbar}
\end{align}
\medskip

\subsection{\bf \texorpdfstring{\mathversion{bold}$E$}{\em E}
-polynomials}\label{defNSMac} There are two equivalent 
definitions of the {\em nonsymmetric Macdonald polynomials\,}, 
denoted by $E_b=E_b(X;q,t)$ for $b\in P$. They belong to 
$\Q(q,t)[X_a,a\in P]$ and, using the pairing $\lan\ ,\ \ran$, 
can be introduced by means of the conditions 
\begin{align}
&E_b\!-\!X_b\in \Si_+(b)\!\equal\!
\oplus_{c\succ b}\Q(q,t) X_c,\,
\langle E_b, X_{c}\rangle\!=\! 0 \hbox{\, for\,} 
P\!\ni\! c\succ b.
\label{macd}
\end{align}
They are well defined because the  pairing is nondegenerate (for 
generic $q,t$) and form a basis in $\Q(q,t)[P]$. 

This definition is due to Macdonald (for $k_{\sht}=k_{\lng}\in 
\Z_+ $), who extended Opdam's nonsymmetric polynomials 
introduced in the differential case in \cite{O2} (Opdam mentions 
Heckman's unpublished lectures in \cite{O2}). The general 
(reduced) case was considered in \cite{C4}. 

\subsubsection{\sf Using
\texorpdfstring{\mathversion{normal}$Y$}{\em Y}-operators} 
Another approach to $E$\~polynomials is based on the 
$Y$\~operators. We continue using the same notation $X,Y,T$ for 
these operators acting in the polynomial representation. 

\begin{proposition}
The polynomials $\{E_b, b\in P\}$ are unique (up to 
proportionality) eigenfunctions of the operators $\{L_f\equal 
f(Y_1,\ldots, Y_n)\},$ where $f\in \Q[X]$, acting in 
$\Q_{q,t}[X]:$ 
\begin{align}
&L_{f}(E_b)\ =\ f(q^{-b_\#})E_b\, \hbox{\ for\ }\,
b_\#\equal b- u_b^{-1}(\rho_k),
\label{Yone} \\
& X_a(q^{b})\ =\
q^{(a,b)}\ ,\hbox{where\ } a,b\in P,\
 u_b=\pi_b^{-1}b, \label{xaonb}
\end{align}
$u_b$ is from Proposition \ref{PIOM},\ 
$b_\#=\pi_b(\!(-\rho_k)\!)$. \label{YONE} 
\end{proposition}
\qed 

The coefficients of the Macdonald polynomials are rational 
functions in terms of $q,t_\nu$ (here either approach can be 
used). Note that $b_\#=b-\rho_k$ for $b\in P_-$ and 
$b_\#=b+\rho_k$ for generic $b\in P_+$ (such that $(b,\al_i)>0$ 
for $i=1,\ldots, n$). 
\smallskip
 
\subsubsection{\sf Standard identities}
We will need the following evaluation and norm formulas. One 
has: 

\begin{align}
&E_{b}(q^{-\rho_k}) \ =\ q^{(\rho_k,b_-)}
\prod_{[\al,j]\in \la\,'\,(\pi_b)}
\Bigl(
\frac{
1- q_\al^{j}t_\al X_\al(q^{\rho_k})
 }{
1- q_\al^{j}X_\al(q^{\rho_k})
}
\Bigr),
\label{ebebs}\\
&\la\,'\,(\pi_b)\ =\
\{[\al,j]\ |\  [-\al,\nu_\al j]\in \la(\pi_b)\}.
\label{jbseto}
\end{align}
Explicitly (see (\ref{lambpi})), 
\begin{align}
\la\,'\,(\pi_b)\ =& \{[\al,j]\, \mid\,\al\in R_+,
\label{jbset}
\\
&-( b_-, \al^\vee )>j> 0
\iif  u_b^{-1}(\al)\in R_-,\notag\\
&-( b_-, \al^\vee )\ge j > 0 \iif
u_b^{-1}(\al)\in R_+ \}.
\notag \end{align}
Formula (\ref{ebebs}) is the nonsymmetric version of the 
Macdonald {\em evaluation conjecture\,} from \cite{C4}. The {\em 
norm-formula\,} is as follows: 

\begin{align}\label{normepols}
&\lan E_b,E_c\ran \!=\!
\de_{bc}\!\!\!\prod_{[\al,j]\in \la'(\pi_b)}
\Bigl(
\frac{
(1\!-\!q_\al^j t_\al^{-1} X_\al(q^{\rho_k}))
(1\!-\!q_\al^j t_\al X_\al(q^{\rho_k}))
}{
(1\!-\!q_\al^j X_\al(q^{\rho_k})) 
(1\!-\!q_\al^j X_\al(q^{\rho_k}))
}
\Bigr).
\end{align}

For later reference, let 
\begin{align}\label{ghqt}
&g_b(q,t)\equal E_b(q^{-\rho_k}),\
h_b(q,t)\equal \lan E_b,E_b\ran \for b\in P. \\
&\hbox{Assuming that\, } t_\nu\to 0 \hbox{\, for all\, }
\nu\in \nu_R \,(\hbox{we set\,} t\to 0),\notag\\
\label{limeval}
&\lim_{t\to 0} q^{-(\rho_k,b_-)}g_b= 1,\ 
h_b^0\,\equal\,\lim_{t\to 0} h_b\,=\!
\prod_{[\al_i,j]\in \la'(\pi_b)}
(1\!-\!q_i^j).
\end{align}
i.e. the last product is over $[-\al,\nu_{\al}j]\in\la(\pi_b)$ 
with simple $\al$. Also: 
\begin{align}\label{ghgt0}
&\lan \overline{E}_b E_c^{\dag\ast}\mu_\circ(t\!\to\!0)
\ran\!=\!\de_{bc}h^0_b
\hbox{\, for\, }\overline{E}_b\!=\!E_b(t\!\to\! 0), 
E_c^{\dag\ast}\!=\!E_c^\ast(t\!\to\! 0).
\end{align} 
See formula (3.42) from \cite{ChO}. 

\medskip
\setcounter{equation}{0} 
\section{\sc Theta-products via 
\texorpdfstring{\mathversion{bold}$E$}{\em E} -polynomials} 
\subsection{\bf Mehta-Macdonald identities}
\subsubsection{\sf Basic notations}
Our approach is based on the difference Mehta-Macdonald formulas 
for {\em standard theta-functions\,}. Given a root system $R$ as 
above and the corresponding $P,Q$, they generally depend on the 
choice of a character $v: \Pi=P/Q\to \C^*$. The group of such 
characters will be denoted by $\Pi'$; the trivial character will 
be $1'$. Let 
\begin{align}\label{zetauv|}
\ze_v(X_a)=v(a)X_a T_w Y_b \for
a\in P.
\end{align} 
For a character $v\in \Pi'$, we set 
\begin{align}
&\theta_v(X)\ \equal\ 
\sum_{b\in P} v(b) q^{(b,b)/2}X_b=
\ze_v(\theta), \where \theta\equal \theta_{1'}.
\label{gauser}
\end{align}

The characters $v$ play here the role of the classical {\em 
theta- characteristics\,}. Definition (\ref{gauser}) is directly 
related to that from \cite{ChFB}, though we used a somewhat 
different setting there. Namely, theta-functions were introduced 
using the partial summations in the series for $\th$, where the 
images of $b$ were taken from some subsets  $\varpi\subset\Pi$. 
Using $v$ instead of $\varpi$, we obtain two different basis in 
the same space (of theta-functions). The usage of $\varpi$ has 
some advantages for the $PSL(2,\Z)$\~modularity and when the 
{\em string functions\,} are considered \cite{KP}.
One of the most 
important observations in \cite{ChFB} was that the 
Rogers-Ramanujan formulas give almost immediate justification of 
the {\em level-rank duality\,} for related string functions. See 
also Lemma 3.1 from \cite{ChD} concerning the ``topological 
DAHA-vertex". 

\medskip

We will use Theorem 5.1 from \cite{C5} (the first formula) and 
Theorem 3.4.5 from \cite{C101}. Enhancing them by the characters 
$v\in \Pi'$ follows Theorem 3.2 from \cite{ChD}, where 
(\ref{epep}) below was proven; the justification of 
(\ref{epepst}) is quite similar. 

We will need {\em nonsymmetric spherical polynomials\,} 
 $\e_b\equal E_b/E_b(q^{-\rho_k})$ for $b\in P$.
See \cite{C101}, and formula (6.30) from \cite{C102}. Using 
these polynomials and the definition of $b_\#$ from 
(\ref{Yone}), the {\em duality theorem\,} states that 
\begin{align}\label{ebdual}
\e_b(q^{c_{\#}})=\e_c(q^{b_{\#}}),\hbox{\,where\,}  
b_\# = b-u_b^{-1}(\rho_k),\, b,c\in P.
\end{align}

We will also use that $\mu^*_{\circ} = \mu_{\circ}$, where 
$\mu_\circ$ is understood as a Laurent series $1+\ldots\,$ in 
terms of $X_b$ with {\em rational\,} $q,t$\~coefficients. This 
readily give the relations: 
\begin{align}\label{ebdual2}
\lan \e_b, \e_c \ran = \lan \e_b^*, \e_c^* \ran =
\de_{bc}\lan \e_b, \e_b \ran=
\de_{bc}\lan \e_b, \e_b \ran^\ast
\for b, c \in P.
\end{align}
One can see the $\ast$\~invariance of $\lan \e_b, \e_b \ran$ 
directly from (\ref{normepols}). 
\smallskip

\subsubsection{\sf The key expansions}
Recall that $X_b^*=X_b^{-1}=X_{-b},$ for $b\in P$, $ q^*=q^{-1}, 
t^*=t^{-1}$  and $u_b(b)=b_-\in P_-$. The following formulas are 
the key for us. 

\begin{theorem}\label{EPEP}
For $b,c\in P$, 
\begin{align}
v(b+c)\lan \e_b \e_c \theta_v\mu_\circ \ran & \,=\,
q^{b_-^2/2+c_-^2/2 -(b_-+c_-,\rho_k)}\,
\e_c(q^{b_\#})\,\lan \theta\mu_\circ  \ran,
\label{epep}
\\
v(b-c)\langle \e_b \e_c^* \theta_v\mu_\circ  
\rangle & \,=\,
q^{b_-^2/2+c_-^2/2 -(b_-+c_-,\rho_k)}\,
%q^{(b_\#,b_\#)/2+(c_\#,c_\#)/2 -(\rho_k,\rho_k)}
\e^*_c(q^{b_\#})\,\lan \theta\mu_\circ  \ran.
\label{epepst}
\end{align}
Here the coefficients of the Laurent series for $\mu_\circ$ are 
naturally expanded in terms of positive powers of $\,q\,$ and 
the proportionality factor is a $q$\~generalization of the Mehta 
-Macdonald integral: 
\begin{align}
&\lan \theta\mu_\circ  \ran\ =\
\prod_{\al\in R_+}\prod_{ j=1}^{\infty}\Bigl(\frac{
1-t_\al^{-1} q_\al^{(\rho_k,\al^\vee)+j}}{
1-      q_\al^{(\rho_k,\al^\vee)+j} }\Bigr).
\label{mehtamu} 
\end{align}
\vskip -0.7cm \sq 
\end{theorem}

For later reference, let us make $t\to 0$ in (\ref{mehtamu}):  
\begin{align}\label{mugau}
\lan \theta\overline{\mu}_\circ  \ran\ =\ 
\lan \overline{\mu}_\circ \ran^{-1}\ =\ 
\!\prod_{i=1}^{n}\prod_{j=1}^{\infty}
(1-q_i^{j}) \hbox{\, when \, } t\to 0.
\end{align}

Switching here (and below) to 
\begin{align}\label{veetheta}
\hat{\theta}_v\equal
\theta_v/\lan \theta\mu_\circ \ran,
\end{align}
the formulas above can be interpreted as the following 
expansions: 
\begin{align}\label{thexpan0}
&\e_c\hat{\theta}_v=
\!\sum_{b\in P}\frac{\lan \e_b \e_c \hat{\theta}_v\mu_\circ \ran}
{\lan \e_b,\e_b\ran}\e^*_b=\!
\sum_{b\in P}\frac{q^{\frac{b_{\!-}^2}{2}+
\frac{c_{\!-}^2}{2} -(b_{\!-} \!+ c_{\!-},\rho_k)}}
{v(b+c)\lan \e_b,\e_b\ran}\e_b(q^{c_\#})\e^*_b,\\
\label{thexpan}
&\e^*_c\hat{\theta}_v=
\sum_{b\in P}\frac{\lan \e_b \e^*_c \hat{\theta}_v\mu_\circ \ran}
{\lan \e_b,\e_b\ran}\e^*_b=
\!\!\sum_{b\in P}\frac{q^{\frac{b_{\!-}^2}{2}+
\frac{c_{\!-}^2}{2} -(b_{\!-} \!+ c_{\!-},\rho_k)}}
{v(b-c)\lan \e_b,\e_b\ran}\e^*_b(q^{c_\#})\e^*_b,\\
\label{thexpan1}
&\e_c\hat{\theta}_v=
\sum_{b\in P}\frac{\lan \e^*_b \e_c \hat{\theta}_v\mu_\circ \ran}
{\lan \e_b,\e_b\ran}\e_b=
\sum_{b\in P}\frac{q^{\frac{b_{\!-}^2}{2}+
\frac{c_{\!-}^2}{2} -(b_{\!-} \!+ c_{\!-},\rho_k)}}
{v(c-b)\lan \e_b,\e_b\ran}\e^*_b(q^{c_\#})\e_b.
\end{align}
Here (\ref{thexpan}) and (\ref{thexpan1}) follow from 
 (\ref{epepst}), where we use the duality. Formula 
(\ref{thexpan1}) in the symmetric variant was the starting point 
for paper \cite{ChFB}. 

Proposition 3.6 and Theorem 3.7 from \cite{ChD} contain a formal 
theory of iterations of these relation. They were stated there 
for symmetric Macdonald polynomials. Change the summations from 
$P_+$ to $P$, replace $-c-\rho_k$ for $c\in P_+$ by $c_\#$, and 
use $\frac{b_{-}^2}{2}+ \frac{c_{-}^2}{2} -(b_{-} \!+ 
c_{-},\rho_k)$ instead of $\frac{b^2}{2}+ \frac{c^2}{2} +(b+ 
c,\rho_k)$ to transfer the iteration formulas in \cite{ChD} to 
the $\e$\~polynomials. 

The formulas for iterations with generic $t$ are quite involved 
because there are no explicit formulas for 
 $\e^*_b(q^{c_\#})$. However
these quantities become some relatively simple $q$\~monomials in 
the limit $t\to 0$ and also when $t=1$. Let us discuss a little 
the case $t=1$, called ``free theory". One has:  $\mu_\circ=1, 
\lan \mu \ran=1$, $\e_b=E_b=X_b$ for any $b\in P$, and 
$\e_b(q^{c_\#})=q^{(b,c)}$. The expansion above reads 
$X_c\theta_v= \sum_{b\in P} q^{(b_- -c_-)^2/2}\, v(b-c)\,X_b$, 
which is simply the result of substitution $b\mapsto b-c$ in the 
formula for $\theta_v$. In spite of its simplicity, this formula 
plays the key role in the theory of theta-functions and 
Gaussians.

\medskip
\subsection{\bf Using
\texorpdfstring{\mathversion{bold}$E$}{\em E} -dag polynomials} 
The two limiting cases $t\!\to\! 0$ and $t\!\to\! \infty$ of the 
DAHA theory  are of significant importance: 
\begin{align}\label{E-dag}
\overline{E}_b(X;q)\!=\!E_b(X;q,t\!\to\!0),\, 
E^\dag_b(X;q)\!=\!E_b(X;q,t\!\to\!\infty),\, b\in P.
\end{align}
These polynomials are well defined; see \cite{ChO}. We will 
simply call them {\em E-dag\,} and {\em E-bar\,} polynomials. 
The former are nonsymmetric generalizations of the $q$\~Hermite 
polynomials, which coincide with the level-one Demazure 
characters in the twisted setting \cite{San,Ion1}. The latter 
are more recent; they were studied in \cite{ChO}. The 
coefficients of $E^\dag_b$, which are from $\Z[q^{-1}]$, were 
conjectured there to be in $\Z_+[q^{-1}]$. Furthermore, it was 
conjectured that 
\begin{align}\label{ncbcoef}
E_c^\dag\!= 
\!\!\sum_{b\in W(c)} %\ni c \succq b} 
q^{n_c(b)} X_b\!\!\mod\! \Si_+(c), 
\!\where n_c(b)\in\Z_{-},\ c\in P,
\end{align}
where we take only $b$ in the summation such that $X_b$ is 
present in $E_c^\dag$. See (\ref{macd}) for the definition of 
$\Si_+(c)$. 

The first conjecture was (generalized and) verified in 
\cite{OS}. The presentation (\ref{ncbcoef}) was checked in 
\cite{ChFE} for $c=c_-$ and verified in \cite{NNS} in the 
simply-laced case. It was mentioned in \cite{ChFE} that it 
formally follows from the case $c=c_-$; we will provide the 
proof below. 
\medskip

\subsubsection{\sf The numbers
\texorpdfstring{\mathversion{bold}$n_c(b)$}{\em n(b)}} By 
constriction, $n_c(c)=0$ and $q^{n_c(c_-)}=0$ unless $c=c_-$. 
Let us use formula (3.52) from \cite{ChO}. 

We assume in the next calculation that 
$(c,\al_i)=(c_-,u_c(\al_i))<0$ for $1\leq i\leq n$. Then 
$u_c(\al_i)\in R_+$. Also, $\pi_c=s_i\pi_{s_i(c)}$ is reduced 
($l(\pi_{c})=l(s_i)+l(\pi_{s_i(c)})$) and, equivalently, the 
product $u_{s_i(c)}=u_{c}s_i$ is reduced. One has: 
\begin{align}
\overline{E}_{s_i(c)}^\dag\ &=\
\begin{cases}
(1-q^{(c,\al_i)})^{-1}
(T_i^\dag)'(\overline{E}_{c}^\dag) 
\text{\,\, if\, \ } u_c(\al_i) &\hbox{ is simple},\\
(T_i^\dag)'(\overline{E}_{c}^\dag) 
&\hbox{ otherwise},\\
\end{cases}
\label{daginteb}\\
(T_i^\dag)'\ &=\ T_i^\dag-1=
\frac{X_{\al_i}}{X_{\al_i}\!-\!1}\,(s_i\!-\!1),\ 
T_i^\dag\equal t_i^{-1/2}T_i(t\to\infty).
\label{tdagprime}
\end{align}
Recall that $q^{(b,\al_i)}=q_i^{(b,\al_i^\vee)}$, 
$q_i=q^{\nu_i}=q^{(\al_i,\al_i)/2}$. For any $b\in P$, 
\begin{align}\label{tdagmomom}
&(T_i^\dag)'(X_b)\!\!\mod \Si_+(b) =
\begin{cases}
X_{s_i(b)},&\text{\, if\, } (b,\al_i)<0\\
\ -X_b, &\text{\, if\, } (b,\al_i)>0\\
\ \ \ \ \ 0,&\text{\, if\, } (b,\al_i)=0.
\end{cases}
\end{align}

Therefore: 
\begin{align}\label{ncbrel}
&(T_i^\dag)'(E^\dag_c) \!\!\mod \Si_+(b)=
\sum_{b\in W(c)}^{(b,\al_i)>0} (q^{n_c(s_i(b))}-q^{n_c(b)})X_b,\\
&q^{n_{s_i(c)}(b)}=
\begin{cases} 
0, &\text{\ \,unless\,\, } 
(b,\al_i)> 0,\\
\frac{q^{n_c(s_i(b))}-q^{n_c(b)}}
{1-q^{(c,\al_i)}}, &\text{\ \,if\, } 
u_c(\al_i) \hbox{ is simple},\\
q^{n_c(s_i(b))}-q^{n_c(b)}, &\text{\ otherwise},
\end{cases}\end{align}
where we impose $(b,\al_i)>0$ in the latter two equalities; 
$q^{n_c(b)}$ may be $0$ for such $b$. The role of simplicity 
$u_c(\al_i)$ is as follows: 
\begin{align}\label{ncbrels}
&\hbox{provided } (b,\al_i)\!>\! 0,\,\,
u_c(\al_i) \hbox{ is simple }\!\Longleftrightarrow
h_c^0\!\neq\! h_{s_i(c)}^0.
\end{align}

Employing the {\em monomiality claim\,}: 
\begin{align}\label{ncsi-rel}
&q^{n_{s_i(c)}(b)}=
\begin{cases}
0, &\text{\, if\, } (b,\al_i)\!\le\! 0
\hbox{\,\,\, or\,\,\,  } q^{n_c(s_i(b))}\!=\!q^{n_c(b)},\\
q^{n_c(s_i(b))},
&\text{\, if\, } (b,\al_i)\!>\! 0
\hbox{\, and\, } q^{n_c(s_i(b))}\!\neq\!q^{n_c(b)}. 
\end{cases}
\end{align}
Moreover, in the latter case,  $ 
q^{n_c(b)}\!=\!q^{n_c(s_i(b))+(c,\al_i)}$ for simple 
$u_c(\al_i)$ and $q^{n_c(b)}\!=\!0$ if $u_c(\al_i)$ is not 
simple; this root is positive due to $(c,\al_i)\!<\!0$. 

Using relations (\ref{ncsi-rel}), we can obtain all $q^{n_c(b)}$ 
for any $c\in W(c_-)\ni b$ assuming that they are known for 
$c=c_-$; see \cite{ChFE} for the latter case. Indeed, any 
element from the orbit $W(c)$ can be obtained from $c=c_-$ by 
consecutively applying proper $s_i$ under the negativity 
condition above: $c'=s_i(c), c''=s_{i'}c'$ for 
$(c',\al_{i'})<0$, and so on. 
\smallskip

\subsubsection{\sf An example for
\texorpdfstring{$A_2$}{\em A2}} It makes some sense to provide a 
simple example (the simplest beyond $A_1$). For $A_2$, let 
 $c_-=-\rho=-\om_1-\om_2$, $s_i=s_1, s_{i'}=s_2$. 
Thus $c=-\om_1-\om_2$, $c'=s_1(c)=\om_1-2\om_2, 
c''=s_2(c')=2\om_2-\om_1.$ Recall that $X_i=X_{\om_1}$. The 
first polynomial we provide below is 
 $\,E_c^\dag\,$. It is  followed by $\rightsquigarrow$ and
the  corresponding $\,\sum_{b\in W(c_-)} q^{n_{s_i(c)}(b)}X_b\,$ 
calculated via formula  (\ref{ncsi-rel}) for $c,s_i$. Note that 
it does coincide with the contribution of $X_b$ for $b\in W(c')$ 
in $E_{c'}^\dag$, which is provided next. Similarly,  $\sum_b 
q^{n_{s_{i'}(c')}(b)}X_b$ from  (\ref{ncsi-rel}) for $c',s_{i'}$ 
is given after $\rightsquigarrow$, which coincides with the 
corresponding portion of $E_{c''}^\dag$. The latter is the last 
polynomial we provide. 

{\small 
\begin{align*}
&\frac{1}{X_1 X_2}+\frac{X_2}{q X_1^2}+\frac{1}{q}
+\frac{2}{q^2}+\frac{X_1}{q X_2^2}
+\frac{X_1^2}{q^2 X_2}
+\frac{X_1 X_2}{q^2}+\frac{X_2^2}{q^2 X_1}
\, \rightsquigarrow\\
&\frac{X_1}{X_2^2}\!+\!\frac{X_1^2}{q X_2},\ \, \,
\frac{X_1}{X_2^2}\!+\!\frac{X_1^2}{q X_2}\!+\!\frac{1}{q}
\,\rightsquigarrow\,
\frac{X_2^2}{X_1}\!+\!\frac{X_1 X_2}{q},\ \,\, 
\frac{X_2^2}{X_1}\!+\!\frac{X_1 X_2}{q}\!+\!1\,.
\end{align*}
} 
\smallskip

\subsubsection{\sf Proof of monomiality}
The monomiality claim (\ref{ncbcoef}) was convenient to use when 
checking (\ref{ncsi-rel}). This is not necessary; the 
monomiality can be actually deduced from the same argument and 
the case $c=c_-, c'=s_i(c)$. 

\begin{proposition}\label{PROPNCB}
Let $\zeta_c(b)$ be the coefficient of $X_b$ for $b\in W(c)$ in 
$E_c^\dag$; we set $\zeta_-(b)=\zeta_{c_-}(b)$. Then given 
$c_-$, $\zeta_c(b)$ is either $0$ or coincides with 
$\zeta_-(u_c(b))$. More exactly, if $c'=s_i(c)$ and 
$(\al_i,c)<0$, then relations (\ref{ncsi-rel}) hold and one can 
proceed by induction. In particular, we arrive at the 
monomiality statement from (\ref{ncbcoef}). 
\end{proposition}

{\it Proof.} We argue by induction with respect to $l=l(u_c)$ 
following Proposition 7.4 from \cite{ChFE} and the next 
mini-section there ``On embeddings of dag-polynomials". They 
provide the expansion (\ref{ncsi-rel}) for $c=c_-$ and the 
monomiality for $c'=s_i(c_-)$, i.e. the induction step $l=1$. We 
will also use Corollary 3.6 from \cite{ChO}, which establishes 
that $\zeta_c(b)$ belong to $\Z[q^{-1}]$. To simplify the 
reasoning we use the positivity of $\ze_c(b)$ from \cite{OS} 
(which can be actually avoided here). \vskip 0.3cm 

For any $c$, let $c'=s_i(c)$. In the case of simple 
$u_c(\al_i)$, we have $u_c(\al_i)=\al_j$ for $1\le j\le n$. Then 
either both terms in $\ze_c(s_i(b))-\ze_c(b)$ vanish or neither 
of them. Indeed, if only one is zero, then it would contradict 
to $\ze_{c'}(b)\in \Z[q^{-1}]$ for $c'=s_i(c)$. If both vanish 
then $\ze_{c'}(b)=0$. If both are nonzero, the we obtain by the 
induction claim for $c$: 
\begin{align*}
&\ze_{c'}(b)=\frac{\ze_c(s_i(b))-\ze_c(b)}{1-q^{(c,\al_i)}}
=\frac{\ze_{-}(u_cs_i(b))-\ze_{-}(u_c(b))}
{1-q^{(c,\al_i)}}\\
=&\frac{\ze_{-}(s_ju_c(b))-\ze_{-}(u_c(b))}
{1-q^{(c_-,\al_j)}}, \hbox{\ where we use \ } u_cs_i=s_ju_c.
\end{align*}
Here the product $u_c s_i$ is reduced by construction and 
$(u_c(b),\al_j)=(b,\al_i)>0$; so we can use step  $l=1$ (for 
$c'=s_j(c_-)$).  
 
If $u_c(\al_i)$ is not simple, then we need to check that 
$\ze_-(u_c(b))$ from $\ze_{c'}(b)= 
\ze_c(s_i(b))-\ze_c(b)=\ze_-(u_cs_i(b))-\ze_-(u_c(b))\,$ is $0$ 
unless these two terms coincide. Note that there is no 
denominator now. This follows from some (minor) development of 
the method from \cite{ChFE} for $c=c_-$. Alternatively, we can 
simply use for this step that $\ze_{c'}(b)\in \Z_+[q^{-1}]$ from 
\cite{OS}; the presence of nonzero $-\ze_c(b)=-\ze_-(u_c(b))$ 
readily contradicts the positivity proven there. This gives the 
required.\sq 

\comment{ 
\begin{align}\label{ncbrelx}
&(T_i^\dag)'(E^\dag_c) \!\!\mod \Si_+(b)=
\sum_{b\in W(c)}^{(b,\al_i)>0} (q^{n_c(s_i(b))}-q^{n_c(b)})X_b,\\
&q^{n_{s_i(c)}(b)}=
\begin{cases} 
0, &\text{\ \,unless\,\, } 
(b,\al_i)> 0,\\
\frac{q^{n_c(s_i(b))}-q^{n_c(b)}}
{1-q^{(b,\al_i)}}, &\text{\ \,if\, } 
u_b(\al_i) \hbox{ is simple},\\
q^{n_c(s_i(b))}-q^{n_c(b)}, &\text{\ otherwise},
\end{cases}\end{align}
gives that $\zeta_c(s_i(b))=0$ unless $(b,\al_i)<0$. If 
additionally $u_c(\al_i)$ is not a simple root, then either 
$\zeta_c(s_i(b))=\zeta_c(b)$ or $\zeta_c(b)=0$ (this was not 
checked). 

It can be formally deduced from the case $c=c_-$. The second 
case in (\ref{ncbrel}) can be managed by a reduction to $q=1$ 
similar to the argument in \cite{ChFE}. However this reduction 
cannot be used for the last case there. Thus we need to verify 
the following lemma, which will follow from Theorem 
\ref{RR-THETA} below. 

\begin{lemma}\label{LEMNCB}
Let $\zeta_c(b)$ be the coefficient of $X_b$ in $E_c^\dag$.  
Assuming that $(c,\al_i)<0$, $\zeta_c(s_i(b))=0$ unless 
$(b,\al_i)<0$ (which was checked). If additionally $u_c(\al_i)$ 
is not a simple root, then either $\zeta_c(s_i(b))=\zeta_c(b)$ 
or $\zeta_c(b)=0$ (this was not checked).\sq 
\end{lemma}
} 
\smallskip

\subsection{\bf Three major expansions}
\subsubsection{\sf Non-spherical formulas} 
In the theory of $E^\dag$\~polynomials it is convenient to use 
the following (obvious) identity: 
\begin{align} \label{Edag0}
&E_b^\dag(X;q)= E^*_b(X^{-1};q^{-1},t\to 0), \ b\in P;
\hbox{\, \ see (\ref{innerpro})}.
\end{align}

Let us consider the limit $t\to 0$ in formulas 
(\ref{thexpan0},\ref{thexpan},\ref{thexpan1}). First, we need to 
restate them in terms of $E$\~polynomials using $g_b(q,t)$ from 
(\ref{ghqt}) and employing the duality in the middle formula: 

\begin{align}\label{thexpann0}
&\frac{E_c}{g_c}\hat{\theta}_v=
\sum_{b\in P}\frac{q^{\frac{b_{\!-}^2}{2}+
\frac{c_{\!-}^2}{2} -(b_{\!-} \!+ c_{\!-},\rho_k)}}
{v(b+c)\lan E_b,E_b\ran}\frac{g_b g_b^*}{g_b g_b^*}
E_b(q^{c_\#})\,E^*_b,\\
\label{thexpann}
&\frac{E^*_c}{g_c^*}\hat{\theta}_v=
\sum_{b\in P}\frac{q^{\frac{b_{\!-}^2}{2}+
\frac{c_{\!-}^2}{2} -(b_{\!-} \!+ c_{\!-},\rho_k)}}
{v(b-c)\lan E_b,E_b\ran}\frac{g_b g_b^*}{g_c^* g_b^*}
E^*_c(q^{b_\#})\,E^*_b, \\
\label{thexpann1}
&\frac{E_c}{g_c}\hat{\theta}_v=
\sum_{b\in P}\frac{q^{\frac{b_{\!-}^2}{2}+
\frac{c_{\!-}^2}{2} -(b_{\!-} \!+ c_{\!-},\rho_k)}}
{v(c-b)\lan E_b,E_b\ran}\frac{g_b g_b^*}{g_b^* g_b}
E^*_b(q^{c_\#})\,E_b.
\end{align}

\subsubsection{\sf The limit at zero}
Now let $t\to 0$ (i.e. $t_\nu\to 0$ for all $\nu$). We will use 
(\ref{limeval}) and $h_b^0=\lim_{t\to 0}h_b$. The last formula 
is somewhat simpler to analyze. Note the cancelation of 
$g_b$\~factors; also,  $1/g_c$ on the l.h.s. cancels 
$q^{-(c_-,\rho_k)}$ on the r.h.s. One obtains: 
\begin{align}
\label{thexpannn1}
&\overline{E}_c\hat{\theta}_v=
\sum_{b\in P}\frac{q^{\frac{b_{-}^2}{2}+
\frac{c_{-}^2}{2}}}
{v(c\!-\!b)h^0_b}
\,\bigl(\lim_{t\to 0} q^{-(b_{-},\rho_k)}
E^*_b(q^{c_\#})\bigr)\,\overline{E}_b.
\end{align}

%%%CHANGED:
Using that $c_\#=u_c^{-1}(c_- -\rho_k)$, the monomials $X_{-a}$ 
from $E_b^*$ that do not vanish in the limit $t\to 0$ upon the 
evaluation at $q^{c_\#}$ are for $a\in W(b)$ and must satisfy 
$$(b_-,\rho_k)=(a,u_c^{-1}(\rho_k))=
(u_c(a),\rho_k)$$ for generic $k$ (before the limit).
This gives the relation $a=u_c^{-1}(b_-)$. 

Let us comment on our usage of $q$ and $k$ here. We assume that 
$0<q_{\al}<1$ and that $k_{\nu}\to +\infty$ (they are generic 
real numbers); then $t_\nu\to 0$, which is exactly what we need 
here.  Thus $X_{-a}$ from $E_b^*$ contributes to the limit 
(\ref{thexpannn1}) only if 
$$\lim_{k_\nu\to \infty} q^{-(b_{-},\rho_k)}
X_{-a} (q^{c_\#}) \neq 0,
$$
and we arrive at the relation above.

\comment{ S.K. EDITED: Here the limit is taken under the 
assumption that all the $k_{\nu}$'s are generic real positive 
numbers (tends to $\infty$ by $t_{\al} \to 0$ with $q_{\al}$ 
fixed but formally understood as $|q_{\al}| < 1$). By the fact 
that $E_b^*(X;q,t \to 0)$ exists and $\lim_{t\to 0} 
q^{-(b_{-},\rho_k)} = 0$, the monomials $X_{-a}$ from $E_b^*$ 
contributes to the limit part of (\ref{thexpannn1}) only if 
$$\lim_{t\to 0} q^{-(b_{-},\rho_k)}
X_{-a} (q^{c_\#}) \neq 0.$$

Using that $c_\#=u_c^{-1}(c_- -\rho_k)$, this implies that 
$$(b_-,\rho_k) \ge (a,u_c^{-1}(\rho_k))=
(u_c(a),\rho_k).$$
On the other hand, we have 
$$(b_-,\rho_k) \le (a,u_c^{-1}(\rho_k))=
(u_c(a),\rho_k)$$
for every $a \in \si_- ( b )$, and the equality is attained only 
when $a=u_c^{-1}(b_-)$ (though it might happen that 
$u_c^{-1}(b_-) \in \si ( b )$). } 

Now let us use the expansion 
\begin{align}\label{ebstar}
E_b^*(X; q,t\!\to\!0)\!=\!E^\dag(X^{-1};q^{-1})\!=\!
\!\sum_{a\in W(b)} X_{-a}q^{-n_b(a)}
\!\!\!\mod \!\Si_+(b),
\end{align}
where $n_b(a)$ are from (\ref{ncbcoef}). One has: 
\begin{align}\label{alim1}
&\lim_{t\to 0} q^{-(b_{-},\rho_k)}
E^*_b(q^{c_\#})=
q^{-n_b(a)}X_{-a}(q^{u_c^{-1}(c_-)})\\
&=q^{-n_b(a)}q^{-(a,u_c^{-1}(c_-))}=
q^{-n_b(a)-(b_-,c_-)},\notag\\
&\frac{q^{\frac{b_{-}^2}{2}+
\frac{c_{-}^2}{2}}}
{v(c\!-\!b)h^0_b}
\,\bigl(\lim_{t\to 0} q^{-(b_{-},\rho_k)}
E^*_b(q^{c_\#})\bigr)=
\frac{q^{\frac{(c_{-}\!-b_{-})^2}{2}\,-n_b(a)}}
{v(c_{-}\!\!-\!b_{-})h^0_b}.\notag
\end{align} 
As above, we take only $a$ such that $X_a$ is present in 
$E_b^\dag$; we technically set $q^{-n_b(a)}=0$ otherwise. Also 
we can obviously replace here $v(c)$ by $v(c_-)$. Thus 
(\ref{thexpannn1}) becomes in the limit: 
\begin{align}
\label{thexpanfin1}
&\overline{E}_c\hat{\theta}_v=
\sum_{b\in P}\frac{q^{\frac{(c_--b_-)^2}{2}-\,n_b(u_c^{-1}(b_-))}}
{v(c_-\!\!-\!b_-)h^0_b}\overline{E}_b.
\end{align}

Assuming that $c=c_-$ we come to the setting of \cite{ChFB}. 
Indeed, $u_c=$id and $q^{-n_b(u_c^{-1}(b_-))}\!=\!0$ unless 
$b=b_-$. For $b=b_-\in P_-$, the polynomials 
$\overline{E}_{b_-}$ coincide with the corresponding {\em 
symmetric\,} bar-polynomials, and therefore the resulting 
decomposition is exactly as it was in \cite{ChFB}. 
Interestingly, this decompositions is  quite non-trivial when 
$c\not\in P_-$;  the numbers $n_b(c)$ are involved. 
\smallskip

\subsubsection{\sf The remaining cases}
Actually, the calculation is very similar in the other two 
cases. Let us proceed with (\ref{thexpann}). After the 
cancelations of $g$ as above, now $q^{-(\rho_k,c_-)}$  replaces 
$q^{-(\rho_k,b_-)}$ on the r.h. To proceed, let us set 
$E_c^{\dag\ast}\equal E^\dag_c(X^{-1};q^{-1})$. Then 
\begin{align}
\label{thexpannn}
&E_c^{\dag\ast}\hat{\theta}_v=
\sum_{b\in P}\frac{q^{\frac{b_{-}^2}{2}+
\frac{c_{-}^2}{2}}}
{v(b\!-\!c)h^0_b}
\,\bigl(\lim_{t\to 0} q^{-(c_{-},\rho_k)}
E^*_c(q^{b_\#})\bigr)\,E^{\dag\ast}_b.
\end{align}

The roles of $b$ and $c$ within $\lim\bigl(\cdots\bigr)$ are now 
exactly the opposite to those above. Therefore the monomial 
$X_{-a}$ from $E^*_c$ contributes to the limit if and only if 
$a=u_b^{-1}(c_-)$ and the final formula becomes: 
\begin{align}
\label{thexpanfin}
&E^{\dag\ast}_c\hat{\theta}_v=
\sum_{b\in P}\frac{q^{\frac{(b_--c_-)^2}{2}-\,
n_c(u_b^{-1}(c_-))}}
{v(b_-\!\!-\!c_-)h^0_b}E^{\dag\ast}_b.
\end{align}
\medskip

The analysis of the remaining (first) formula is also quite 
close to that for (\ref{thexpan1}). The cancelation of the 
$g$\~factors is exactly the same. The only change is that the 
polynomials $E_b$ replace $E_b^\ast$ there. Proposition 3.1 of 
\cite{ChO} states that 
\begin{align}\label{ebstar0}
E_b(X; q,t\!\to\!0)\!=\!
\!\sum_{a\in W(b)} \varsigma_b(a) X_{a}
\!\!\!\mod \!\Si_+(b),
\end{align}
where $\varsigma_b(a)=1$ if $u_a\ge u_b$ in the sense of the 
Bruhat order and $0$ otherwise. This replaces more involved 
(\ref{ebstar}). For the (unique) $X_{a}$ in $E_b$ (instead of 
$X_{-a}$ in $E_b^\ast$) contributing to the limit one now has: 
$$(b_-,\rho_k)=-(a,u_c^{-1}(\rho_k))=
-(u_c(a),\rho_k).$$

Therefore, $a=u_c^{-1}w_0(b_-)$, where we use that 
$-w_0(\rho_k)=\rho_k$, and 
\begin{align}\label{alim0}
&\lim_{t\to 0} q^{-(b_{-},\rho_k)}
E_b(q^{c_\#})\!=\!
X_{a}(q^{u_c^{-1}(c_-)})\!
=\!q^{(a,u_c^{-1}(c_-))}\!=\!
q^{-(b_-,c^\iota_-)}
\end{align} 
for $c^{\,\iota}\equal -w_0(c)$. Note that $u_{c^{\iota}}= 
w_0u_cw_0$ and $(c_-)^{\iota}=(c^{\,\iota})_-$. Finally: 
\begin{align}
\label{thexpanfin0}
\overline{E}_{c}\hat{\theta}_v&=
\sum_{b\in P}\varsigma_b(u_c^{-1}w_0(b_-))
\frac{q^{\frac{(b_--\,c^\iota_-)^2}{2}}}
{v(b_-\!\!+\!c_-)h^0_b}E^{\dag\ast}_b,\notag\\
\overline{E}_{c^\iota}\hat{\theta}_v&=
\sum_{b\in P}\varsigma_b(w_0u_c^{-1}(b_-))
\frac{q^{\frac{(b_--\,c_-)^2}{2}}}
{v(b_-\!\!-\!c_-)h^0_b}E^{\dag\ast}_b,
\end{align}
where we changed $c$ to $c^{\,\iota}$ in the second formula. 
%and use that $v(c^\iota)=v(-c)$. 

As an example, let us consider the case $c\in P_-$.  Then {\em 
all\,} $b\in P\,$ will occur in this summation because $u_c=$id. 
%\medskip

\subsection{\bf Iteration formulas}
Let us summarize the analysis above in the following theorem. 

\begin{theorem} \label{RR-THETA}
Let us fix $v\in\Pi'=Hom(\Pi,\C^*)$. We use the function 
$n_c(b)$ from (\ref{ncbcoef}); recall that we set 
$q^{-n_b(a)}\!=\!0$ if $X_a$ is not present in $E_b^\dag$. Also, 
$\varsigma_b(a)=1$ if $u_a\ge u_b$ in the sense of the Bruhat 
order in $W$ and $0$ otherwise, $c^{\,\iota}=-w_0(c)$ and 
$h_b^0=\!\!\prod_{[\al_i,j]\in \la'(\pi_b)} (1\!-\!q_i^j) \for 
b,c\in P$;  see (\ref{limeval}) and (\ref{jbset}) . One has: 
\begin{align}
\label{thexfin0}
&\overline{E}_{c^\iota}\hat{\theta}_v=
\sum_{b\in P}\varsigma_b(w_0u_c^{-1}(b_-))
\frac{q^{\frac{(b_--c_-)^2}{2}}}
{v(b_-\!\!-\!c_-)h^0_b}E^{\dag\ast}_b,\\
\label{thexfin}
&E^{\dag\ast}_c\hat{\theta}_v=
\sum_{b\in P}\frac{q^{\frac{(b_--c_-)^2}{2}-\,n_c(u_b^{-1}(c_-))}}
{v(b_-\!\!-\!c_-)h^0_b}E^{\dag\ast}_b,\\
\label{thexfin1}
&\overline{E}_c\hat{\theta}_v=
\sum_{b\in P}\frac{q^{\frac{(c_--b_-)^2}{2}-\,n_b(u_c^{-1}(b_-))}}
{v(c_-\!\!-\!b_-)h^0_b}\overline{E}_b.
\end{align}
In particular for $c=0$ (notice $b\in P_-$ in the second sum 
below): 
\begin{align}
\label{thefin0}
&\hat{\theta}_v=
\sum_{b\in P}\frac{q^{\frac{b^2}{2}}}
{v(b_-)h^0_b}E^{\dag\ast}_b=
\sum_{b\in P_-}\frac{q^{\frac{b^2}{2}}}
{v(-\!b_-)h^0_b}\overline{E}_b\,.
\end{align}
\vskip -1.5cm \sq 
\end{theorem}
\medskip

\subsubsection{\sf Some remarks}
{\em (i)\,} Formula (\ref{thexfin1}) for $v=1'$ is actually not 
new. It can be deduced from formula (4.43) from \cite{ChO}. The 
notation $\tga^{\ominus}$ there becomes $\theta$ in this paper, 
and $a_{b,u_c^{-1}}$ used in (3.43) is the coefficient of 
$X_{-u_c^{-1}(b_-)}$ in $E^*_b(t\to 0)=E_b^{\dag\ast}$; see 
Proposition 4.3 from \cite{ChO}. \vskip 0.2cm 

{\em (ii)\,} Let us discuss some symmetries of the expansions 
above. First of all, (\ref{thexfin}) formally follows from 
(\ref{thexfin1}). Indeed, the coefficient of 
$E_b^{\dag\ast}/h_b^0$ in this expansion equals $\lan 
E_c^{\dag\ast}\hat{\theta}_v \overline{E}_b 
\overline{\mu}_\circ\ran$. The coefficient of 
$\overline{E}_c/h_c^0$ in the r.h.s. of  (\ref{thexfin1}) is 
given by the same expression: $\lan \overline{E}_b\hat{\theta}_v 
{E}^{\dag\ast}_c \overline{\mu}_\circ\ran$. 

Similarly, 
%one obtains
%the symmetry $c\leftrightarrow b$ in (\ref{thexfin0}).
the coefficient of $E_b^{\dag\ast}/h_b^0$ in (\ref{thexfin0}) is 
$\lan \overline{E}_{c^\iota} \hat{\theta}_v \overline{E}_b 
\overline{\mu}_\circ\ran$, which is symmetric under 
$c^{\,\iota}\leftrightarrow b$. This readily gives a general 
relation: 
$$\varsigma_b(u_c^{-1}w_0(b_-))=
\varsigma_c(u_b^{-1}w_0(c_-)),\ b,c\in P.
$$

Let us outline a direct justification of the latter relation.  
Setting $u_c^{-1}w_0(b_-)=a=u_a^{-1}(b_-)$, Proposition 3.4  
from \cite{ChO} states that $X_a$ occurs in $\overline{E}_b$ 
(always with the coefficient $1$) if and only if $u_b\le u_a$ in 
for the Bruhat order in $W$. Equivalently, $u_b\le 
w_0u_c$ due to the minimality of $u_b$ modulo the centralizer of 
$b_-$ on the left. Generally, if $uv^{-1}=w_0$ and $l(u')<l(u)$ 
for $u'=us_\al$, which corresponds to deleting one simple 
reflection from the reduced decomposition of $u$, then 
$u'(v')^{-1}=w_0$ for $v'=vs_\al$, $l(v')>l(v)$, and therefore 
$v'>v$. Thus $u_b\le w_0u_c\Leftrightarrow w_0u_b\ge u_c$, which 
gives the required. 

\vskip 0.2cm 

{\em (iii)\,} The second expansion from (\ref{thefin0}), namely 
the expansion in terms of $\overline{E}_b$, is a special case of 
(\ref{thexfin1}). The first formula there, which is for 
$\hat{\th}_v$ via $E_b^\dag$, formally follows from the second  
due to the identity: 
\begin{align}
\label{thefin00}
\sum_{b\in W(c)}\frac{q^{\frac{b^2}{2}}}{h_b^0}
E^{\dag\ast}_b\ =\ 
\frac{q^{\frac{b_-^2}{2}}}{h^0_{b_-}}\,
\overline{E}_{b_-^\iota} \hbox{\, for\, } 
\ b_-^{\,\iota}\!=\!(-w_0(b))_-\!=\!-w_0(b_-), 
\end{align}
where $q^{b^2/2}=q^{b_-^2/2}$ and therefore both $q$\~powers can 
be omitted. Equivalently, $P^\dag_{b_-}=$ $ \sum_{b\in 
W(c)}E^{\dag}_b\,(h_{b_-}^0/ h_b^0)^\ast\ $ for the 
corresponding symmetric Macdonald polynomial $P_{b_-}^\dag= 
P_{b_-}(t\!\to\!\infty)$. As an example, let us provide the 
latter identity for $A_1$. Setting $X\!=\!X_{\om_1}$ for $n\ge 
0$\,: 
\begin{align*}
&E^\dag_{-n} + (1-q^{-n})E^\dag_n=P_{-n}^\dag,
\hbox{\, where\, }  E_n=E_{n\om_1};\hbox{\, letting\, }  n=2:\\
&(X^{-2}\!+\!\frac{X^2}{q^2}\!+\!
\frac{1\!+\!q}{q^2})+(1\!-\!\frac{1}{q^2})(X^2\!+\!1)
\!=\!P^\dag_{-2}\!=\!X^2\!+\!X^{-2}\!+\!\frac{1\!+\!q}{q}\,.
\end{align*}

Formula (\ref{thefin00}) results from (2.57) or general formula 
(4.19) from Proposition 4.2 in \cite{ChO}: 
$E_{w(b)}^\ast/h_{w(b)}\ =\ t^{l(w)/2}\,\tilde{\Psi}'_w 
(E_{b}^\ast/h_{b}) \for b=b_-, $   in the notation there. We 
obtain that the l.h.s. of (\ref{thefin00}) is $W$\~invariant and 
therefore proportional to $P_b^\ast=P_{\iota(b)}$ in the limit 
$t\to 0$; the coefficient of proportionality is simple. One can 
also make $t\to \infty$ in formula (3.3.15) from \cite{C101}, 
which expresses $P_{b_-}$ in terms of $E_{w(b_-)}$. 
%\smallskip

\subsubsection{\sf Using intertwiners}
It is instructional to check that applying the intertwiners to 
(\ref{thexfin1}) results in relations for $q^{n_c(b)}$ obtained 
above. Let $\zeta_c(b)=q^{n_c(b)}$ be the coefficient of $X_b$ 
in $E_c^\dag$. Then for $v=1'$, 
\begin{align}\label{thexfin11}
&\overline{E}_c\hat{\theta}=
\sum_{b\in P}q^{\frac{(c_--b_-)^2}{2}}
\zeta^*_b(u_c^{-1}(b_-))\overline{E}_b/h_b^0.
\end{align}

The first relation in (\ref{ncbrel}) gives that for 
$(b,\al_i)<0$, 
\begin{align}\label{balic}
\{\zeta_{s_i(b)}(u_c^{-1}(b_-))\!\neq\! 0\}&\Rightarrow
\{0\!<\!(u_c^{-1}(b_-),\al_i)\!=\!(b_-,u_c(\al_i))\}\\
&\Rightarrow \{u_c(\al_i)\!\in\! -R_+\}\Leftrightarrow
\{(c,\al_i)\!>\!0\}.\notag
\end{align}
Recall that the equality $(c,\al_i)=0$ here is equivalent to 
$u_c(\al_i)=\al_j$ for simple $\al_j$ such that $(\al_j,c_-)=0$, 
which contradicts to $(b_-,u_c(\al_i))>0$. See Proposition 
\ref{BSTAL} and Lemma \ref{VANISHL}. 

Following Section 4.4 from \cite{ChO}, 
\begin{align}\label{intbar}
\overline{T}'_i(\overline{E}_b)=
\begin{cases}
&\overline{E}_{s_i(b)} \hbox{\, if\, } (b,\al_i)>0, \\
&\overline{E}_b \hbox{\, if\, } (b,\al_i)\le 0, \\
\end{cases}
\end{align}
where $\overline{T}'_i=1+\overline{T}_i= 
1+(X_{\al_i}-1)^{-1}(s_i-1)\ (1\le i\le n)$. Explicitly: 
\begin{align}
\overline{T}_i'(X_b)=
\begin{cases}
X_{s_i(b)}+X_b, & \text{if }(b,\al_i)>0,\\
X_b, & \text{if }(b,\al_i)=0,\\
0, & \text{if }(b,\al_i)<0.
\label{intmodsi}
\end{cases}
\end{align}

Applying $\overline{T}_i'$ to (\ref{thexfin11}): 
\begin{align}\label{thexfin111}
&\overline{T}'_i\left(\overline{E}_c\hat{\theta}\right)=
\overline{T}'_i\left(\overline{E}_c\right)\hat{\theta}\\
=
\sum_{(b,\al_i)>0}&q^{\frac{(c_--b_-)^2}{2}}
\zeta^*_{b}(u_c^{-1}(b_-))\frac{\overline{E}_{s_i(b)}}{h_b^0}+
\sum_{(b,\al_i)\le 0}q^{\frac{(c_--b_-)^2}{2}}
\zeta^*_{b}(u_c^{-1}(b_-))\frac{\overline{E}_{b}}{h_b^0}.\notag
\end{align}
When $(c,\al_i)\le 0$, we arrive at the identities: 
\begin{align*}
&\frac{\zeta^*_b(u_c^{-1}(b_-))}{h_b^0}
+\frac{\zeta^*_{s_i(b)}(u_c^{-1}(b_-))}{h_{s_i(b)}^0}
=\frac{\zeta^*_{s_i(b)}(u_c^{-1}(b_-))}{h_{s_i(b)}^0} 
\hbox{\,\, if\,\, }(b,\al_i)> 0.
\end{align*}
Since $(s_i(b),\al_i)<0$ and $(c,\al_i)\le 0$, the relation\, 
$\zeta_{b}(u_c^{-1}(b_-))=0$\, follows from (\ref{balic}). 

Now let $(c,\al_i)>0$. Then $u_c=u_{s_i(c)}s_i$ (the product is 
reduced) and we obtain: 
\begin{align*}%\label{zetabc}
&\frac{\zeta^*_b(u_c^{-1}(b_-))}{h_b^0}
\!+\!\frac{\zeta^*_{s_i(b)}(u_c^{-1}(b_-))}{h_{s_i(b)}^0}
\!=\!\frac{\zeta^*_{s_i(b)}(s_iu_c^{-1}(b_-))}{h_{s_i(b)}^0}
\hbox{\,\, if\,\, }(b,\al_i)> 0.
\end{align*}
Collecting the terms with $s_i(b)$ in the r.h.s., we obtain the 
last two relations from (\ref{ncbrel}). 
%\medskip

\subsubsection{\sf Main Theorem}\label{SEC:MAINTH}
Let us now iterate the formulas from Theorem \ref{RR-THETA}. We 
set ${\mathbf v}=\{v_1,\ldots,v_p\}\in \Pi'$, 
$\hat{\theta}_{\mathbf v}=\hat{\theta}_{v_1}\cdots 
\hat{\theta}_{v_p}$, and ${\mathbf b}=\{b_k\in P, 1\le k\le 
p\}$. We will use the following system of notations: 

\begin{align*}
&m_c(b)\equal\, -n_c(u_b^{-1}(c_-))\,=-n_c(u_b^{-1}u_c(c)),\\
&\xi_b(c)\,\equal\, \varsigma_b(u_c^{-1}w_0(b_-))\,=\,
\varsigma_b(u_c^{-1}w_0 u_b(b)),\\
&(\al_i^\vee,b)^{\malt}\,\equal\, -(\al_i^\vee,b_-) 
\hbox{\, when\, } 
u_b^{-1}(\al_i)\!\in\! R_+\\ 
&\hbox{and\, \ } (\al_i^\vee,b)^{\malt}\equal
-(\al_i^\vee,b_-)-1 \hbox{\, otherwise}.
\end{align*}
Note that $(\al_i^\vee,b)^{\malt}\ge 0$; it is needed in the 
norms in the denominators, where we use (\ref{limeval}) and 
(\ref{ghgt0}): 
\begin{align}\label{limev}
\lan E_b E_c^{\dag\ast}\overline{\mu}_\circ\ran
=\de_{bc}h_b^0=\de_{bc}\prod_{j=1}^{(\al_i^\vee,b)^{\malt}}
(1-q_i^{j});
\end{align}
see (\ref{lambpi}). We will also set $b^-=b_-$ to improve the 
visibility of the formulas below. 

Recall that $q^{-m_c(b)}$ is the coefficient of 
$X_{u_b^{-1}(c_-)}$ in $E_c^\dag$. The monomiality claim is that 
$n_c(b)\in -\Z_+$ (unless this coefficient is $0$). Switching to 
$m_c(b)=-n_c(b)$ is quite natural here. We set 
$m_c(b)=+\!\infty$ and $q^{+\!\infty}=0$ if such $X$ is not 
present in $E_c^\dag$. Since $m_c(b)$ is a non-negative integer 
otherwise, the sum of $m$\~terms is $\infty$ if and only if at 
least one of the terms is $\infty$ ($+\infty$ to be exact). 
\smallskip

\begin{theorem} \label{MAINTHM}
For an arbitrary sequence $\mathbf v$ and $c\in P$, 
\begin{align}\label{pggdag}
E_c^{\dag\ast}\,\hat{\theta}_{\mathbf v}&=
\sum_{\mathbf b}\frac{
q^{\left((c^--b^-_1)^2+(b^-_1-b^-_2)^2+\ldots+
(b^-_{p-1}-b^-_p)^2\right)\!/2}}
{\prod_{i=1}^n\prod_{k=1}^{p-1}
\prod_{j=1}^{(\al_i^\vee,b_k)^{\malt}}
(1-q_i^{j})}\,\\
&\times \frac{q^{m_c(b_1)+m_{b_1}(b_2)+
\ldots+m_{b_{p-1}}(b_p)}}
{v_1(b_1\!-\!c)\,v_2(b_2\!-\!b_1)\cdots
v_p(b_p\!-\!b_{p\!-\!1})}
\,\frac{E_{b_p}^{\dag\ast}}{h_{b_p}^0},\notag\\
\label{pggbar}
\overline{E}_c\,\hat{\theta}_{\mathbf v}&=
\sum_{\mathbf b}\frac{
q^{\left((c^--b^-_1)^2+(b^-_1-b^-_2)^2+\ldots+
(b^-_{p-1}-b^-_p)^2\right)\!/2}}
{\prod_{i=1}^n\prod_{k=1}^{p-1}
\prod_{j=1}^{(\al_i^\vee,b_k)^{\malt}}
(1-q_i^{j})}\,\\
&\times 
\frac{q^{m_{b_1}(c)+m_{b_2}(b_1)+\ldots+m_{b_{p}}(b_{p-1})}}
{v_1(c\!-\!b_1)\,v_2(b_1\!-\!b_2)\cdots
v_p(b_{p\!-\!1}\!-\!b_p)}
\,\frac{\overline{E}_{b_p}}{h_{b_p}^0}.\notag
\end{align}
Recall that in these formulas and those below any $b_k,c$ inside  
$v(\ )$ can be replaced by $b^-_k,c^-$. We obtain that the 
coefficient of  $E_{b_p}$ in the first or the second expansion 
depends only on $b_p$, the unordered set $\mathbf v$ and initial 
$c$. The same holds for the following {\sf mixed} expansion: 
\begin{align}\label{pggmix}
\overline{E}_{c^{\iota}}\,\hat{\theta}_{\mathbf v}&=
\sum_{\mathbf b}\frac{
q^{\left((c^--b^-_1)^2+(b^-_1-b^-_2)^2+\ldots+
(b^-_{p-1}-b^-_p)^2\right)\!/\,2+(b^-_r,b^-_{r+1})}}
{\prod_{i=1}^n\prod_{k=1}^{p-1}
\prod_{j=1}^{(\al_i^\vee,b_k)^{\malt}}
(1-q_i^{j})}\,\\
&\times 
\frac{q^{m_{b_1}(c)+m_{b_2}(b_1)+\ldots+m_{b_{r}}(b_{r-1})}\,
\xi_{b_{r+1}}(b_{r})}
{v_1(c\!-\!b_1)\,v_2(b_1\!-\!b_2)\cdots
v_r(b_{r\!-1}\!-\!b_r)\,
v_{r\!+1}(b_{r}\!+\!b_{r\!+1})}\notag\\
&\times 
\frac{q^{m_{b_{r+1}}(b_{r+2})+m_{b_{r+2}}(b_{r+3})
+\ldots+m_{b_{p-1}}(b_{p})}}
{v_{r\!+2}(b_{r\!+2}-\!b_{r\!+1})\,
v_{r\!+3}(b_{r\!+3}-\!b_{r\!+2})\cdots
v_{p}(b_{p}-\!b_{p\!-1})}
\,\frac{E^{\dag\ast}_{b_p}}{h_{b_p}^0},\notag
\end{align}
where $0\le r\le p-1$. In this expansion, we switch from 
$\overline{E}$ to $E^{\dag\ast}$ at place $r+1$ using 
(\ref{thexfin0}). \vskip -0.5cm \sq 
\end{theorem}

%\medskip
\subsubsection{\sf Comments}
The following lemma is important to understand how far the 
summations above are from those over $P_-$. It readily results 
from the implication $\{b\in P_-, m_c(b)\neq \infty\} 
\Rightarrow \{c \in P_-, m_c(b)=0\}$. Here we use that $u_b=$id 
if and only if $b\in P_-$. We will set $c=b_0$ for the sake of 
uniformity. 

\begin{lemma}\label{mbcLEM}
Assuming that $b_r\in P_-$ in a nonzero term (a product) from 
the summation in (\ref{pggdag}), all previous $b_{s}$ $(0\le 
s\le r)$ must be then from $P_-$ in this product and 
$m_{b_{s-1}}(b_{s})=0$ for\, $1\le s'\le r$.  Similarly, if some 
$b_r$ belongs to $P_-$ in a nonzero product in (\ref{pggbar}), 
then $b_{s}\in P_-$ for\, $s\ge r$ and $m_{b_{s}}(b_{s+1})=0$ 
for\, $r\le s\le p-1$. \sq 
\end{lemma} 
\smallskip

Let us focus on (\ref{pggdag}). We set 
\begin{align}\label{Xidag}
\Xi_{p,\mathbf v}^{\,c,a}\equal 
\lan E_c^{\dag\ast}\hat{\theta}_{\mathbf v}
E_a\overline{\mu}_\circ\ran,
\end{align}
which is the coefficient of $E_a^{\dag\ast}/h_a^0$ in 
(\ref{pggdag}). 

Note that a somewhat different definition of theta-functions was 
used in \cite{ChFB}. Namely, we set for a collection  
$\boldsymbol{\varpi}=\{\varpi_k\subset \Pi, 1\le k\le p\}$: 
\begin{align}
&\theta_\varpi(X)\ \equal\ 
\sum_{b\in \varpi+Q} q^{(b,b)/2}X_b,\ \,
\hat{\theta}_{\boldsymbol{\varpi}}\equal
\frac
{\prod_{k=1}^p\hat{\theta}_{\varpi_k}}
{\lan \theta\overline{\mu}_\circ \ran^p}.
\label{thetapi}
\end{align}

Then the following modification of (\ref{pggdag}) is necessary: 
\begin{align}\label{pggdagpi}
E_c^{\dag\ast}\,\hat{\theta}_{\boldsymbol{\varpi}}&=
\sum_{\mathbf b}\frac{
q^{\left((c^--b^-_1)^2+(b^-_1-b^-_2)^2+\ldots+
(b^-_{p-1}-b^-_p)^2\right)\!/2}}
{\prod_{i=1}^n\prod_{k=1}^{p-1}
\prod_{j=1}^{(\al_i^\vee,b_k)^{\malt}}
(1-q_i^{j})}\,\\
&\times q^{m_c(b_1)+m_{b_1}(b_2)+
\ldots+m_{b_{p-1}}(b_p)}
\,E_{b_p}^{\dag\ast}/h_{b_p}^0,\notag\\
\hbox{where\,\ }
c\!-\!b_1&\!\in\!\varpi_1\!+\!Q,
\,b_1\!-\!b_2\!\in\!\varpi_2\!+\!Q,
\ldots,
b_{p\!-\!1}\!-\!b_p\!\in\! \varpi_p\!+\!Q.\notag
\end{align}
Accordingly, we must switch from $\Xi_{p,\mathbf v}^{\,c,a}$ in 
(\ref{Xidag}) to 
\begin{align}\label{Xidagpi}
\Xi_{p,\boldsymbol{\varpi}}^{\,c,\,a}\equal 
\lan E_c^{\dag\ast}\hat{\theta}_{\boldsymbol{\varpi}}
E_a\overline{\mu}_\circ\ran.
\end{align}
\vskip 0.3cm 

Let us now discuss Theorem 2.3 from \cite{ChFB}, which addresses 
the $PSL(2,\Z)$\~modularity, adjusting it to the 
$E^\dag$\~expansions from the theorem. The modularity 
(generally) occurs for $c\!=\!0$ and minuscule $a$. 

\begin{corollary}\label{MODUL-NS}
We assume that $v_k$ from $\mathbf v$ are trivial in 
$(P\cap Q^\vee)/Q$. Then 
$\Xi_{p,\mathbf v}^{\,0,a}\!\equal\! \lan 
\hat{\theta}_{\mathbf v} E_a\overline{\mu}_\circ\ran$ is a 
modular function of weight $0$ for $a\!=\!0$ or for minuscule 
$-a\in P_+$ with respect to some congruence subgroup of 
$PSL(2,\Z)$. Due to Lemma \ref{mbcLEM}, only $b_k\in P_-$ 
contribute to these summations. Namely, setting $\mathbf b'= 
\{b_1,\ldots,b_{p-1}\}\in P_-^{p-1}$: 
\begin{align}\label{coefdag}
\Xi_{p,\mathbf v}^{\,0,a}&=
\sum_{\mathbf b'}\frac{
q^{\left((b_1)^2+(b_1-b_2)^2+\ldots+
(b_{p-2}-b_{p-1})^2+(b_{p-1}-a)^2\right)\!/2}}
{\prod_{i=1}^n\prod_{k=1}^{p-1}
\prod_{j=1}^{-(\al_i^\vee,b_k)}
(1-q_i^{j})}\,\\
&\times \frac{1 }
{v_1(b_1)\,v_2(b_2\!-\!b_1)\cdots
v_{p\!-\!1}(b_{p\!-\!1}\!-\!b_{p\!-\!2}) 
v_p(a\!-\!b_{p\!-\!1})}.\notag
\end{align}
Using $\hat{\theta}_{\boldsymbol{\varpi}}$, the corresponding 
$\Xi_{p,\boldsymbol{\varpi}}^{\,c,\,a}$ for $c=0$ reduce to 
modular Rogers-Ramanujan type sums from Theorem 2.3 in 
\cite{ChFB} upon the restrictions $\varpi_k=\varpi_k+ Q^\vee/Q\, 
(1\le k\le p)$. 
\end{corollary}
{\it Proof.} We use relation (\ref{limev}) and then  follow 
\cite{ChFB}, mainly Section 2.3 there.  See also 
\cite{And,An,VZ,War,Za} and more recent \cite{CGZ,GOW}. \sq 
\smallskip

\subsubsection{\sf The case of 
\texorpdfstring{$A_1$}{\em A1}} We set $X=X_{\om_1}, 
t=t_{\sht}$, and denote $n\om_1$ for $n\in \Z$ simply by $n$; 
thus $s(n)=-n$ in this notation for $s=s_1$. One has: 
$q^{(n\om_1)^2/2}=q^{n^2/4}$ and 
$h_n^0=\prod_{j=1}^{|n|'}(1-q^j)$, where $|0|'=0$, 
$|\!-n|'=|n|=|n+1|'$ if $n\ge 0$. Then $u_{n+1}=s$, $u_{-n}=$id 
and the coefficient of $X^n$ in $E^\dag_{-n}$ is $q^{-n}$ for 
$n\ge 0$. Therefore for $b,c\in \Z$, 

\[
%\label{mcba1}
m_{c}(b)\!=\! -n_c\left(u_b^{-1}(\!-|c|)\right)\!=
\left\{
\begin{array}{ccc}
\ 0\,  &\hbox{for} &c\le 0, b\le 0 \hbox{\, or\, } c>0,b>0,\\
 |c|   &\hbox{for} &c\le 0, b>0,\\
 \infty&\hbox{for} &c>0, b\le 0.\,
\end{array}
\right.
\]
%\end{align}

Recall that we set $q^{\!\infty}=0$.   Finally, the characters 
$v$ for $A_1$ are $(\pm1)^n$. 

Let $\mathbf n'=\{n_1,\cdots,n_{p-1}\}$, $n_0=c,n_p=a$, and 
$\tilde{\mathbf n}^+=\{c,\mathbf n',a\}$. Then given any 
$v=\{v_1,\ldots,v_p\}$ and $a,c\in \Z$, we arrive at: 
\begin{align}\label{coefdag1}
\Xi_{p,\mathbf v}^{\,c,a}=&
\sum_{\mathbf n'}\frac{
q^{\left((|c|-|n_1|)^2+(|n_1|-|n_2|)^2+\ldots+
(|n_{p-2}|-|n_{p-1}|)^2+(|n_{p-1}|-|a|)^2\right)\!/4}}
{\prod_{k=1}^{p-1}\prod_{j=1}^{|n_k|'}
(1-q^{j})}\,\\
&\times \frac{q^{m_c(n_1)+m_{n_1}(n_2)+\ldots
+m_{n_{p-2}}(n_{p-1})+m_{n_{p-1}}(a)} }
{v_1(n_1-c)\,v_2(n_2\!-\!n_1)\cdots
v_{p\!-\!1}(n_{p\!-\!1}\!-\!n_{p\!-\!2}) 
v_p(a\!-\!n_{p\!-\!1})}\,.\notag
\end{align}
Note that for $a\le 0$, the terms are nonzero if and only if the 
set $\{n_0,\cdots, n_{p-1}\}$ is from $P_-$ and $c\le 0$, which 
matches Corollary \ref{MODUL-NS}. Similarly, if $c\!>\!0$, then 
$\{n_1,\cdots, n_{p}\}$ must be all from $1+\Z_+$. Generally, 
nonzero terms are exactly for $\tilde{\mathbf n}$ such that 
$\tilde{\mathbf n}=\mathbf n^-[r]\cup \mathbf n^+[r+1]$, where 
$$
\mathbf n^-[r]=\{n_0,\ldots,n_r\}\subset -\Z_+,\ 
\mathbf n^+[r+1]=\{n_{r+1},\ldots,n_{p}\}\subset 1+\Z_+
$$ 
for $-1\le r\le p$. I.e. the sequence $\tilde{\mathbf n}$  must 
have no single transition from strictly positive $n_k$ to 
non-negative $n_{k+1}$. We will call the corresponding $\mathbf 
n'$ {\em good\,} and set $\eta(\tilde{\mathbf n})=0$ unless 
$c\!\le\! 0$ and $a\!>\!0$; in the latter case let 
$\eta(\tilde{\mathbf n})=|n_r|$. Using this analysis: 
\begin{align}\label{coefdag11}
\Xi_{p,\mathbf v}^{\,c,a}=&\!\!\!
\!\!\sum_{\hbox{\tiny good\,}\mathbf n'}\!\!\frac{
q^{\left((|c|-|n_1|)^2+(|n_1|-|n_2|)^2+\ldots+
(|n_{p-2}|-|n_{p-1}|)^2+(|n_{p-1}|-|a|)^2\right)\!/4}}
{\prod_{k=1}^{p-1}\prod_{j=1}^{|n_k|'}
(1-q^{j})}\,\\
&\times \frac{q^{\eta(\tilde{\mathbf n})}}
{v_1(n_1\!-\!c)\,v_2(n_2\!-\!n_1)\cdots
v_{p\!-\!1}(n_{p\!-\!1}\!-\!n_{p\!-\!2}) 
v_p(a\!-\!n_{p\!-\!1})}\,.\notag
\end{align}

\setcounter{equation}{0} 
\section{\sc Demazure slices}\label{sec:Demazure}
In this section, we create representation-theoretical
tool to interpret formula (\ref{thefin0})
from Theorem \ref{RR-THETA}, which is actually the key
in there.

We use the results from \cite{FKM} and \cite{KL}, 
generalize them and provide some new proofs.  
Our approach is based on the associated 
graded of the filtration of integrable highest weight modules by 
its  thin (usual) and thick (upper) Demazure submodules. 

Following \cite{Kas94, AKT, Kat}, we compare 
the associate graded of the thick Demazure filtration
in a level one integrable modules of an
affine Lie algebra $\mathfrak g$ with the Weyl modules of 
(twisted) current algebras. Due to
\cite{FL,FMS} and the expansion formula from 
\cite{ChFB}, this construction results in 
symmetric Macdonald polynomials at $t\!=\!0$, which is 
Theorem \ref{G-DEM-CH}. Concerning the second one, we need a
proper version of the Demazure-Joseph functor 
from \cite{J} together with Section 4 from \cite{FKM}.
This connects the thick associated graded with 
the non-symmetric Macdonald polynomials at $t = \infty$ 
(Corollary \ref{E-D-DAG}); the $\mathrm{Ext}$-interpretation 
of the pairing (\ref{limev}) between the Macdonald 
polynomials at $t\!=\!0$ and $t\!=\!\infty$ is used here. 
\vskip 0.2cm

{\sf List of basic modules and functors.}
\label{sec:notations}
For the convenience of readers, we provide a list of basic 
modules and functors to be used in Sections \ref{sec:Demazure} 
and \ref{sec:exp-theta}, with the links to corresponding
subsections. 

The elements $u,w,\,\ldots\,$ will be now from $\widehat{W}$
(not from $W$ as above);
$[w]$ stands for the image of $w\in \widehat{W}$ in $\Pi$.
Through this part of the paper,
$M^{\vee}$ denotes the restricted dual of a 
$\widetilde{\mathfrak h}$-semisimple module $M$ with 
finite-dimensional weight spaces; for finite-dimensional
$M$, usual dual $M^{\ast}$ is sufficient. 

Note that $\mathrm{gr} ^w L$ in the table below
is a $\widetilde{\mathfrak b}^-$-module, 
while $\mathrm{Gr} ^b L$ is generally a 
$\mathfrak g_{\le 0}$-module. 

\comment{
\extra{I ADDED ${}_{[b]}$ and level-one. 
Please check: --- Sorry. $b \in \widehat{W}$, and 
hence I set $L ( \Lambda )$ to make things free from 
elements of Weyl groups. Also, the omission of scripts 
in $E^{\dag *}/h^0$ is INTENSIONAL as the actual 
expression in Corollary \ref{E-D-DAG} requires a 
couple of twists (that is just a notational issue and 
should be ignored here).}
}

\begin{align*}
& L\,=\,L ( \Lambda _{[w]}),\,\ L ( \Lambda_{[b]} ) 
 & \text{level-one integrable modules: } 
& \ref{uDm}\\
& L^w, L^b\! =\! L^{\pi_b} \subset\! L ( \Lambda_{[b]} ) 
& \text{thick(upper) Demazure modules: } & \ref{uDm}\\
& \mathrm{gr} ^w L\ (w\in \widehat{W}),\, \mathrm{Gr} ^b L 
& \text{associated graded (pieces) of $L$: } 
& \ref{subsub:demfilt}\\
& \mathbb D_b\!=\!\mathrm{gr} ^{\pi_b} L\,\text{ for } b\in P 
& \text{modules with {\sf gch}\,$=q^{\frac{b^2}{2}}\, 
\frac{w_{0} (E^{\dag *}_{b^{\iota}})}{h_{b}^0}$: } 
& \ref{subsub:demfilt}\\
& \mathscr D_i (i\!\ge 0),\, \mathscr D_w (w\!\in\!\widehat{W})
& \text{Demazure-Joseph functors: } 
& \ref{subsub:Jos}\\
& \mathbb W_b ( \subset \mathrm{Gr} ^{b_+} L ),\, b\in P  
& \text{generalized global Weyl modules: } 
&  \ref{subsub:formula}\\
& D_b^{\vee}  ( \subset L ( \Lambda_{[b]} )), b\in P 
& \text{dual of thin Demazure modules: } 
& \ref{subsub:orth}\\
& W_b \,\,( \text{ covered by } \mathbb W_b\, ) 
& \text{generalized local Weyl modules: } 
&  \ref{subsub:end}\\
& L ( b\!+\! k \Lambda_0 ) \text{ for } b\in P_+ 
& \text{level-$k$ integrable modules: } 
&\ref{sec:General}.\\
\end{align*}

%\begin{align*}
%&D_b \ : \ \hbox{\, thin Demazure module\ } &-\ 
%Formula (\ref{ })\\
%&\mathbb D_b\ : \ ....&-\ \hbSection \ref{}  and so on \\
%...
%\end{align*}.

\subsection{\bf Demazure slices}

\subsubsection{\sf Thick Demazure modules}\label{uDm}
We identify the cosets in  
$\widehat{W} / W$ with their minimal length 
representatives in $\widehat{W}$. Namely, each $w\in \widehat{W}$
can be uniquely represented in the form $w=bu$ for 
$b\in P, u\in W$. Then  $\pi_b$ is such a minimal representative
of $w$ in the notation from Proposition \ref{PIOM}. 
Following (\ref{hWPi}), we set from now on:
\begin{equation}
\widehat{W} \ni w=bu \mapsto [w] = \om_r \for r\in O
\hbox{\,\, such that\,\, } b-\om_r\in Q, 
\label{W-MINISCHULE}
\end{equation}
i.e. the images of $w$ and minuscule 
$\om_r$ (or $0$) coincide in $\widehat{W}/\tilde{W}=
\Pi$.  To simplify the notations,
we will use $w$ till the end of the
paper without ``hat" for the elements
in $\hat{W}$. Also, $\al$ will be generally {\em affine\,}
roots, unless in $\tal=[\al,j\nu_\al]$. 
Let $\widetilde{\mathfrak g}$ be the affine Kac-Moody 
algebra over $\C$ with $\widetilde{R}$ as its set of real roots 
and the degree operator $\mathsf{d}$, which corresponds
to $\mathsf{d}$ from (\ref{dform}). Its Cartan subalgebra
will be denoted by $\widetilde{\mathfrak h}$. See
\cite{Kac}, Chapter 7 and 8 here and below.

For each $\al \in \widetilde{R}$, the corresponding root 
space will be denoted by 
$\widetilde{\mathfrak g}_{\al}$. We also set 
$\mathfrak{sl} ( 2 )_{\al} \equal \widetilde{\mathfrak g}_{\al} 
\oplus [ \widetilde{\mathfrak g}_{\al}, 
\widetilde{\mathfrak g}_{-\al} ] \oplus 
\widetilde{\mathfrak g}_{- \al}$, that is a Lie subalgebra 
of $\widetilde{\mathfrak g}$ isomorphic to $\mathfrak{sl} ( 2 )$. 
For $0 \le i \le n$, we also set $\mathfrak{sl} ( 2 )_{i} \equal
\mathfrak{sl} ( 2 )_{\al_i}$. 
%\extra{K. I edited around and 
%kick out several definitions of $\mathfrak{sl} ( 2 )_{?}$.}
The triangular decomposition of $\widetilde{\mathfrak g}$ is:
$
\widetilde{\mathfrak g} = \widetilde{\mathfrak n} \oplus 
\widetilde{\mathfrak h} \oplus \widetilde{\mathfrak n}^-,
$
where the set of $\widetilde{\mathfrak h}$\~weights of 
$\widetilde{\mathfrak n}$ is $\widetilde{R}_+$ completed by  
the positive imaginary roots. We also set:
$\widetilde{\mathfrak b}^{+} \equal 
\widetilde{\mathfrak h} \oplus \widetilde{\mathfrak n}$,  
$\widetilde{\mathfrak b}^{-} \equal 
\widetilde{\mathfrak h} \oplus \widetilde{\mathfrak n}^-$.

The degree operator $\mathsf{d}$ is normalized
so that the non-positive degree 
part $\widetilde{\mathfrak g} _{\le 0}$ of 
$\widetilde{\mathfrak g}$ contains 
$\widetilde{\mathfrak b}^{-}$; let
$\widetilde{\mathfrak g}' _{\le 0} \equal [\widetilde{\mathfrak 
g} _{\le 0}, \widetilde{\mathfrak g} _{\le 0}]$. The
$\mathsf{d}$-zero-part of 
$\widetilde{\mathfrak g}' _{\le 0}$, 
denoted by $\mathfrak g$, is the simple Lie algebra 
corresponding to $R$. One has:
$$
\mathfrak g = \mathfrak n \oplus \mathfrak h \oplus 
\mathfrak n^-, \where \mathfrak{n} = \medoplus_{\, \al \in R_+\,} 
\widetilde{\mathfrak g}_{\al},\ 
\mathfrak{h}=\widetilde{\mathfrak h}\cap\mathfrak{g}. 
$$ 
We will also use 
$\mathfrak b^- \equal \widetilde{\mathfrak h} \oplus \mathfrak n^-$ 
and $\mathfrak g_0 \equal \mathfrak g + \widetilde{\mathfrak h}$.
The latter is the $\mathsf{d}$-zero-part of 
$\widetilde{\mathfrak g}$.
% One has: $\mathfrak h \cong P \otimes _{\Z} 
%\C \subset \widetilde{\mathfrak h}$. 

The affine roots $\al^\vee=\al/\nu_\al$ will be considered
in this and the next sections as coroots, i.e. 
the corresponding elements of $\mathfrak {h}$. 
Accordingly, $\widetilde{\al} = [\al, j \nu_{\al}] \in 
\widetilde{R}$\, and\, $\widetilde{\al}^\vee=\tal/\nu_\al$
will be interpreted as
\begin{equation}
\tal=\al+ j\nu_{\al}\, \de,\ \,\,
\widetilde{\al}^{*} = \al^{\vee} + j K\, \in\,
 \widetilde{\mathfrak h}
 \cap 
[\widetilde{\mathfrak g}, \widetilde{\mathfrak g}],
%\widetilde{\al}^{*} \in 
%\widetilde{\mathfrak h},
\label{aff-coroot}
\end{equation}
where $K$ is the standard central element of 
$\widetilde{\mathfrak g}$, $\delta \in \widetilde{\mathfrak h}^*$ 
is the (positive) primitive imaginary root. 

The basic level-one fundamental weight 
$\Lambda_0 \in \widetilde{\mathfrak h}^*$ is defined
as follows:
$$\La_0(\al_i^{*})= \delta_{i0}, \hbox{\,\, equivalently,\,\, }
\Lambda_0 ( K )=1, \Lambda_0 (\mathfrak h)=0.
$$
Here and further we will frequently use the identification
$\La(\tal^*)=\tal^*(\La).$ Also, note that $\de(K)=0$ and
$$
\widetilde{\al}^{*} ( \beta+z\de ) = \alpha^{\vee} ( \beta ) \for 
\widetilde{\al} = [\al, j \nu_{\al}] \in 
\widetilde{R},\  \beta \in R,\  z\in \C.
$$

For each $w \in \widehat{W}$, we set $\Lambda_w \equal  
w ( \Lambda_0 ) \in \widetilde{\mathfrak h}^*$. One has: 
\begin{equation}
\Lambda_w = w (\!(0)\!) + \Lambda _0 \hbox{\, mod\, } \Q\delta
\hbox{\,\, for the affine action from (\ref{afaction}).}
\label{Lambda}
\end{equation}

For $[w]$ from (\ref{W-MINISCHULE}),
$\Lambda_{[w]}$ is a level-one fundamental weight of 
$\widetilde{\mathfrak g}$.  The exact formula from
\cite{Kac} is:
$$
\Lambda_{w} = c + \Lambda_0-\frac{c^2}{2} \delta 
\in \Lambda_{[w]} + Q + \Z \delta.
$$
Thus we have $c=\om_r$ when $w = [w]=\om_r$ for some $r\in O$.

Let $\mathsf v_{[w]}$ be a unique up to proportionality
$\widetilde{\mathfrak n}$\~fixed 
vector of the corresponding integrable
irreducible level-one $\widetilde{\mathfrak g}$\~module denoted
by $L ( \Lambda_{[w]} )$.
For each $w \in \widehat{W}$, there exists a unique (up to
a scalar) vector $\mathsf v_{w} \in 
L ( \Lambda_{[w]} )$ of $\widetilde{\mathfrak h}$\~weight 
$\Lambda _w$. One has:
%By translating the highest weight 
%condition by the braid group action, we have} 
%\extra{ it is unique without next:} 
\begin{equation}
\widetilde{\mathfrak g}_{\al} \mathsf v_{w} = \{ 0 \} 
\text{\, for any\, }  \al \in 
w (\widetilde{R}_+).\label{EXT-INT}
\end{equation}

The {\em thick (upper)
Demazure modules\,} $L ^w$ are defined as
follows: 
$$
L ^w \equal  U ( \widetilde{\mathfrak n}^- ) \mathsf v_{w} \subset 
L ( \Lambda_{[w]} ), \where w\in \widehat{W};
$$
they are $\widetilde{\mathfrak b}^-$\~modules. 
%\extra{:I ADDED --- Thanks!}

By (\ref{Lambda}), the set of vectors $\{\mathsf v_{w}\}_{w \in 
\widehat{W}}$ is in bijection with $P$ via the
map $w\mapsto w(\!(0)\!)$. We restrict it to the elements
of the form $\pi_b\in \widehat{W}$; recall that 
$\pi_b$ map to $b$ under this map. Finally we set:
$$
L^b \equal  L^{\pi_b},\ \mathsf v_{b} \equal  
\mathsf v _{\pi_b} \in 
L ( \Lambda_{[\pi_b]} ), \where b\in P.
$$
Any thick Demazure modules are of the form  $L ^b$
for proper $b \in P$. 

By the triangular decomposition, each $L^w$ is a 
$\widetilde{\mathfrak h}$\~semisimple module. We have the Bruhat 
order $\le $ on $\widehat{W}$ defined as $v \le w$ if and only 
if $v \in \mathcal B ( w )$; see Propositions \ref{PIOM},
\ref{BSTAL}. 
Recall that if
 $W ( b ) = W ( c ) \subset P$ for $b, c \in P$ , then
$b \ll c$ by the 
partial order from (\ref{order}) if and only if $u_b < u_c$ 
for the Bruhat order.

\begin{lem}\label{COVER}
Let $w \in \widehat{W}$. For each $\al \in \widetilde{R}_+$, we 
have 
$$
\mathsf v_{s_\al w} \in \begin{cases} 
\widetilde{\mathfrak g}_{- \al}^{\,\,\al^{*} ( \Lambda _ w)}\, 
\mathsf v_{w} & 
(\al^{*} ( \Lambda _ w) > 0)\\ \C \mathsf v_{w} & (\al^{*} 
( \Lambda _ w) = 0)\\ 
\widetilde{\mathfrak g}_{\al}^{\,- \al^{*} 
( \Lambda _ w)}\, \mathsf v_{w} & 
(\al^{*} ( \Lambda _ w) < 0).\end{cases}.$$
Here and further 
$\widetilde{\mathfrak g}_{\al}^{\,m}$ 
($m \ge 1$) is the $m$-th multiplicative power of 
$\widetilde{\mathfrak g}_{\al}$ 
in the universal enveloping algebra 
$U ( \widetilde{\mathfrak g} )$.
In particular, the spaces in the right-hand side 
are always non-zero and one-dimensional.
\end{lem}

{\it Proof.}
Since $L ( \Lambda _{[w]} )$ is an integrable 
$\widetilde{\mathfrak g}$\~module, the 
$\widetilde{\mathfrak h}$\~eigenvalue spaces 
for any weights from $\widetilde{W} 
( \Lambda _{[w]} )$ are of dimension one. Also,
if $\mathsf v_{w}$ is non-zero, then so is 
$\mathsf v_{s_\al w}$; consider the action of 
$\mathfrak{sl} ( 2 )_{\al}$. This 
gives the required.  \sq 

\begin{cor}\label{DEM-INCL}
For each $w, u \in \widehat{W}$, we have $w \le u$ if and only if 
$L^u \subset L^w$. 
\end{cor}

{\it Proof.} We first prove the ``only if" 
part of the corollary.
By \cite[Section 2.2]{BB}, there exist 
$w \le x < u$ such that $u^{-1} x$ is a reflection 
and $\ell ( x ) = \ell ( u ) - 1$. Here and below
we use the Bruhat order; see also \cite{Hu} here and below.
Then
Lemma \ref{COVER} implies that $L^u \subset L^x$;
continuing we obtain the ``only if" part. 

The ``if" part is as follows. 
We assume $L^u \subset L^w$ and need to prove that $u \ge w$.
The inclusion $L^u \subset L^w$ gives that $\mathsf v_{u} \in 
L^{w}$. 
Then we use that the module $L^{w}$ is stable with respect to the
action of $\mathfrak{sl} ( 2 )_i$ corresponding 
to $\al_{i}$, where $0 \le i \le n$, assuming 
that $s_{i} w > w$. The relation (\ref{EXT-INT}) is
applied and the PBW theorem.

Let us provide some details.
One has that  $\mathsf v_{s_{i}u}, \mathsf v_{u}$ 
belong to a single $\mathfrak{sl} ( 2 )_i$\~string 
by Lemma \ref{COVER}. Therefore  
$\mathsf v_{u} \in L^w$ implies that  
$\mathsf v_{s_{i}u}, \mathsf v_{u}$ belong to a single 
$\mathfrak{sl} ( 2 )_i$\~string in
$L^{\min \{ s_i w, w \}}$ and that
$$\mathsf v_{s_{i}u}, \mathsf v_{u} \in 
\begin{cases} L^{s_{i} w} & 
(s_{i} w < w)\\ L^{w} & (s_{i} w > w). 
\end{cases}$$

Now let us assume that $w \not\le u$ and $\mathsf v_{u} \in 
L^{w}$ and prove that this is impossible.
Here $u$ can be taken minimal  satisfying these two
conditions for the Bruhat order on $\widehat{W}$. 
One has either $\ell ( u ) = 0$ or 
$\ell ( u ) > 0$. In the latter case, there exists 
$0 \le i \le n$ such that $\ell ( s_i u ) = \ell ( u ) - 1$.

If $\ell(u)=0$, then $\mathsf v_{u}$ has to be  
the highest weight vector of $L ( \Lambda_{[u]} )$. Hence, 
$\mathsf v_u \in L^w$ implies $L^w = L^u$ and 
therefore $w = u$. 
Thus, this case cannot occur due to the assumption $w 
\not\le u$. 

If $\ell ( s_i u ) = \ell ( u ) - 1$, then $s_{i} u < u$. 
The minimality of $u$ gives that  
$s_{i} u > \min \{ w, s_{i} w \}$ for any $0 \le i \le n$.
There are two possibilities here:
$s_{i} w > w$ or $s_{i} w < w$. If $s_{i} w > w$, then 
$u > s_i u > w$, which contradicts 
to $w \not\le u$. 
If $s_{i} w < w$, then  $s_i u > s_{i} w$
and  again $u > w$ due to \cite[Proposition 2.2.7]{BB}
(see also \cite{Hu}).

Finally, we  
conclude that there is no pair $u, w \in \widehat{W}$ such
that $w \not\le u$ and $\mathsf v_{u} \in L^{w}$,
which proves the ``if" part. \sq

\subsubsection{\sf Filtrations of 
\texorpdfstring{$L ( \Lambda_{[b]} )$}{L}}
\label{subsub:demfilt}

The quantum analogue of $L^w$ admits a global base in the sense 
of Kashiwara \cite[Section 4, Proposition 4.1]{Kas94}. 
The Littelmann path model \cite{Lit95} and the interpretation of 
its initial and final directions from
\cite[Theorem 6.23]{AKT} results
in the following theorem. We note that \cite{Kas94} and 
\cite[Corollary 2.18]{Kat} give its independent proof. 

\begin{thm}\label{DEM-DIST}
For each $S \subset \widehat{W}$, there exists  
$S' \subset \widehat{W}$ such that 
$$\bigcap_{w \in S} L^{w} = \sum_{u \in S'} L^{u}.$$
\end{thm}
\vskip -0.7cm
\sq

For each $w \in \widehat{W}$, let 
$\mathrm{gr}^{w} L \equal  L^{w} / \sum_{u > w} L^{u},$
where the quotient is well-defined due to Corollary 
\ref{DEM-INCL}. 

Let us use Proposition \ref{BSTAL}, $(i)$.
See also \cite[(2.7.5) and (2.7.11)]{M3}.
One has that 
$u \ge w$ implies $u (\!(0)\!) \preceq w (\!(0)\!)$ 
for $u, w \in \widehat{W}$.
%under the notation of Section \ref{DefOrder}. 
For each $b \in P_+$, we set 
$$\mathrm{Gr}^{b} L \equal  \left( L^{b} + 
\sum_{c \prec b_-} L^{c} 
\right) / \sum_{c \prec b_-} L^{c}.$$
We will denote the image of $\mathsf v_b$ in $\mathrm{Gr}^{b} L$ 
by the same 
letter. The $\mathfrak h$\~weight of $\mathsf v_{b}$ is $b$; see
(\ref{Lambda}). 

\begin{cor}\label{GR-PROP}
For every $u, w \in \widehat{W}$, the vector subspace 
$$\left( L^{u} \cap L^{w} \right) / \sum_{x > w} L^{x} \subset 
\mathrm{gr}^{w} L$$
is either $\mathrm{gr}^{w} L$ or $\{ 0 \}$. 
\end{cor}

{\it Proof.}
By Theorem \ref{DEM-DIST}, we have 
\begin{equation}
L^{u} \cap L^{w} = \sum_{x \in S} L^{x},\label{INT}
\end{equation}
where $S \subset \widehat{W}$. Corollary \ref{DEM-INCL} gives 
that $x \ge w$ for each $x \in S$ and, moreover,
$x > w$ when $u \not\le w$. 

The space $L^{u} \cap L^{w}$
belongs to $\sum_{x > w} L^{x}$ 
when $u \not\le w$, so the quotient
is  $\{ 0 \}$ in this case.
Otherwise, $L^{w} \subset L^{u}$ and the quotient
is $\mathrm{gr}^{w} L$ as 
stated. \sq 

Let us introduce the $\widetilde{\mathfrak b}^-$\~modules
$\mathbb D_b \equal \mathrm{gr}^{\pi_b} L$, 
where the latter are
$\mathrm{gr}^{w} L$ above for %\extra{: I CHANGED THIS --- Thanks!}
$w=\pi_b$. By Proposition \ref{PIOM}:
$$\{\mathrm{gr} ^w L \}_{w \in \widetilde{W}} = 
\{ \mathbb D_b \}_{b \in P}.$$

\begin{cor}\label{G-GR-PROP}
For every $b \in P_+$, the vector space $\mathrm{Gr}^{b} L$ 
admits a $\widetilde{\mathfrak g}_{\le 0}$\~action. Moreover, 
 there is a natural finite filtration of
$\mathrm{Gr}^{b} L$ by $\{\mathbb D_c \}_{c \in W(b)}$, 
%\extra{By some or all?--I think 
%``\underline{occurs} with multiplicity one" is stronger 
%than just all of them appears. Nevertheless, I have edited.}
where every $\mathbb D_c$ occurs exactly once. 
We will say that $\mathrm{Gr}^{b} L$ is filtered by 
$\{\mathbb D_c \}_{c \in W(b)}$ with multiplicities
one. %\extra{: I ADDED THIS--- Thanks}
\end{cor}

{\it Proof.}
We have $\mathfrak n\, \mathsf v_b = 0$ by construction. 
In particular, 
$U ( \mathfrak g ) \mathsf v_b$ is a finite-dimensional 
$\mathfrak g$\~module with the highest weight $b$. Therefore, 
the PBW theorem 
implies that the $\widetilde{\mathfrak b}^-$\~action on $L^{b}$ 
extends to the $\widetilde{\mathfrak g}_{\le 0}$\~action. 
Then we use Corollary \ref{GR-PROP}, which gives that
the $\widetilde{\mathfrak b}^-$\~module $\mathrm{Gr}^{b} L$ 
has a natural finite filtration by $\{ \mathrm{gr}^{w} L \}_{w \in 
\widehat{W}}$. Since $w \in \widehat{W}$ such that 
$\mathrm{gr}^{w} L$ appears in this filtration must 
satisfy $w \ge \pi_b$ and $w \not\ge \pi_c$ for every $b \succ c 
\in P_+$, we obtain the required. \sq 

\subsubsection{\sf Characters of Demazure slices}

Recall that the $U ( \widetilde{\mathfrak g} _{\le 0} )$\~module 
$\mathrm{Gr}^{b} L$ has a cyclic vector 
$\mathsf v_{b}$ ($b \in P_+$) of 
$\mathfrak h$\~weight $b$. 

\begin{prop}\label{DEF-W}
Let $b \in P_+$. Consider the cyclic $U ( \widetilde{\mathfrak 
g}'_{\le 0} )$\~module $\mathbb W_{b}'$ generated by the 
cyclic vector $v$ 
subject to the following relations: 
\begin{enumerate}
\item $H v = b ( H ) v$ for each $H \in \mathfrak h$; 
\item $\widetilde{\mathfrak g}_{\al} v = 0$ for each $\al 
    \in \widetilde{R} \cap ( R_+ + \Z_{\le 0} \delta )$; 
\item $\widetilde{\mathfrak g}_{-\al} ^{\,\,\al^{*} ( b )\!
+\! 1}\, v = 0$ when $\al = [\beta, 0]$ or 
$[\beta, \nu_{\beta} ]$ 
    for some $\beta \in R_+$. 
\end{enumerate}
Then, $\mathbb W_{b}'$ maps surjectively onto 
$\mathrm{Gr}^{b} L$ as 
$U ( \widetilde{\mathfrak g}'_{\le 0} )$\~modules. 
\end{prop}

{\it Proof.}
Setting,
$q ( b + \Lambda_0 + m \delta ) \equal \frac{b^2}{2} - m$,
the weights $\Lambda_b$ for $b \in P$ satisfy the following 
(hyperbolic) equation: $q ( b + \Lambda_0 + m \delta )=0$.
%on $\widetilde{\mathfrak h}^*$. 
Moreover, they are exactly solutions of this equations
from all $\widetilde{\mathfrak h}$\~weights of 
$\medoplus_{\,\pi \in \Pi} 
L ( \Lambda_{\pi} )$, which generally satisfy the inequality 
$q ( b + \Lambda_0 + m \delta )\le 0$ (i.e. are inside the
corresponding hyperboloid).

Thanks to Lemma \ref{COVER} and the construction of 
$\mathrm{Gr}^{b} L$, the only $\widetilde{\mathfrak h}$\~weights 
$\Lambda$ of $\mathrm{Gr}^{b} L$ satisfying
$q ( \Lambda ) = 0$ 
are $\Lambda_{w b}$ for $w \in W$. One has $q<0$ for 
all the other weights in  $\mathrm{Gr}^{b} L$.
Hence, the cyclic vector 
$v=\mathsf v_{b} \in \mathrm{Gr}^{b} L$ 
satisfies the conditions
$(1),(3)$ above and the following modification of $(2)$:

\ (2')\  $\widetilde{\mathfrak g}_{\al} v = 0$ when $\al = 
    [\beta, 0], [\beta, - \nu_{\beta} ]$ for some $\beta \in 
    R_+.$

Therefore, it suffices to check that  condition $(2)$ results 
from $(1), (3)$ and $(2')$ to prove $(1,2,3)$ from the proposition.
%\extra{So it is somewhat different from those in the claim.
%I don't understand what you are doing below. Contradiction
%to what? Please fix this. --- The condition and its proof are 
%edited in the above.}

Upon the restriction to
$\mathfrak{sl} ( 2 )_{\al} \otimes_\C \C [z]$-calculation, we 
obtain that
\begin{equation}
\widetilde{\mathfrak g}_{\al + n \nu_\al \delta} 
\cdot v = 0 \hskip 3mm n 
\le 0\label{loop}
\end{equation}
for each $\al \in R_+$; see e.g. \cite{FMO}. Therefore,
  
$$\widetilde{\mathfrak g}_{\al} v = 0 \hskip 5mm \text{for each} 
\hskip 5mm \al \in \widetilde{R} \cap ( R_+ + \Z_{\le 0} 
\delta ),$$
which gives the required. The surjectivity
$\mathbb W_{b}'\to \mathrm{Gr}^{b} L$ readily follows.
\vskip -.5cm
\sq 

\smallskip

Next, we use that the $\widetilde{\mathfrak g}'_{\le 0}$\~action 
in $\mathbb W_{b}'$ can be naturally extended
to the $\widetilde{\mathfrak g}_{\le 0}$\~action 
by the formula $K v = 0$. The grading of the cyclic vector $v$
is $0$; we put $\mathsf{d} v = 0$. 

Let $\mathfrak B$ be the category of finitely generated 
$U ( \widetilde{\mathfrak b}^{-} )$\~modules with semisimple
$\widetilde{\mathfrak h}$\~action and such that every weight 
space is finite dimensional with its weight in $P \oplus \Z 
\Lambda_{0} \oplus \frac{1}{e} \Z \delta \subset 
\widetilde{\mathfrak h}^*$. As in Section \ref{SEC:POL}, here
$e$ is the minimal positive 
integer satisfying  $e(P,P) /2\subset \Z$.

For each $M \in \mathfrak B$, we set (formally):
$$\mathsf{gch} \, M \equal\!\! \sum_{c - m \delta \in P \oplus 
\frac{1}{e} \Z \delta} 
q^{m} X_{c} \cdot \dim \, \mathrm{Hom}_{\mathfrak h 
\oplus \C \mathsf{d}} 
( \C_{c - m \delta}, M ).$$

We put $f \le g$ for two polynomials 
$f, g \in \Z[q^{\pm 1/e}] [X_b,b\in P]$ for $X_b$ 
from (\ref{Xdex}) if this inequality holds 
{\em coefficient-wise\,}, i.e. for all pairs of
corresponding (integer) coefficients of the monomials 
$q^{m/e} X_b$ ($m \in \Z, b\in P$) in $f$ and $g$.  

\begin{thm}\label{G-DEM-CH}
For each $b \in P_+$ and $\mathbb W'_{b}$
from Proposition \ref{DEF-W}, we have 
$$q^{- \frac{b^2}{2}} \mathsf{gch} \, \mathrm{Gr}^{b} L = 
\mathsf{gch} \, \mathbb W'_{b} = 
\frac{\overline{E}_{b_{-}}}{h^0_{b_{-}}}.$$ 
\end{thm}

{\it Proof.}
We use \cite[Definition 2]{FL}, \cite[3.6]{CFS}, and 
\cite[(3.3)]{FMS}. Proposition \ref{DEF-W} implies that 
$\mathbb W_{b}'$ is a quotient of the 
{\em global Weyl modules\,} there.  Then
\cite[Proposition 4.3]{CI} implies the inequality
\begin{equation}
\mathsf{gch} \, \mathbb W_{b}' \le 
\frac{\overline{E}_{b_{-}}}{h^0_{b_{-}}}.
\label{INEQ-WE}
\end{equation}

Since $\mathsf{v} _{b}$ has 
$\mathsf{d}$-degree $- \frac{b^2}{2}$ in $L ( \Lambda_{[b]} )$, 
(\ref{INEQ-WE}) results in
\begin{align*}
&\mathsf{gch} \, L ( \Lambda_{[b]} ) = 
\sum _{c \in ( b + Q ) \cap P_{+}} 
\mathsf{gch} \, \mathrm{Gr}^{b} L \\
&\le \sum _{c \in ( b + Q ) \cap P_{+}} 
q^{\frac{c^2}{2}} \mathsf{gch} \, \mathbb W_{c_{-}}' \le 
\sum _{c \in ( b + Q ) \cap P_{+}} \frac{q^{\frac{c^2}{2}} 
\overline{E}_{c_{-}}}{h^0_{c_{-}}}.
\end{align*}
Here we employ
Proposition \ref{DEF-W} and (\ref{INEQ-WE}). Using now 
(\ref{thefin0}), 
$$\mathsf{gch} \, L ( \Lambda_{[b]} ) = 
\sum _{c \in ( b + Q ) \cap P_{+}} 
\frac{q^{\frac{c^2}{2}} \overline{E}_{c_{-}}}{h^0_{c_{-}}}.$$
Therefore the inequality in (\ref{INEQ-WE}) 
is actually an equality. \sq 

\subsection{\bf Demazure slices and 
\texorpdfstring{\mathversion{bold}$E_b^{\dag}\,(b\in P)$\,}{\em 
E-polynomials}} 

\subsubsection{\sf Demazure-Joseph functors}\label{subsub:Jos}

The main reference here is \cite{J}. It is for 
semi-simple Lie algebras, but the construction there can
be extended to our (affine, twisted) case. This is what we 
are going to do now.

For each $b \in P$ and $k \in \Z, m \in (1/e) \Z$, let
$$M ( b + k \Lambda_0 + m \delta ) \equal  
U ( \widetilde{\mathfrak g}) 
\otimes_{U ( \widetilde{\mathfrak b}^+)} \C _{b + 
k \Lambda_0 + m \delta},$$
where $\C_{b + k \Lambda_0 + m \delta}$ is the natural 
$\widetilde{\mathfrak b}^-$\~module for the
$\widetilde{\mathfrak h}$\~weight
$b + k \Lambda_0 + m \delta$. 
This is the usual definition of {\em Verma modules\,}
of $\widetilde{\mathfrak g}$. 

\begin{prop}[\cite{CG}]\label{BPROJ}
For any $b \in P$, $k \in \Z,\, m \in (1/e) \Z$,\,
 the Verma module 
$M ( b + k \Lambda_0 + m \delta ) $ viewed as a 
$\widetilde{\mathfrak b}^-$\~module, is the projective cover of 
$\C _{b + k \Lambda_0 + m \delta}$ in $\mathfrak B$. \sq 
\end{prop}

For each $0 \le i \le n$, we set $\mathfrak p_i^- \equal  
\widetilde{\mathfrak g} _{\al_i} \oplus \widetilde{\mathfrak 
b}^-$. A $\mathfrak p_i^-$\~module is said to be 
$\mathfrak{sl} ( 2 )$\~integrable if it is 
$\widetilde{\mathfrak h}$\~semisimple 
and is a direct sum of finite-dimensional modules of 
$\mathfrak{sl} ( 2 )_i$. 

Let us introduce some
general terminology.
For an abstract Lie algebra $\mathfrak L$ and its 
finite-dimensional
reductive Lie subalgebra 
$\mathfrak r$, $U ( \mathfrak r )$\~semisimple 
$U ( \mathfrak L )$\~modules will be called
 $(\mathfrak L, \mathfrak r)$\~modules. 

Given $0 \le i \le n$ and a $( \widetilde{\mathfrak b}^-, 
\widetilde{\mathfrak h} )$\~module $M$, consider the $U ( 
\widetilde{\mathfrak b}^- )$\~module $\mathscr D_i ( M )$ 
obtained as the maximal $\mathfrak{sl} ( 2 )_i$\~integrable 
quotient of $U ( \mathfrak p_i^- ) \otimes _{U ( 
\widetilde{\mathfrak b}^- )} M$. It is straightforward to see 
that $\mathscr D_i ( M )$ is $\widetilde{\mathfrak 
h}$\~semisimple. The correspondence $M\mapsto
\mathscr D_i(M)$ gives rise to a 
functor, which is usually called the {\em Demazure-Joseph functor}. 

\begin{thm}[\cite{J}]\label{JOS}
The functors $\{ \mathscr D_{i} \}_{0 \le i \le n}$ satisfy
the following.
\begin{enumerate}
\item Each $\mathscr D_{i}$ is right exact.
\item For $i,j \in I$ such that $(s_{i}s_{j})^{m} = 1$, one
    has 
$$\overbrace{\mathscr D_{i} \mathscr D_{j} \cdots}^{m} \cong 
\overbrace{\mathscr D_{j} \mathscr
D_{i} \cdots}^{m}.$$
\item There is a natural morphism $\mathrm{Id} \to 
    \mathscr D_{i}$.
\item For a $\mathfrak{sl} ( 2 )$\~integrable 
%\extra{Please do not use $\mathfrak sl$. It does not 
%look nice to change the font from the first letter to the 
%second letter.} 
$\mathfrak p_i^-$\~module $M$: $\mathscr D_{i} ( M ) \cong 
    M$. In particular, $\mathscr D_{i}^{2} \cong 
    \mathscr D_{i}$.
\item For a $\mathfrak{sl} ( 2 )$\~integrable $\mathfrak 
    p_i^-$\~module $M$ and a $\widetilde{\mathfrak 
    b}^-$\~module $L$:
$$\mathscr D_i ( L \otimes M ) \cong \mathscr D_i ( L ) 
\otimes M.$$
\end{enumerate}
\end{thm}
\vskip -0.7cm \sq

By Theorem \ref{JOS}, we can consider the {\em left derived\,} 
functor $\mathbb L \mathscr D_i$ in the category of $( 
\widetilde{\mathfrak b}^-, \widetilde{\mathfrak h} )$\~modules. 

For each $w \in \widetilde{W}$ with a reduced expression $w = 
s_{i_1} s_{i_2} \cdots s_{i_\ell}$, let 
$$\mathscr D _w \equal  \mathscr D_{i_1} \circ \mathscr D_{i_2} 
\circ \cdots 
\circ \mathscr D_{i_\ell}.$$
Thanks to the same theorem, $\mathscr D _w$ does not depend on 
the particular choice of the reduced expression of $w$. 

For $M$ with finite-dimensional $\widetilde{\mathfrak h}$\~weight 
spaces, we denote by $M^{\vee}$ its restricted dual (i.e. the 
direct sum of the duals of the weight spaces). Clearly, 
$M^{\vee}$ is again a $( \widetilde{\mathfrak b}^-, 
\widetilde{\mathfrak h} )$\~module with finite-dimensional weight 
spaces, and we have $M \cong ( M^{\vee} ) ^{\vee}$. We set 
$\mathscr D^{\sharp}_i \equal  \vee \circ \mathscr D_i \circ \vee$. 

\begin{thm}[\cite{FKM} Proposition 5.7]\label{D-ADJ}
For any two $\widetilde{\mathfrak h}$\~semisimple $U ( 
\widetilde{\mathfrak b}^- ) $\~modules $M$ and $N$ and 
for $0 \le i \le n$, one has:
$$\mathrm{Ext}^p_{( \widetilde{\mathfrak b}^-, 
\widetilde{\mathfrak h} )} 
( \mathbb L \mathscr D_i ( M ), N ) \cong \mathrm{Ext}
^p_{( \widetilde{\mathfrak b}^-, \widetilde{\mathfrak h} )} 
( M, \mathbb R \mathscr D_i^{\sharp} ( N ) ),$$
where $\mathrm{Ext}^p_{( \widetilde{\mathfrak b}^-, 
\widetilde{\mathfrak h} )} ( \,\cdot\,, \,\cdot\, )$ is the 
relative extension $($see e.g. 
\cite[I\!I\!I]{Kum}$)$.\sq 
\end{thm}
\vskip 0.2cm

Recall for the following Theorem and below,
that  $\prec$ is the orderings defined
in (\ref{succ}) and $<$ is the Bruhat order
from Proposition \ref{BSTAL}.
%\extra{K. Why it is inserted? see also the corrections of 
%notation in the next subsubsection.}

\begin{thm}\label{UDC}
Let $w \in \widetilde{W}$ and $0 \le i \le n$. If $s_i w 
< w$, then $\mathscr D_i ( L^w ) = L ^{s_i w}.$
If $s_i w > w$, then one has:
$\mathscr D_i ( L^w ) = L ^{w}.$
Moveover, $\mathbb L^{<0} \mathscr D_i ( L ^{w} ) = 
\{ 0 \}$. 

\end{thm}

{\it Proof.}
Using quantum group,  this result follows from its
analog of \cite[Proposition 3.3.4]{Kas94};
see also Theorem \ref{DEM-DIST} above. 
%This is not explicit in \cite[Section 4]{Kas94}), but follows
%the same lines. 
In the geometric approach from \cite[Theorem 2.15]{Kat}, 
on can obtain this claim as a corollary of 
\cite[Lemma 3.1 and Proposition 3.2]{KS}. \sq 

\subsubsection{\sf More on Demazure-Joseph functors}
\label{subsub:formula}

For any $b \in P_+$, let $\mathbb W_c$ be the image of 
$L^{c}$ in $\mathrm{Gr}^b L$, where $c \in W ( b )$. 
It is a $\widetilde{\mathfrak b}^-$\~module. 
%\extra{: I ADDED --- Thanks!}

By Corollary \ref{DEM-INCL} and Corollary \ref{G-GR-PROP}, 
$\mathbb W_c \cong 
\mathbb D_c$ as $\widetilde{\mathfrak b}^-$\~modules 
when $c = b_-$. Also,
$\mathbb W_c \cong \mathbb W'_c \otimes 
\C _{\Lambda_0 - \frac{c^2}{2}\delta}$ when $c = b_+$;
combine here Theorem \ref{G-DEM-CH} 
and Proposition \ref{DEF-W}.

\begin{lem}\label{E-U-STR}
For each $b \in P_+$, the module $\mathbb W_b$ admits a finite 
filtration by $\mathbb D_c$ (as constituents)
with $c \in W(b)$ such that each of them appears exactly once;
we say that $\mathbb W_b$ is filtered by 
$\{ \mathbb D_c \}_{c \in W ( b )}$ 
with multiplicities one.  
%\extra{: I edited. --- I complimented a little.}
%\extra{some $c$ or all? -- I have edited.}
\end{lem}

{\it Proof.}
Apply Corollary \ref{G-GR-PROP} to $\mathbb W_c \cong 
\mathrm{Gr}^b L$. 
\sq 
\vskip 0.1cm

\begin{prop}\label{GWEYL}
Let $c \in P$ and $0 < i \le n$. Then
$$\mathscr D_i ( \mathbb W_c ) = \begin{cases} 
\mathbb W_{s_i ( c )} & 
(s_i ( c ) \succeq c),\\ \mathbb W_{c} & (s_i ( c ) 
\not\succeq c). \end{cases}$$
Moveover, $\mathbb L^{<0} \mathscr D_i ( \mathbb W_c ) = 
\{ 0 \}$. 
\end{prop}

{\it Proof.}
By the $\mathfrak g$\~invariance, $\mathscr D_i ( 
\mathrm{Gr}^b L ) \cong \mathrm{Gr}^b L$ and $\mathbb L^{-1} 
\mathscr D_i ( \mathrm{Gr}^b L ) = \{ 0 \}$ for every $b \in 
P_+$. Then we use that by construction,  
$$\mathbb W_c \cong ( L ^{c} + M ) / M,$$
where $M \subset L^c$ is the sum of $L^{b}$ such that $b \prec c$ 
and $b \in P_+$. Applying $\mathscr D_i$ to the short exact 
sequence 
$$0 \to M \to ( L ^{c} + M ) \to \mathbb W_c \to 0,$$
and utilizing Theorem \ref{UDC}, we obtain the exact sequence 
$$0 \to \mathbb L^{-1} \mathscr D_i ( \mathbb W_c )\to 
M \to ( L ^{c'} + M ) 
\to \mathscr D_i ( \mathbb W_c ) \to 0,$$
where $c' = s_i ( c )$\,\, if\,\, $s_i ( c ) \succeq c$, 
and $c' = c$\,\, if\,\, 
$s_i ( c ) \prec c$. Therefore $\mathbb L^{-1} \mathscr 
D_i ( \mathbb W_c ) = \{ 
0 \}$ and either 
$\mathscr D_i ( \mathbb W_c ) \cong \mathbb W_c$ 
(for $s_i ( c ) \preceq c$) or $\mathbb W_{s_ic}$ 
(for $s_i ( c ) \succeq c$) 
as required. \sq

\begin{cor}\label{aGWEYL}
For $b \in P_-$, one has: $\mathscr D_{w_0} ( \mathbb 
D _b ) \cong \mathbb W_{w_0 ( b )}$. \sq 
\end{cor}

\begin{cor}\label{L-EXACT}
Let $b \in P$ and $0 \le i \le n$. If $s_i ( b ) \succ b$, 
then we have a short exact sequence of 
$\widetilde{\mathfrak b}^-$\~modules: 
$$0 \rightarrow \mathbb D _b \rightarrow \mathscr D_i 
( \mathbb D _b ) \rightarrow 
\mathbb D _{s_i ( b )} 
\rightarrow 0.
$$
If $s_i ( b ) \preceq b$, then 
$\mathscr D_i ( \mathbb D _{b} ) = \{ 0 \}.$
Moveover, $\mathbb L^{<0} \mathscr D_i ( 
\mathbb D _{b} ) = \{ 0 \}$ in each of these
two cases. 
\end{cor}

{\it Proof.} Using Corollary \ref{G-GR-PROP},
$\mathbb D _b \cong L^{b} / \sum_{u (\!( 0 )\!) 
\prec b} L ^{u}$. Let
 $S \equal\{ u \in \widehat{W} \mid u \not\le 
\pi_b, s_i u  \not\le \pi_b 
\}$ and $M \equal \sum_{u \in S} L^u$. It is 
straightforward to 
see that $\mathscr D_i ( M ) \cong M$. Also, if 
$u \not\le \pi_b$ and $s_i \pi_b < \pi_b$, then we have $s_i u 
\not\le s_i \pi_b$ and $u \not\le s_i \pi_b$ by \cite[Theorem 
2.2.2]{BB}. Therefore, we arrive at the following
short exact sequence: 
$$0 \rightarrow M \rightarrow ( L^{b} + M ) \rightarrow 
\mathbb D _b \rightarrow 0.$$
Its image under $\mathscr D_i$ here, 
is the following short exact sequence:
$$0 \rightarrow M \rightarrow ( L^{b'} + M ) \rightarrow 
\mathscr D_i ( \mathbb D _b ) \rightarrow 0.$$
See Proposition \ref{GWEYL}. Here 
$b' = b$ if $s_i ( b ) \prec b$ and 
$b'=s_i ( b )$ if $s_i ( b ) \succ b$. 
Combining these two exact sequence and using 
Theorem \ref{JOS}, we obtain the required. \sq

\subsubsection{\sf Orthogonality relations}\label{subsub:orth}

For a $\Z$-graded vector space 
$V = \medoplus_{\, i \in \Z} 
V_{i}$, we define 
$$\mathrm{end} ( V ) \equal \bigoplus_{j \in \Z} \prod_{i \in \Z} 
\mathrm{Hom} ( V_{i}, V_{i + j} ) 
\subset \mathrm{End} ( V ),$$
which is a ring. If $V$ has a $\widetilde{\mathfrak 
b}^{-}$\~module structure, then 
the $\mathsf{d}$-grading can be used here as 
the $\Z$-grading. Accordingly, we define
$$\mathrm{end}_{\widetilde{\mathfrak b}^{-}} ( V ) \subset 
\mathrm{end} ( V ),$$
the subring of all endomorphisms
commuting with the action of $\widetilde{\mathfrak b}^{-}$.

Each $\mathbb W_{b}$ is 
cyclically generated by a $\widetilde{\mathfrak h}$\~eigenvector. 
For instance, all elements of 
$\mathrm{end}_{\widetilde{\mathfrak b}^{-}} ( \mathbb W_{b} )$
commute with $\mathfrak g$ assuming that $\mathbb W_{b}$ 
is $\mathfrak g$\~invariant, which occurs for $b \in P_{+}$. 

For each $b \in P$, the (restricted) dual of the level one 
integrable representation $L ( \Lambda _{[b]} )^{\vee}$ has
a unique extremal weight vector $\mathsf v^*_{b}$ 
(up to a scalar) that 
is dual to $\mathsf v_b = \mathsf v_{\pi_b} \in 
L ( \Lambda _{[b]} )$. The 
$\widetilde{\mathfrak h}$-weight of $\mathsf v^*_b$ is 
$- \Lambda_{\pi_b}$.

We define the usual (thin) Demazure module for $b \in P$ to be
$$D_b \equal U ( \widetilde{\mathfrak b}^- ) 
\mathsf v^*_b \subset L ( \Lambda _{[b]} )^{\vee}.$$
These modules are referred to as Demazure modules in \cite{Kum}.

\begin{prop}[\cite{San, Ion1, CI}]\label{D-CHARS}
For any $b \in P$ and the reduced decomposition $b = 
s_{i_1} s_{i_2} \cdots s_{i_\ell} \pi$ in $\widehat{W}$ 
%\extra{Is it the decomposition of $b$ or $\pi_b$? --- 
%Changing $\pi_b$ with $\pi_b w$ ($w \in W$) gives the 
%same result. Hence, I think we can use both choices 
%if we want. Nevertheless, we cannot avoid to replace 
%during the proof of Theorem \ref{EXT-EL} unless we take 
%it in the maximal length representative, that is neither 
%$b$ nor $\pi_b$, I think.} 
with $\pi \in \Pi$ (the reduced decomposition
of any $bw$ for $w\in W$ can be taken here), we have an 
isomorphism of $\widetilde{\mathfrak b}^-$\~modules
$$D_b \cong \mathscr D_{i_1} \circ \cdots \circ 
\mathscr D_{i_{\ell}} 
( \C _{- \pi \Lambda_0} ).$$
Moreover, we have the followings: 
\begin{enumerate}
\item  $\mathsf{gch} \, D_b = q^{- 
    b^2/2} w_{0}\overline{E}_{b^{\iota}}$ 
(called below the character equality);
\item for $b \in P_+$, the module $\mathbb W_{b}$ admits a 
decreasing separable filtration as $\widetilde{\mathfrak 
    g}_{\le 0}'$\~modules whose associated graded is 
isomorphic to a direct sum of $D_{-b} \otimes 
    \C_{2 \Lambda_0}$ (i.e. $\mathbb W_{b}$ is a 
self-extension of $D_{-b} \otimes 
    \C_{2 \Lambda_0}$). In addition, 
    $\mathrm{end}_{\widetilde{\mathfrak b}^{-}} ( \mathbb 
    W_b )$ is a polynomial ring. 
\end{enumerate}
\end{prop}

%\extra{What is self-extension? --- modified explanation.}
%\extra{I don't understand: possibly infinitely repeated.
%Can you please exactly explain what do you mean here. --- I did.}

{\it Proof.}
The character equality is from \cite{Ion1}. The second claim 
follows from \cite[Corollary 2.10]{CI}.\sq

\begin{prop}[\cite{FKM} Appendix A]\label{EP-PAIRING}
For $M \in \mathfrak B$ and any $b \in P$,
\begin{align*}
\bigl< w_{0} ( \overline{E}_{b^{\iota}} &
\bigl(\bar{\mu}/\lan \bar{\mu} \ran \bigr) 
\mathsf{gch} \, M \bigr> \\
=&\,\,q^{- \frac{b^2}{2}} 
\sum_{p \ge 0,\, m, k \in \Q} (-1)^p q ^{- m} 
\dim \, \mathrm{Ext}^p _{( \widetilde{\mathfrak b}^-,
 \widetilde{\mathfrak h} )} ( M \otimes_{\C} 
\C _{m \delta + k \Lambda_0}, D_b^{\vee} ),
\end{align*}
where $\left< \cdot \right>$ denotes the constant term.
Recall that $D_b^{\vee}$ is the restricted dual (coinciding
with the full dual since  $\dim \, D_b < \infty$). %\sq 
\end{prop} 

%\extra{What is star in  $D_b^*$ here and in quite a few
%places incl. $\mathbb W$ below. Please explain and recall
%below once or twice. --- Maybe it is better to change uniformly 
%to $D_b^{\vee}$. ${}^*$ is the full dual, while ${}^{\vee}$ 
%is the restricted dual. They are same here as 
%$\dim \, D_b < \infty$.}

\comment{
\extra{You restricted $k$ and made
$m\in (1/e)\Z$ in the definition of  $\mathfrak B$, not all
$k,m$ can appear, am I right? Anyway
please check and comment. --- 
I think it is OK as all the weight spaces are zero unless 
it is relevant to what we have discussed. In any case, 
I have modified overall to improve little.}}

{\it Proof.} This claim is essentially
from \cite{FKM} for the $ADE$ systems. 
The same approach is applicable to 
arbitrary twisted $\tR$. Namely, we 
use (\ref{mubar}) and 
(\ref{constermbar}); they result in  $\mathsf{gch} \, \mathfrak 
n^- = w_{0} \mu_{\circ} ( t \to 0 )$. Then Proposition 
\ref{BPROJ} provides that the proof from 
\cite[Appendix A]{FKM} works in this generality. \sq 
 
\vskip 0.2cm
The following is a generalization of  
\cite[Theorem 5.12]{FKM}: from $\mathfrak g$ of types 
$\mathsf{A,D,E}_6,\mathsf{E}_7$ to arbitrary twisted ones.

\begin{thm}\label{U-EXT-ORTH}
For any $b, c \in P$ and $m, k \in \Q$, one has:
$$\dim \, \mathrm{Ext}^p _{( \widetilde{\mathfrak b}^-, 
\widetilde{\mathfrak h})} ( \mathbb D_b \otimes_{\C} 
\C_{m \delta + k \Lambda_0}, D_{c}^{\vee} ) = \delta_{p,0} 
\delta_{b,c} \delta_{m,0} \delta_{k,0}.
$$
\end{thm}

{\it Proof.}
As above, this claim is basically an extension of that from 
\cite[Theorem 5.12]{FKM}. 
We will not give a full proof; the following
modifications of \cite[Section 5]{FKM} are necessary.

First of all, the claim holds for $k=0$, since 
the action of $\widetilde{\mathfrak b}^-$ does not 
change the $K$-eigenvalue. Then Proposition \ref{GWEYL}
and Corollary \ref{L-EXACT}
are applied  to extend the proof of Theorem 
5.12 and Proposition 5.9 in \cite{FKM}. We also use 
\cite[Theorem 4.4]{CI}.
\sq 

%\extra{Non-understandable: what is a substitute of what;
%please fix! --- It is impossible to FIX (unless we reproduce 
%\cite[Section 5]{FKM} that would take at least few pages). 
%I have just added explanation of the nature of this ``proof".}

\begin{cor}\label{EP-ORTH}
For any $b, c \in P$, one has:
$$\left< ( \mathsf{gch} \, \mathbb D_b ) w_{0} ( 
\overline{E}_{c^{\iota}} \bar{\mu} / \lan \bar{\mu} \ran ) 
\right> = \delta_{b,c}.$$
\end{cor}

{\it Proof.}
Combine Proposition \ref{EP-PAIRING} and Theorem 
\ref{U-EXT-ORTH}.\sq

\begin{cor}[\sf Demazure slices and $E^{\dag}$-polynomials]
\label{E-D-DAG}
For any $b \in P$, 
$$\mathsf{gch} \, \mathbb D_b = q^{\frac{b^2}{2}} 
\frac{w_{0} (E^{\dag *}_{b^{\iota}})}{h_{b}^0}.$$
\end{cor}

{\it Proof.}
One has: $\mathsf{gch} \, \mathbb D_b \le 
\mathsf{gch} \, \mathbb W_{b_+}$
due to Theorem \ref{G-DEM-CH} and Corollary \ref{G-GR-PROP}.
The orthogonality relations from (\ref{ghgt0}) with the 
polynomials 
$\{ \overline{E}_b \}_{b \in P}$, provided
by Corollary \ref{EP-ORTH} for 
$\mathsf{gch} \, \mathbb D_b$,
and the comparison of the coefficients of the leading 
monomials give the required. 

We mention that all $\mathfrak h$\~weights of 
$\mathbb W_{b_+}$ belong to $\sigma _-( b )$ for
$\si_-(b)$ from (\ref{cones})
by \cite[Proposition 4.3]{CI}. 
Therefore  
$\mathsf{gch} \, \mathbb D_b \in \sum_{b_- \preceq c} 
\Q ( q ) X_c$. The inequality for $c$ here is weaker
than  $b\preceq c$ in the definition/construction of the
$E$\~polynomials, but using this fact
is not necessary anyway in our approach (the orthogonality
relations with $\overline{E}_b$ are sufficient).
\sq 

\comment{
Indeed, $h_{b} = h_{\iota ( b )}$ 
for the norms from (\ref{ghqt}) and 
$\{ E^{\dag}_b \}_{b \in P}$ is the only family of 
polynomials in $\Q (q) [X_a, a \in P]$ in terms of the
monomials $X_c$ with $c \succeq b_-$ satisfying the 
orthogonality relations from (\ref{ghgt0})
with the polynomials $\{ \overline{E}_b \}_{b \in P}$.
} 

%\extra{You need only the orthogonality relations, which
%we know. Is it correct? --- Yes. I also added ``, which we know 
%by Corollary \ref{E-D-DAG},"}

\begin{cor}\label{TR-PROJ-U}
For $b \in P$, the module $\mathbb D_b$ is projective 
in the category of  $\widetilde{\mathfrak h}$\~semisimple
$U ( \widetilde{\mathfrak b}^- )$\~modules 
with the weights $c$ satisfying $c \succeq b$. 
\end{cor}

{\it Proof.}
We use Theorem \ref{U-EXT-ORTH} and that $X_c$ occur
in $E_b$ only for $c\succeq b$. Thus
$
\mathrm{Ext}^{>0} _{( \widetilde{\mathfrak b}^-, 
\widetilde{\mathfrak h})} 
( \mathbb D_b \otimes_{\C} \C_{m \delta + 
k \Lambda_0}, \C _{b'} )
 = \{0\}$
for every $b \preceq b'$ and $m, k \in \Q$; here
the vanishing property  holds of course
for any $m, k \in \C$. \sq
 
\setcounter{equation}{0} 
\section{\sc Filtrations of tensor products}\label{sec:exp-theta}
In this section, the identification of the characters
of Demazure slices with $E$-dag polynomials will be used
to address Theorem \ref{MAINTHM}. We employ Theorem \ref{LZD},
which is related to similar results in \cite{Kas05,Kat16} 
in the $ADE$ case (it establishes a connection with the 
intertwiners from \cite{ChO}). The vanishing theorem 
from \cite{KL} and, its variant,
Corollary \ref{TEXT-VAN} are also important to us;
see also Theorems \ref{EXT-EL} and \ref{EXT-UL}. These facts
are essentially sufficient to interpret the existence
of the expansion from (\ref{pggmix}) representation-theoretically.
At the end, we discuss the remaining theta-function 
expansions from Theorem \ref{MAINTHM}. The exposition
is compressed in this section, with quite a few 
references to prior works and some details omitted (especially
if the corresponding results are known in the $ADE$ case).
Recall that we provided the list of main modules
under consideration in the beginning of
Section \ref{sec:Demazure}.

\subsection{\bf The 
\texorpdfstring{\mathversion{bold}$\mathbb W_c$}{\em W}-modules}
\subsubsection{\sf Demazure operators}\label{subsub:demop}
Recall that we have the {\em Demazure operators\,}:
$$
T_i ^{\dag} : \C (\!(q^{1/e})\!) [X_b] \ni f \mapsto 
\frac{f - X_{\alpha_{i}}s_{i} ( f )}{1 - X_{\alpha_{i}}},
\where X_{\alpha_0} = q^{-1} X_{- \vartheta},
$$
for each $0 \le i \le n$ defined in (\ref{tdagprime}), 
where $T_i$ is from (\ref{Demazx}). 

\begin{lem}
We assume that $M \in \mathfrak B$ has only finitely many 
distinct $\mathfrak h$\~weights and the central element $K \in 
\widetilde{\mathfrak h} \subset \widetilde{\mathfrak g}$ acts 
trivially (is zero) in $M$. Then
$\mathsf{gch} \, \mathscr D_{i} ( M ) = T_i ^{\dag} 
( \mathsf{gch} \, M )$
for each $0 \le i \le n$. 
\end{lem}

{\it Proof.} We can follow here \cite[5.4 Lemma]{J},
using that our $\mathfrak h$\~weights are 
bounded by the assumption.
 \sq

\subsubsection{\sf Level-zero theory}
Here, the action of $K$ is assumed to be trivial. The 
Lie algebra $\widetilde{\mathfrak g}_{\le 0}$ has 
one-dimensional representations $\C_{k \Lambda_0}$ and 
$\C_{m \delta}$ for $k\in \Z, m \in (1/e) \Z$ for $e$
as above: $e(P,P)/2=\Z.$   
Tensoring with such one-dimensional modules will be 
called {\em character twist}. 
The trivial (zero) $K$-action can be always achieved by 
an appropriate character twist.

The following is an extension of the 
corresponding result from \cite{Kat16}.

\begin{thm}\label{LZD}
For each $b \in P$ and $0\le i\le n$, we have
$$\mathsf{gch} \, \mathbb W_{s'_i (b)} = \begin{cases} 
q ^{- \delta_{i0} \vartheta ^{\vee} ( b )}
T_i ^{\dag} ( \mathsf{gch} \, \mathbb W_b ) & 
(\alpha_i ^{*} ( b ) \le 0),\\ \mathsf{gch} \, \mathbb W_b &
 ( \alpha_i ^{*} ( b ) > 0) \end{cases},$$
where $s_i' = s_i$ $(i \neq 0)$ and 
$s_0' = s_{\vartheta}$ $(i = 0)$. 
Also, $\mathbb L^{<0} \mathscr D_i ( 
\mathbb W_b ) =\{0\}$ for every $0 \le i \le n$. 
\end{thm}

{\it Proof.}
The case $i \neq 0$ is Proposition \ref{GWEYL}. 
Hence, assume $i = 0$; we will omit some details in the
proof below. 

When $\alpha_0 ^{*} ( b ) \ge 0$, it suffices to show that 
$\mathbb W_b$ is an integrable 
$\mathfrak p_0^-$\~module; see  
Theorem \ref{JOS}. To establish this, we will check that  
the $\widetilde{\mathfrak b}^-$\~action on 
$\mathbb W_b$ can be enhanced to a $\mathfrak p_0^-$\~action.
%\extra{ IS IT WHAT YOU ARE DOING BELOW? --- 
%Yes. This is a short summary of the 
%$\alpha_0 ^{*} ( b ) \ge 0$-case of our proof.}

For every $b \in P$, there is an embedding 
$\mathbb W_b \subset \mathbb W_{b_+}$ by definition
of $\mathbb W_b$.
%accordingly  $\mathsf v_b$ and $\mathsf v_{b_+}$ are 
%$\widetilde{\mathfrak b}^-$\~cyclic 
%vectors in this modules. 
The set of $\mathfrak h$\~weights of $\mathbb W_b$ 
is contained in the convex hull of $W(b) = W(b_+)$; 
see Theorem \ref{G-DEM-CH}. In particular, 
$\alpha_0 ^{*} ( b ) \ge 0$ implies that
\begin{equation}
\mathfrak g_{- \alpha_0}^{\alpha_0 ^{*} ( b ) + 1} 
\mathsf{v}_b = \{0\}, \hbox{\, for the
$\widetilde{\mathfrak b}^-$\~cyclic vector\, \,} \mathsf v_b\in 
\mathbb W_b.
\label{v-vanish-0}
\end{equation}

By Proposition \ref{DEF-W} and Theorem \ref{G-DEM-CH}, 
we obtain that $\mathsf{v}_b$ satisfies 
the relations from Proposition \ref{DEF-W} 
upon the application of $(w_0 u_b)^{-1}$ to the roots there.
Therefore $\alpha_0 ^{*} ( b ) \ge 0$ gives that
\begin{equation}
\mathfrak g_{- \alpha_0}^{\alpha_0 ^{*} ( b )} 
\mathsf{v}_b \neq \{0\}.\label{v-nonvanish-0}
\end{equation}

Next, $\mathsf{v}_{s_{\vartheta} b}$ 
has the same $\mathfrak h$\~weight as that of
$\mathfrak g_{- \alpha_0}^{\alpha_0 ^{*} ( b )} \mathsf{v}_b$,
but its $\mathsf{d}$-degree is different by  
$\vartheta ^{\vee} ( b )$. 
%\extra{: I EDITED A BIT. PLEASE DOUBLE-CHECK. --- 
%The expression is OK}
Thus Theorem \ref{free}
below implies that the isomorphism 
$\C \mathsf{v}_{s_{\vartheta}(b)} \cong 
\mathfrak g_{- \alpha_0}^{\alpha_i ^{*} ( b )} 
\mathsf{v}_b$ induces the following embedding of
$\widetilde{\mathfrak b}^-$\~modules: 
$$\mathbb W_{s_{\vartheta}(b)} \hookrightarrow 
\mathbb W_b.$$
%enhancing that of $\mathbb W_{b_+} \hookrightarrow 
%\mathbb W_{b_+}$. 
Moreover, we can equip
$\sum_{k = 0}^{\alpha_0 ^{*} ( b )} 
\mathfrak g_{- \alpha_0}^{k} \mathsf{v}_b$ with 
a structure of $\mathfrak{sl} ( 2 )_0$-module; see  
(\ref{v-nonvanish-0}), and (\ref{v-vanish-0}). 
Since $\mathsf{v}_b$ is a cyclic vector of $\mathbb W_b$, 
the PBW theorem implies that $\mathbb W_b$ admits an 
integrable $\mathfrak p_0^-$\~action. In particular, 
$$\mathscr D_0 ( \mathbb W_b ) \cong 
\mathbb W_b, \hskip 3mm \mathbb L^{<0} 
\mathscr D_0 ( \mathbb W_b ) \cong \{ 0 \}$$
whenever $\alpha_0 ^{*} ( b ) \ge 0$ by Theorem \ref{JOS}. 
Therefore, the case $\alpha_0 ^{*} ( b ) \ge 0$ of the 
assertion is proved. 
%\extra{ DID YOU PROVED THAT 
%$\mathbb W_b$ generates a $\mathfrak{sl} ( 2 )_0$-invariant 
%subspace of $\mathbb W_b$, which was stated as the claim you
%need to verify this case? PLEASE FIX THE LOGIC! --- 
%I just formally fixed the claim in the heading. 
%The previous assertion was equivalent to the current one.}

\vskip 0.2cm

Let us now consider the case 
$\alpha_0 ^{*} ( b ) \le 0$.  
%We will use 
%\cite[Lemma 4.4]{Kat16} (details will be omitted).
%It is needed to 
%to obtain the lower estimate for the character from (\ref{SUM}), 
%and to verify the existence of the filtration  
%filtration resulting  
%$\mathbb L^{<0} \mathscr D_0 = \{ 0 \}$. With the 
%help of numerical verification of an equality (\ref{ESUM}) 
%\extra{whose proof is omitted as my proof requires a 
%lot of preparation}, we deduce the desired equalities.
%In case $b \in P$ satisfies $\alpha_0 ^{*} ( b ) \le 0$, 
%then we swap the roles of $b$ and $s_{\vartheta} b$ 
%to obtain an embedding
We have now: 
$$\mathbb W_{b} \hookrightarrow 
\mathbb W_{s_{\vartheta}(b)},$$
where the latter space is an integrable 
$\mathfrak p_0^-$\~module with finitely many distinct 
$\mathfrak h$\~weights. This induces the following
homomorphism of $\widetilde{\mathfrak b}^-$\~modules:
%\extra{: I CHANGED THIS. CHECK PLEASE. :--- I have checked.}
$$\mathscr D_0 ( \mathbb W_b ) 
\stackrel{\eta}{\longrightarrow}  \mathbb 
W_{s_{\vartheta}(b)}, \hbox{\, where\, } 
\mathbb W_{s_{\vartheta}(b)}\ni
\mathsf{v}_{s_{\vartheta}(b)}\in 
\eta(\mathscr D_0 ( \mathbb W_b ) );$$
use (\ref{v-vanish-0}).  Therefore, $\eta$ is surjective
and 
\begin{equation}
T^{\dag}_{0} ( \mathsf{gch} \, \mathbb W_b ) \ge 
q ^{- \vartheta ^{\vee} ( b )} \mathsf{gch} \, 
\mathbb W_{s_{\vartheta} ( b )}\label{SUM}
\end{equation}
by \cite[Lemma 4.4]{Kat16}.

Employing now Lemma \ref{E-U-STR} and Corollary \ref{GR-PROP}, 
we obtain that $\mathbb W_c$ ($c \in P$) has a filtration 
by $\mathbb D_{c'}$ (as always, for appropriate 
character twists) for  $c' \in W(c)$ such that $c' \le c$; 
these modules occur (as constituents) with
multiplicities one. By Corollary \ref{E-D-DAG}, 
the inequality (\ref{SUM}) is equivalent to
\begin{equation}
\sum_{c \in W b} \varsigma_c ( b ) T^{\dag}_{0} 
( \frac{w_0 E_{c^{\iota}}^{\dag *}}{h^0_c} ) \ge 
\sum_{c \in W b}  q ^{- \vartheta ^{\vee} ( b )} 
\varsigma_{c}
 ( s_{\vartheta} ( b ) ) 
\frac{w_0 E_{c^{\iota}}^{\dag *}}{h^0_c};\label{ESUM}
\end{equation}
see Theorem \ref{RR-THETA}.
Thus (\ref{ESUM}) is in fact an equality, and therefore
(\ref{SUM}) is an equality too.  We obtain that
$\mathscr D_0 ( \mathbb W_b ) 
\cong \mathbb W_{s_{\vartheta}(b)}$ up to a character twist 
and $\mathbb L^{< 0}\mathscr D_0 ( \mathbb W_b ) = \{ 0 \}$ 
in the case $\alpha_0 ^{*} ( b ) \le 0$. Here we use again
\cite[Lemma 4.4]{Kat16}. \sq
 
%\extra{The passage to any $c$ from $c_+$ in the last par.
%is not convincing to me; I think the referee will not 
%accept this. Maybe I miss something. }

%\extra{THE PROOF ABOVE IS NOT UNDERSTANDABLE TO ME. IF 
%IT IS DIFFICULT TO FIX, PLEASE CONSIDER
%OMITTING IT COMPLETELY. IN THIS CASE, PLEASE WRITE 1-2 PAR
%ON THE TOOLS YOU USE.} 

\subsubsection{\sf The ring
\texorpdfstring{\mathversion{normal}$\mathrm{end} 
(\mathbb W_c )$}{\em end(W) }}\label{subsub:end}

\begin{thm}[\cite{FL}, \cite{FMS},\cite{CI} 
Corollary 2.4]\label{free}
Let $c \in P_+$. The ring $\mathrm{end} ( \mathbb W_c )$ 
is isomorphic to a graded polynomial ring,
namely, to:
$$\bigotimes _{i=1}^n \C [X_{i,1}, \ldots, X_{i,m_i}]
^{\mathfrak S_{m_i}}, \hbox{\,\, where\, \,} 
 m_i\!= \!\al^{\vee}_i ( c_+ ),\, \ \deg X_{i, j}\!=\!1,
$$
\vskip -0.2cm
\noindent
$\mathfrak S_m$ are the symmetric groups. Moreover, 
$\mathrm{end} ( \mathbb W_c )$ is isomorphic 
to $\mathrm{Hom} _{\mathfrak h} ( \C_{b}, \mathbb W_c )$ 
for every $b \in W(c)$.\sq
\end{thm}

Using the results from the previous section, we can extend
this theorem to any $c\in P$.

\begin{cor}\label{ENDW}
For any $c\in P$, the ring $\mathrm{end} 
( \mathbb W_c )$ does not depend on
$c \in W(b)$ and therefore it is isomorphic to the
graded polynomial ring from Theorem \ref{free}.
%\extra{ which polynomial ring exactly? --- 
%I added, but I am rather doubtful 
%if this helps understanding....}
\end{cor}

{\it Proof.}
The case $c = c_+$ is Theorem \ref{free}. This 
ring is isomorphic to 
$\mathrm{Hom}_{\mathfrak h} ( \C_{c}, \mathbb 
W_c ) \subset \mathbb W_c$ for such $c$.  
%\extra{what is $c$-part?--- I rewrote.}

Now let us consider an arbitrary $c\in P$. We will 
proceed by induction. Namely, we assume this claim 
for $c$ and prove it for $s_i ( c )$ for $0 \le i \le n$ 
provided that either 
the corresponding length increases for $i \neq 0$, 
or that the corresponding length decreases for $i = 0$.
%\extra{K. This is a tricky point as this order is 
%the semi-infinite/generic Bruhat order. It differs from 
%usual Bruhat order in the sense that the length can be negative. 
%Hence, I added a sentence before the end of the proof.}

%\extra{correct? --- Yes.}.
For any $0 \le i \le n$, $\mathscr D_i$ is a functor and 
therefore  induces the following homomorphism of algebras: 
\begin{equation}
\mathrm{end} ( \mathbb W_{c} ) \rightarrow \mathrm{end} 
( \mathscr D_i ( \mathbb W_{c} ) ) = \mathrm{end} 
( \mathbb W_{s_i ( c )} ),
\label{end-functorial}
\end{equation}
where $s_i ( c ) \gg c$\, for $i \neq 0$ and 
$s_{\vartheta} ( c ) \ll c$\, if\, $i = 0$; see
(\ref{order}).

Moreover,  
%$\mathrm{Id} \rightarrow \mathscr D _i$ 
Theorem \ref{JOS} gives a homomorphism
of $\widetilde{\mathfrak b}^-$\~modules 
%\extra{: I ADDED b-modules --- Thanks!} 
$\mathbb W_{c} \to \mathscr D _i ( \mathbb W _c )$, 
which is an inclusion. The latter follows directly
from the definition of 
$\mathbb W_c$ for $i \neq 0$ and results from  
Theorem \ref{LZD} when $i = 0$.

%\extra{ %The end of the proof must be rewritten. I don't
%%understand what you are doing. --- 
%%Substantial editing was done
%I STILL DON'T UNDERSTAND HOW YOU OBTAIN THE EMBEDDING
%$\mathbb W_c \hookrightarrow \mathbb W_c$. --- 
%I complimented here. In the below, I have simplified 
%the proof.}

For every $0 \le i \le n$ and $c \in W(b)$, the reduction to 
$\mathfrak{sl} ( 2 )_i$ gives that:
\begin{equation}
\mathbb W_c \supset \mathrm{Hom}_{\mathfrak h} 
( \C_{c}, \mathbb 
W_c ) \cong \mathrm{Hom}_{\mathfrak h} 
( \C_{s_i(c)}, \mathbb 
W_{s_i(c)} ) \subset \mathbb W_{s_i(c)},\label{wt-trans}
\end{equation}
with a possible character twist of $\mathbb W_{s_i(c)}$ when 
$i = 0$. These are maps of linear spaces.

Since $\mathbb W_c$ is cyclic, 
%$\widetilde{\mathfrak b}^-$\~module, $\mathrm{end} 
%( \mathbb W_c )$ is completely determined by its behavior 
%on the $\mathfrak h$-weight $c$-part of $\mathbb W_c$. 
%In particular, we have 
$\mathrm{end} ( \mathbb W_c ) \subset 
\mathrm{Hom}_{\mathfrak h} ( \C_{c}, \mathbb 
W_c )$ for every $c \in P$.
Let us assume that 
$$\mathrm{end} ( \mathbb W_c ) \cong 
\mathrm{Hom}_{\mathfrak h} ( \C_{c}, \mathbb 
W_c ) \and \mathscr D _i ( \mathbb W _c ) \cong 
\mathbb W _{s_i c}.
$$  Then (\ref{wt-trans}) gives that
$$\mathrm{end} ( \mathbb W_c ) \cong 
\mathrm{Hom}_{\mathfrak h} ( \C_{c}, \mathbb 
W_c ) \cong \mathrm{Hom}_{\mathfrak h} 
( \C_{s_i(c)}, \mathbb 
W_{s_i(c)} ) \supset \mathrm{end} ( \mathbb W_{s_i(c)} ).$$
Therefore the image of the map from (\ref{end-functorial}) 
is the whole eigenspace for the $\mathfrak h$\~weight 
$s_i ( c )$ in  $\mathbb W_{s_i(c)}$, and we conclude that
$$\mathrm{end} ( \mathbb W_c ) \cong 
\mathrm{Hom}_{\mathfrak h} ( \C_{s_i(c)}, \mathbb 
W_{s_i(c)} ) \cong \mathrm{end} ( \mathbb W_{s_i(c)} ).$$
 
%\extra{a character twist
%or some other twist of cyclic vectors? --- Rewritten.}
%the $\mathrm{end} ( \mathbb W_c )$-twist of cyclic vectors

The assumption we used here holds for $c\in P$ 
if it holds for $s_i ( c ) \gg c$ 
%\extra{WHICH ORDER DO 
%YOU USE HERE? $<$ means Bruhat. --- Dominance here.}
for some $i \neq 0$ or if it holds for  
$s_{\vartheta} ( c ) \ll c$. 
%Considering the reduced 
%decompositions of $b u_c \in \widehat{W}$ for $b \in P_- \cap Q$ 
%with $( \al_i ^{\vee}, b ) < 0$ for $1 \le i \le n$ (with 
%all $s_0$ replaced with $s_{\vartheta}$), we arrive at 
%$c_-$ from every $c \in P$ by a repeated 
%application of the above simple reflections. 
Hence, we can employ the induction, starting with $c=c_+$.\sq

%\extra{I STILL
%DON'T UNDERSTAND WHAT YOU ARE DOING. PLEASE MAKE AT LEAST THE
%LOGIC (BASIC STEPS) UNDERSTANDABLE OR OMIT THE PROOF. 
%--- I complimented in the above.} 

\vskip 0.2cm
\subsubsection{\sf Introducing {\it W}-modules}\label{SEC:W-mod}
Let 
$W_c \equal \left( \C_0 \otimes_{\mathrm{end} 
( \mathbb W_c )} 
\mathbb W_c \right) \otimes \C _{- \Lambda_0}$ for
$c\in P$; 
the action of $K$ is trivial (zero) in this module.

\begin{cor}\label{rLZD}
Similar to Theorem \ref{LZD}, 
$$\mathsf{gch} \, W_{s_i(b)} = \begin{cases} 
q ^{- \delta_{i0} \vartheta ^{\vee} ( b )}T^{\dag}_{i} 
( \mathsf{gch} \, W_b ) & 
(\alpha_i ^{*} ( b ) \le 0)\\ \mathsf{gch} \, W_b & 
(\alpha_i ^{*} ( b ) > 0), \end{cases}$$
for  any  $0 \le i \le n\,$ and $\,b \in P$. Moreover,  
$\, \mathbb L^{<0} \mathscr D_i ( W_b ) = \{ 0 \}$. 
\end{cor}

{\it Proof.}
We use that $\mathbb W_{s_i(b)}$ and $\mathbb W_b$ are free 
modules over the polynomial ring 
$\mathrm{end} ( \mathbb W_{b_+} )$.
The {\em Koszul resolution\,} of $\C$ considered as an
$\mathrm{end} ( \mathbb W_{b_+} )$\~module therefore 
results in 
$\widetilde{\mathfrak g}_{\le 0}$\~module resolutions of 
$W_{s_i(b)}$ 
and $W_b$ by some complexes whose terms are direct sums of 
$\mathbb W_{s_i(b)}$ and $\mathbb W_b$, respectively. 
%\extra{What is resolution BY something? --- 
%I added little more words.}
Then we apply $\mathscr D_i$ to these
resolutions and deduce the claim 
from Theorem \ref{LZD}. \sq

\begin{proposition}\label{coin}
For any $b \in P_+$, we have an isomorphism 
$W_b \cong D_{- b}$ 
as $\widetilde{\mathfrak g}_{\le 0}'$\~modules,
which may require a character twist.
\end{proposition}

{\it Proof.}
This follows
from  \cite[Theorem 2.7]{CI}. See also Theorem \ref{D-CHARS}.
%\extra{Is this Lemma or Corollary; the prove must be
%rewritten. I suggest: The claim follows
%from  \cite[Theorem 2.7]{CI}. See also Theorem \ref{D-CHARS}. --- 
%Here I just followed. As we lack ``Theorem" from \cite{CI}, 
%I think 
%we can keep ``Lemma".}
\vskip -0.5cm\sq

\subsection{\bf Vanishing theorems}
\subsubsection{\sf General results}\label{sec:General}
Let $L ( b + k \Lambda_0 )$ be the {\em integrable highest 
weight $\widetilde{\mathfrak g}$\~module\,} 
for the highest weight 
$b + k \Lambda _0 \in \widetilde{\mathfrak h}^*$. Here  
$b \in P_+$,
$k\in \Z_{\ge 0}$
and we assume that $\alpha_0 ^{*} ( b + k \Lambda _0 ) 
\in \Z_{\ge 0}$, equivalently,  
$( \vartheta ^{\vee}, b ) \in \Z_{\le k}\,;$
otherwise it does not exist. 
%It exists if and only if $b \in P_+$ and 
%$$\alpha_0 ^{*} ( b + k \Lambda _0 ) \in \Z_{\ge 0} 
%\Leftrightarrow 
%( \vartheta ^{\vee}, b ) \in \Z_{\le k}.$$

\begin{thm}[\cite{KL}]\label{KL-FILT}
For any $b,b' \in P_+$, $k \in \Z_{> 0}$, and $k' \in \Q$ 
such that $( \vartheta ^{\vee}, b ) \le k$, 
%\extra{ there is only pairing (,) --- Corrected.} 
we have 
$$\mathrm{Ext}^p _{( \widetilde{\mathfrak g}_{\le 0}, 
\mathfrak g + 
\widetilde{\mathfrak h} )} ( L ( b + k \Lambda_0 ), 
W_{b'}^{\vee} \otimes_{\C} 
\C _{k' \Lambda_0 + m \delta} ) = \{ 0 \} \for  p > 0,$$
where $\mathrm{Ext}^p _{( \widetilde{\mathfrak g}_{\le 0}, 
\mathfrak g + \widetilde{\mathfrak h} )} 
( \,\cdot\,, \,\cdot\, )$ 
are defined for the relative Lie algebra 
cohomology; see \cite[I\!I\!I]{Kum}. 
%\extra{what is relative extension with $p$, higher ext? --- 
%This is not the first time we use this. In any case, this is 
%a version of the relative Lie algebra cohomology/homology. For 
%this particular case, it can be also represented as a part of 
%the Hochschild-Serre spectral sequence. This is also higher 
%extensions of the category of $( \mathfrak g + 
%\widetilde{\mathfrak h})$\~semisimple 
%$\widetilde{\mathfrak g}_{\le 0}$\~modules. But I just recommend 
%readers to open the book of Shrawan or Borel-Wallach.}
\end{thm}
{\it Proof.}
The corresponding claim from \cite{KL} can be
extended to the twisted case following Section 
\ref{AWG}. We omit details.\sq

\begin{cor}\label{TEXT-VAN}
Let $b,b' \in P_+$, $k \in \N$\, and $\,m,\,k' \in \Q$ 
provided  the inequality 
%\extra{pairing! --- Thanks! Corrected.} 
$( \vartheta ^{\vee}, b ) \le k$. Then 
$$\mathrm{Ext}^p _{( \widetilde{\mathfrak b}^-, 
\widetilde{\mathfrak h} )} 
( L ( b + k \Lambda_0 ), D _{- b'}^{\vee} \otimes 
\C_{k' \Lambda_0 + m \delta} ) 
= \{ 0 \} \hbox{\, for any\, } p > 0.$$
\end{cor}

{\it Proof.}
We know that $W _{c} \cong D_{-c}$ for $c \in P_+$ by 
Lemma \ref{coin}. 
Thus it suffices to see that the higher $Ext$ in the category
of $( \widetilde{\mathfrak g}_{\le 0}, \mathfrak g + 
\widetilde{\mathfrak h} )$\~modules are the same as those for
the action of 
$( \widetilde{\mathfrak b}^-, \widetilde{\mathfrak h} )$. 
This can be seen from the Hochschild-Serre spectral sequence 
(\cite[E.12]{Kum}).
Indeed, 
$$\mathrm{Hom} _{( \mathfrak g, \mathfrak h)} ( M, N ) = 
\bigoplus_{p \in \Z} \mathrm{Ext}^p _{( \mathfrak g, 
\mathfrak h)} ( M, N ) 
\cong \bigoplus_{p \in \Z} 
\mathrm{Ext}^p _{( \mathfrak b^-, \mathfrak h)} 
( M, N )$$
for every $\mathfrak g$\~modules $M$ and $N$ that are direct 
sums of finite-dimensional $\mathfrak g$\~modules.
This finiteness property holds, in particular,
for the minimal projective resolution of 
$L ( b + k \Lambda_0 )$ in the category of 
$\widetilde{\mathfrak h}$\~semisimple 
$\widetilde{\mathfrak g}_{\le 0}$\~modules. \sq

%\extra{ Please rewrite the rest. --- 
%I minimized the explanation here.}

\subsubsection{\sf Filtrations by
\texorpdfstring{\mathversion{bold}$\mathbb W,
\mathbb D$}{\em W,D}-modules }

\begin{thm}\label{CRIT-FILT}
Let $M$ be a $( \widetilde{\mathfrak b}^-, 
\widetilde{\mathfrak h} )$\~module, whose set of 
$\widetilde{\mathfrak h}$\~weights is contained in 
$\bigcup_{j = 1}^m \left( \Lambda_j + \sum_{i = 0}^n 
\Z_{\le 0} \alpha_i\right)$ for a finite subset 
$\{\Lambda_j\}_{j = 1}^m \subset \widetilde{\mathfrak h}^*$. 

(i) A module $M$ admits a decreasing separable filtration such that 
    its associated graded components are of the form
 $\{ \mathbb W_b 
    \}_{b \in P_+}$, possibly with character twists, 
%\extra{what is
%it? --- I made definition just before Theorem \ref{LZD}} 
%\extra{what do you mean by the
%latter? --- Subsumed to the character twist}
 if and only if it is obtained as the 
    restriction of some $( \widetilde{\mathfrak g}_{\le 0}, 
    \mathfrak g + \widetilde{\mathfrak h} )$\~module and, 
additionally, the following holds  for 
    any $c \in P_-$\, and\, $m, k \in \Q$\,:
$$\mathrm{Ext}^1 _{( \widetilde{\mathfrak g}_{\le 0}, 
\mathfrak g + \widetilde{\mathfrak h} )} ( M \otimes 
\C _{m \delta + k \Lambda_0}, D _{c}^{\vee} ) = \{ 0 \}.$$

(ii) Similarly, $M$ admits a decreasing separable filtration 
such that  its associated graded components are of the form 
$\{\mathbb D_b \}_{b \in P}$, possibly with character 
twists, %of $K$ %\extra{: explain --- the same as above},
if and only if for 
    any $c \in P$ and $m, k \in \Q$: 
$$\mathrm{Ext}^1 _{(\widetilde{\mathfrak b}^-, 
\widetilde{\mathfrak h})} 
( M \otimes \C _{m \delta + k \Lambda_0}, 
D _{c}^{\vee} ) = \{ 0 \}.$$

For $(i)$ or $(ii)$, the $\mathrm{Ext}^p$\~spaces
vanish in all positive degrees $p$.
\end{thm}

{\it Proof.}
We will begin with the ``only if" part.
For (ii), it follows from Theorem \ref{U-EXT-ORTH} 
and a (repeated) usage of the long exact sequences.
To prove the ``only if" part in the case
of $(i)$, one needs an 
analogue of this theorem  for the modules 
$\mathbb W_{- c}$ and $D_{b}^{\vee}$ when $b, c \in P_-$.
%\extra{Using  $P_+$, where  $P_-$ is used everywhere can
%be a bit confused; I understand that this is due to $-b$, but
%maybe you can comment on this. --- Maybe we should just correct. 
%I made the correction.}
Let us outline the proof in this case.
%\extra{is it what you are 
%doing below? --- Yes. As far as I look around, containing 
%the next paragraph.}

For $b \in P_-$, the module $D_b$ is $\mathfrak g$\~invariant. 
Hence, $\mathscr D_{i} ( D_b ) \cong D_b$ and $\mathbb 
L^{< 0} \mathscr D_{i} ( D_b ) = \{ 0 \}$ for each 
$1 \le i \le 
n$. Thus, $\mathscr D_{w_0} ( D_b ) \cong D_b$. 
Thanks to Theorem \ref{D-ADJ} and Corollary \ref{aGWEYL}:
$$\{ 0 \} = \mathrm{Ext} ^{>0} _{(\widetilde{\mathfrak b}^-, 
\widetilde{\mathfrak h})} ( \mathbb D_{c} \otimes 
\C _{m \delta}, D_{b}^{\vee}) \cong \mathrm{Ext} ^{>0} 
_{(\widetilde{\mathfrak b}^-, \widetilde{\mathfrak h})} 
( \mathbb W_{w_0(c)} \otimes \C _{m \delta}, D_{b}^{\vee})$$
for any $b, c \in P_-$. Then we use the long exact 
sequences, which gives the ``only if" part of $(i)$. 
%\extra{$b$ was arbitrary in
%(ii)?? not in $P_-$; it was in $P_+$ in (i). Please fix. --- 
%There are two rather canonical choices, and there are no typos. 
%I also edited so that the logic becomes more transparent.} 

Let us come to the ``if" part. We will consider only
$(ii)$ (the first case is similar).  
Let $c \in P\,$ be a maximal element with respect to 
$\preceq$ such that 
$$\bigoplus_{m,k \in \Q} \mathrm{Hom}
_{(\widetilde{\mathfrak b}^-, 
\widetilde{\mathfrak h})} ( M, D _{c}^{\vee} \otimes 
\C_{m \delta + k \Lambda_0} ) \neq \{ 0 \}.$$

Using Proposition \ref{D-CHARS} and Section \ref{defNSMac}, 
every $\mathfrak h$\~weight of $D_c^{\vee}$ satisfies 
$\succeq c$, and $c$ occurs with multiplicity one. 
Recall that for any $c'\in P$, there exists a unique  
simple $\widetilde{\mathfrak b}^-$\~submodule
of $D_{c'}^{\vee}$ and it is isomorphic to $\C_{c'}$.
%\extra{: ref or mini-definition --- changed the notation}
%\extra{ref --- 
%if this refers to the previous thing, it's from the very 
%definition of Demazure modules. Since $D_{c'}$ is cyclic, 
%it follows that $D_{c'}^{\vee}$ is cocyclic, and hence 
%there is unique simple submodule as 
%$\widetilde{\mathfrak b}^-$\~submodules. if this refers to the 
%following sentence, this is by our choice of $c$ in the previous 
%paragraph. I think it is too early to review the definition 
%of $c$.} 
Therefore the  maximality of $c$ results in 
\begin{equation}
\bigoplus_{m,k \in \C} \mathrm{Hom} _{(\widetilde{\mathfrak b}^-, 
\widetilde{\mathfrak h})} ( M, \mathrm{coker} 
( \C_{c} \to D_c^{\vee} ) \otimes 
\C _{m \delta + k \Lambda_0} ) = \{ 0 \}. \label{Hom-Cond}
\end{equation}
%from the maximality of $c$ (by a finite successive application 
%of short exact sequences of $\mathrm{Hom}$).

Let $M_1$ be the maximal quotient
(a $\widetilde{\mathfrak b}^-$\~module) 
of $M$ such that all 
its $\mathfrak h$\~weights $b$ satisfy $b \succeq c$. By our 
maximality assumption, every simple quotient of $M_1$ 
is isomorphic to $\C_c$.
Hence, Corollary \ref{TR-PROJ-U} results in a surjective  
$\widetilde{\mathfrak g}_{\le 0}'$\~homomorphism 
$\phi :(\mathbb D_c ) ^{\oplus r} \rightarrow M_1$ 
%\extra{this I don't understand: --- where??? 
%I have edited little around here.} 
for some $r \in \Z_{> 0} \cup \{ \infty \}$. 
We can further assume that 
every direct 
summand of $(\mathbb D_c )^{\oplus r}$ 
maps non-trivially to $M_1$, changing $r$ if necessary.

If a simple quotient of $\ker \, \phi$ contains $\C_{c'}$, 
%\extra{what is
%head? --- I changed the notation}, 
then $c' \succeq c$ and 
\begin{equation}
\mathrm{Ext}^1 _{(\widetilde{\mathfrak b}^-, 
\widetilde{\mathfrak h})} ( 
M \otimes\C _{m \delta + k \Lambda_0}, \C_{c'} ) 
\neq \{ 0 \} \hskip 5mm 
\text{for some} \hskip 2mm m,k \in \Q.\label{C-Ext1}
\end{equation}
Using that $\mathrm{Ext}^{>0}$ are vanishing and the long 
exact sequences associated with the short exact sequences for 
subquotients of $\{D_{b}^{\vee}\}_{b \succeq c}$:
$$\mathrm{Ext}^{p} _{(\widetilde{\mathfrak b}^-, 
\widetilde{\mathfrak h})} ( 
M \otimes\C _{m \delta + k \Lambda_0}, \C_{c'} ) 
= \{ 0 \} \,\,\,  
\text{for every} \hskip 2mm c' \succeq c, p > 0, m,k \in \Q.$$
We use here (\ref{Hom-Cond}). This contradicts (\ref{C-Ext1}). 
%\extra{why? --- I edited around to clarify the logic.} 
Therefore, 
$$\mathrm{Ext}^1 _{(\widetilde{\mathfrak b}^-, 
\widetilde{\mathfrak h})} ( M, \C _{c'} \otimes 
\C _{m \delta + k \Lambda_0} ) = \{ 0 \} \,\,
\text{for every} 
\hskip 2mm c' \succeq c, m,k \in \Q,$$
%\extra{ remove or explain the rest please: 
%--- Edited mainly above.}
and we obtain an isomorphism 
$( \mathbb D_c ) ^{\oplus r} 
\cong M_1$. % as (a simple quotient of) its kernel 
%yields (\ref{C-Ext1}). 
Then we replace $M$ with $\ker \, ( M \to M_1 )$ and 
proceed by induction; note that our condition on weights is 
stable under taking subquotients. This concludes 
the ``if part" for $(ii)$.
The last assertion follows from Theorem \ref{U-EXT-ORTH}
for both, $(i)$ and $(ii)$. 
\sq 
 
\vskip 0.2cm
Theorem \ref{CRIT-FILT} is actually a level-one version 
of the (non-affine) van der Kallen's criterion \cite{vdK}. 
The proof we suggest here is applicable only to level one, 
since we cannot generalize Corollary \ref{TR-PROJ-U} to 
any levels.

\subsubsection{\sf The passage to {\it W}-modules}
We will use the results above for $W$\~modules from
Section \ref{SEC:W-mod}.

\begin{thm}\label{EXT-EL}
Let $\Lambda$ be a level-one dominant integral weight. Then
for any 
$b,c \in P_+$, $m \in \Q$ and $p > 0$, 
\begin{equation}\label{extpw}
\mathrm{Ext}^p_{(\widetilde{\mathfrak b}^-, 
\widetilde{\mathfrak h})} 
( W_b \otimes L ( \Lambda ), W_c^{\vee} \otimes 
\C_{\Lambda_0 + m \delta} ) = 
\{ 0 \}.
\end{equation}
\end{thm}

{\it Proof.}
Since $W_b, D_b, L ( \Lambda ), W_c^{\vee}, D_b^{\vee}$ are 
indecomposable $\widetilde{\mathfrak b}^-$\~modules, the
$\mathfrak h$\~weights of each of them coincide in
$P/Q$. Therefore $\mathrm{Hom}$ and $\mathrm{Ext}$ 
vanish if the weights in the first and the second
$\mathrm{Ext}^p$\~component from (\ref{extpw}) have different 
images 
in $P/Q$. Thus the right-hand side there is $\{0\}$
when $b + c + 
\Lambda - \Lambda_0 \not\in Q + \Z \delta$. 

%\extra{there is no first assertion: --- Sorry. 
%This is a historical typo.}

For any $b,c\in P_+$ and $p \ge 0$ ($p=0$ is allowed),
%\extra{ do you allow $p=0$ here? --- Yes.},
one has: 
\begin{align*}\mathrm{Ext}^p_{(\widetilde{\mathfrak b}^-, 
\widetilde{\mathfrak h})} & 
( W_b \otimes L ( \Lambda ), W_c^{\vee} \otimes 
\C_{\Lambda_0 + m \delta} ) \\
\cong & \mathrm{Ext}^p_{(\widetilde{\mathfrak b}^-, 
\widetilde{\mathfrak h})} 
( W_b \otimes \C_{- \Lambda_0} \otimes L ( \Lambda ), 
W_c^{\vee} \otimes 
\C_{m \delta} ).
\end{align*}
Since $L ( \Lambda )$ and $\C_{m \delta}$ are integrable 
$\widetilde{\mathfrak g}$\~modules, this gives that
\begin{align*}\mathrm{Ext}^p_{(\widetilde{\mathfrak b}^-, 
\widetilde{\mathfrak h})} & ( \mathscr D_i 
( W_b \otimes \C_{- \Lambda_0} ) 
\otimes L ( \Lambda ), W_c^{\vee} \otimes \C_{m \delta} ) \\
\cong & \mathrm{Ext}^p_{(\widetilde{\mathfrak b}^-, 
\widetilde{\mathfrak h})} ( W_b \otimes \C_{- \Lambda_0} \otimes 
L ( \Lambda ), ( \mathscr D_i ( W_c ) )^{\vee} 
\otimes \C_{m \delta} )
\end{align*}
for any $0 \le i \le n$ and $p \ge 0$. 
%\extra{$\Z_+$ I think? --- I fixed.}
We used Theorem \ref{JOS} and Theorem 
\ref{D-ADJ}.
\vskip 0.2cm

Let $\Lambda' \equal  \Lambda_{[-b]}$ and $v^*$ be the 
lowest weight vector of $L ( \Lambda' )^{\vee}$. Then we 
have an embedding $W_b \otimes \C_{- \Lambda_0} \subset L ( 
\Lambda' )^{\vee}$ of $\widetilde{\mathfrak g}_{\le 0}$\~modules, 
possibly
up to a character twist; see Lemma \ref{coin}. Let us use 
Proposition \ref{D-CHARS}. We need a reduced decomposition 
$\pi_{w_{0}(b)} = s_{i_1} s_{i_2} \cdots s_{i_\ell} \pi$, where 
$i_1,\ldots,i_\ell \in [0,n],\, \pi \in \Pi$. Then :
\begin{equation}
W_b \otimes \C_{- \Lambda_0} \cong \mathscr D_{i_1} \circ 
\mathscr D_{i_2} 
\circ \cdots \circ \mathscr D_{i_\ell} ( \C v^* ).\label{SID}
\end{equation}
Since $\mathscr D_i ( \C v^* ) \cong \C v^*$ when 
$\alpha_i^{*} ( \Lambda' ) = 0$, we can replace $\pi_{w_{0}(b)}$ 
here with the maximal length representative $x\in \hW$ in
the double coset $W \pi_{w_{0}(b)} W$. The relation  
(\ref{SID}) will hold.
%\extra{:I changed $w_0 b$ here by $w_0(b)$ --- 
%Here it is correct.}

%\extra{ I suggest, remove this par:
%Let us rearrange $\ell$ and 
%the sequence $\{i_{j}\}_{j} \in [0,n]^\ell$ if necessary such that 
%it gives a reduced expression of $x$ with keeping (\ref{SID}) 
%hold. --- I don't understand. $x$ can never be equal to 
%$\pi_{w_0 b}$. 
%I need to change within the double coset so that the argument in 
%the below works. Nevertheless, I modified the sentence.}

Obviously $x^{-1}$ is also a maximal length 
representative of the corresponding
element in $W \backslash \widehat{W} / W$. Thus 
Corollary \ref{rLZD} implies that there exists 
$m'\in \Z$ such that 
\begin{align}\label{Dwprime}
\mathscr D_{\pi x^{-1}} ( W_c ) = \mathscr D_{i_\ell} \circ 
\mathscr D_{i_{\ell -1}} 
\circ \cdots \circ 
\mathscr D_{i_1} ( W_c ) \cong W _c \otimes \C_{m' \delta},
\end{align}
$\mathbb L^{<0}\mathscr D_{\pi x^{-1}} 
( W_c )$ are all  $\{ 0 \}$.
%\extra{: I don't understand "homology terms". --- I spelled out.} 
Then Theorem \ref{D-ADJ} implies that
\begin{align*}
\mathrm{Ext}^p_{(\widetilde{\mathfrak b}^-, 
\widetilde{\mathfrak h})} & 
( W_b \otimes \C_{- \Lambda_0}
\otimes L ( \Lambda ), W_c^{\vee} 
\otimes \C_{m \delta} ) \\
\cong & \mathrm{Ext}^p_{(\widetilde{\mathfrak b}^-, 
\widetilde{\mathfrak h})} 
( \mathscr D_{x} ( \C_{- \Lambda_0} ) 
\otimes L ( \Lambda ), W_c^{\vee} 
\otimes \C_{m \delta} ) \\
\cong & \mathrm{Ext}^p_{(\widetilde{\mathfrak b}^-, 
\widetilde{\mathfrak h})} 
( \C_{- \Lambda'} \otimes L ( \Lambda ),  W_c^{\vee} \otimes 
\C_{( m - m') \delta} ).
\end{align*}
Recall that $b,c \in P_+$,\, $m \in \Q$ are arbitrary and 
$m'$ is determined by (\ref{Dwprime}). 

Let us use now that
 $\pi\in \Pi$ above induces an 
automorphism of $\widetilde{\mathfrak g}$  preserving 
$\widetilde{\mathfrak b}^-$ and $\widetilde{\mathfrak h}$ and 
that $\pi \Lambda_0 = \Lambda_{[-b]}$. 
There exists $m''\in \Q$ and some $d\in P$ such that 
\begin{align*}
\mathrm{Ext}^p_{(\widetilde{\mathfrak b}^-,
\widetilde{\mathfrak h})}& 
( \C_{- \Lambda'} \otimes L ( \Lambda ),  W_c^{\vee} \otimes 
\C_{( m - m') \delta} )\\ 
&\cong \mathrm{Ext}^p
_{(\widetilde{\mathfrak b}^-, \widetilde{\mathfrak h})} 
( \C_{- \Lambda_0} \otimes L ( \Lambda'' ),  D_d^{\vee} \otimes 
\C_{- \Lambda_0 + m'' \delta} )
\end{align*}
for the standard (thin) Demazure module $D_d^{\vee} \subset 
L ( \pi^{-1} \Lambda_{[-c]} )$. Here
$\Lambda'' \equal \pi^{-1}\Lambda$ is a level-one 
dominant weight and there is an explicit (quadratic) formula
for $m''$ in terms of $\Lambda''$, initial (level-one
dominant) $\La$ and $\La'$ above.
We used the fact that the definition of
Demazure modules is stable under the diagram automorphisms. 

%\extra{EXPLAIN BETTER WHAT ARE $m''$, $d$. --- I don't think so. 
%In order to persuade you, let me spell out. 
%If we write $\Lambda = a + \Lambda_0, \Lambda'' = a'' + \Lambda_0$ 
%for some $a,a'' \in P_+$, then we have an equality
%$$\pi ( d + \Lambda_0 ) + 
%( m'' - \frac{b^2}{2} - \frac{(a'')^2}{2} ) 
%\delta = c + \Lambda_0 + ( m - m' - \frac{a^2}{2}) \delta.$$
%(By our convention on the degree of $W_c$, 
%there are no a priori cancelation.) Hence we have
%$$\pi ( d + \Lambda _0 ) \equiv c + \Lambda _0 \mod \delta, 
%\hskip 5mm m'' = \frac{(a'')^2 + b^2 - a^2}{2} + m - m'.$$
%It's very difficult for me to imagine that someone feel 
%happy with this. To save the situation, probably we need 
%some better notation system that I do not want to introduce...} 
Next, we use that $\mathscr D_i ( L ( 
\Lambda'' ) ) \cong L ( \Lambda'' )$ for any $0 \le i \le n$. 
One has:
\begin{align*}
\mathrm{Ext}^p_{(\widetilde{\mathfrak b}^-, 
\widetilde{\mathfrak h})} 
( L ( \Lambda'' ),  
D_d^{\vee} \otimes \C_{m'' \delta} ) &
= \mathrm{Ext}^p_{(\widetilde{\mathfrak b}^-, 
\widetilde{\mathfrak h})} 
( L ( \Lambda'' ),  D_d^{\vee} \otimes 
\C_{m'' \delta} )\\
& \cong \mathrm{Ext}^p_{(\widetilde{\mathfrak b}^-, 
\widetilde{\mathfrak h})} 
( \mathscr D_{w_0} ( L ( \Lambda'' ) ),  D_d^{\vee} \otimes 
\C_{m'' \delta} )\\
& \cong \mathrm{Ext}^p_{(\widetilde{\mathfrak b}^-, 
\widetilde{\mathfrak h})} 
( L ( \Lambda'' ),  \mathscr D_{w_0} ( D_d )^{\vee} \otimes 
\C_{m'' \delta} )\\
& = \mathrm{Ext}^p_{(\widetilde{\mathfrak b}^-, 
\widetilde{\mathfrak h})} 
( L ( \Lambda'' ), D_{d_-}^{\vee} \otimes 
\C_{m'' \delta} ).
\end{align*}
Finally, the usage of Corollary \ref{TEXT-VAN} 
completes the proof. Here, as with the proof of the
previous theorem, we omit some technical details. 
\sq

\subsection{\bf Theta-products 
via \texorpdfstring{$\mathbb D_b$}
{Demazure slices}}
\subsubsection{\sf Theta-products 
via \texorpdfstring{$\mathbb W_c$}{\it W}}
The following stratification is an application of the 
theory above.
\begin{cor}\label{MODULEMAIN-W}
Let $\mathbf v = ( \pi_1, \ldots, \pi_p )$ be a sequence of 
elements of $\Pi$. For each $b \in P_+$, the tensor product 
$$W_b \otimes L ( \pi_1 \Lambda_0 ) \otimes \cdots \otimes 
L ( \pi_p \Lambda_0 )$$
admits a decreasing separable filtration by $\{ \mathbb W_c 
\}_{c \in P_+}$ (as constituents), possibly with character twists. 
%\extra{: please explain --- I think this is OK here with 
%this modification.}
Moreover, the multiplicities in this 
filtration are given by formula (\ref{pggbar})
in Theorem \ref{MAINTHM} for $c = c_-$. 
\end{cor}

{\it Proof.} We first prove the existence of the 
filtration by induction on $p$. In the case $p = 1$, 
we apply Theorem \ref{EXT-EL}; the vanishing of  
$\mathrm{Ext}^1$ follows from Theorem 
\ref{CRIT-FILT}. Generally, the existence 
of such a filtration by 
$\{ \mathbb W_c \}_{c \in P_+}$ implies the existence
of the  filtration by modules
$\{ W_c \}_{c \in P_+}$. Therefore, 
we can go from $p = k - 1$ to $p = k$ by taking 
the associated graded. This gives the existence claim. 

The multiplicity claim for this filtration follows then 
from formula (\ref{pggbar}) for $c = c_-$ because the 
graded characters of $\{ \mathbb W_c \}_{c \in P_+}$ 
are linearly independent. \sq 

%\extra{PLEASE CLEARLY EXPLAIN HERE WHAT IS THE
%LOGIC OF THE PASSAGE FROM ABOVE COROLLARY TO THIS
%THM AND THEN TO THE NEXT COROLLARY, WHICH IS THE
%MAIN RESULT OF THIS SECTION. MAKE IT ONE PAR AT LEAST.
%REMIND WHAT IS WHAT, PLEASE. --- I did.} 

\subsubsection{\sf More on vanishing 
\texorpdfstring{$Ext^p$}
{\it Ext}} 
\begin{thm}\label{EXT-UL}
Let $\Lambda_1,\ldots,\Lambda_s$ be a sequence 
level-one dominant integral weight. We set
$$L ( \vec{\Lambda} )\equal L ( \Lambda_1 ) \otimes 
\cdots \otimes L ( \Lambda_s ).$$
For any
$\,b,c \in P$, $m \in \Q$ and $p > 0$, we have 
\begin{equation}
\mathrm{Ext}^p_{(\widetilde{\mathfrak b}^-, 
\widetilde{\mathfrak h})} 
( D_b \otimes L ( \vec{\Lambda} ), D_c^{\vee} 
\otimes \C_{(s-2)\Lambda_0 
+ m \delta} ) = 
\{ 0 \}.
\end{equation}
\end{thm}

{\it Proof.}
Let $\Lambda \equal \sum_{j = 1}^s \Lambda_j$ and 
$\Lambda = \omega + s \Lambda _0 \mod \Z\delta$ 
for proper $\omega \in P_+$. By Corollary 
\ref{MODULEMAIN-W}, we know that $D_{-b} \otimes L 
( \vec{\Lambda} )$ admits a decreasing separable filtration by $\{ 
\mathbb W_c \}_{c \in P_+}$, where $b \in P_+$. 
%\extra{This is corollary \ref{MODULEMAIN-W}
%and I removed minus in $b$. ---- We have $D_{-b} \cong W_b$ for 
%$b \in P_+$. The both of the modules $D_b$ and $D_b^{\vee}$ are 
%not $\mathfrak g$-stable when $b \in P_+$ in general.}
Using
Theorem \ref{RR-THETA} and formula
(\ref{thexfin0}), we conclude that 
$\mathbb W_c$ occurs in $D_{b} \otimes L ( \vec{\Lambda} )$
if and only if
$c \in - b + \omega + Q \mod \Z \delta$.  

Lemma \ref{E-U-STR} gives that each $\mathbb W_c$ admits a 
filtration by $\{\mathbb D_b\}$ with 
$b \in W(c)$ and that every such module occurs with 
multiplicity one. The character twists can be  needed
here; actually one (common) character twist serves  
all $b \in W (c)$. 
%\extra{: I don't 
%understand shift by multiplicity one --- I hope this works}

Now let us use that the
$\widetilde{\mathfrak b}^-$\~submodule 
$$D_{b_{+}} \otimes L ( \vec{\Lambda} ) \subset D_{b_{-}} 
\otimes L ( \vec{\Lambda} )$$
generates the whole $D_{b_-} \otimes L ( \vec{\Lambda} )$ as a 
$\mathfrak g$\~module, which is Corollary \ref{aGWEYL}. By the PBW 
theorem and the fact that $D_{b_{+}} \otimes L ( \vec{\Lambda} )$ 
is $\widetilde{\mathfrak b}^-$\~invariant, we obtain that 
$D_{b_-} \otimes L ( \vec{\Lambda} )$ generates $D_{b_{-}} 
\otimes L ( \vec{\Lambda} )$ as a $\mathfrak n$\~module.
\vskip 0.2cm

Let $b' \in P_+$ be such that 
$\mathbb W_{b'}$ appears in $D_{b_{-}} \otimes 
L ( \vec{\Lambda} )$ for the filtration of 
Corollary \ref{MODULEMAIN-W};
as always in this section, this is up to a suitable 
character twist.  
Then there exists
a pair of $\widetilde{\mathfrak g}_{\le 0}$\~submodules  
$M' \subset 
M \subset D_{b_{-}} \otimes L ( \vec{\Lambda} )$ such that 
$(a)$ \ $M'$ and $M$ 
admit compatible decreasing separable filtrations by 
$\{ \mathbb W_c \}_{c \in P_+}$,  
$(b)$\ $M / M' \cong \mathbb W_{c}$, 
and $(c)$\ $M$ is generated by $( D_{b_{+}} 
\otimes L ( \vec{\Lambda} )) \cap M $ 
as a $\mathfrak n$\~module. This results in the existence of
an irreducible 
$( \mathfrak g + \widetilde{\mathfrak h} )$\~submodule 
$V \subset M$ 
generating  $\mathbb W_{c} \cong M / M'$ as a 
$\widetilde{\mathfrak g}_{\le 0}$\~module, where every
$\mathfrak h$\~eigenvector of $V$ of the weight $c'$
from $W ( c )$ 
generates $\mathbb W_{c'}$ in $\mathbb W_{c} \cong M / M'$. 

%\extra{: Please check what I changed : --- I checked around 
%and added $(c \in P_+)$ in the first sentence.}
Generally, a $\widetilde{\mathfrak b}^-$\~invariant subspace 
$U \subset \mathbb W_{c}$ ($c \in P_+$) generates 
$\mathbb W_{c}$ 
with respect to the action $\mathfrak n$ if it  
contains the largest $\mathfrak h$\~submodule 
$\mathbb U_{c}$ of $\mathbb W_{c}$ that is a direct sum 
of some copies of $\C_{c_-}$. We use here (\ref{ncbcoef}).

By (the proof of) Corollary \ref{ENDW},
$\mathbb U_{c} \subset \mathbb W_{c_-}$ for such $\mathbb U_c$. 
The vector space $\mathbb U_{c}$ contains the cyclic vector of 
$\mathbb W_{c_-}$ and  generates the whole 
$\mathbb W_{c_-}$ considered as a 
$\widetilde{\mathfrak b}^-$\~module.
We conclude that every $\mathfrak h$\~submodule 
$\mathbb U_{c}^+ \subset M$ such that $( \mathbb U_{c}^+ + M' ) 
/ M' \cong \mathbb U_{c}$ generates $\mathbb W_{c_-} \cong 
\mathbb D_{c_-}$ upon the passage to
$\mathbb W_{c} \cong M / M'$.

Let us now use that $M$ is generated by $( D_{b_{+}} 
\otimes L ( \vec{\Lambda} ) \cap M )$ 
as a $\mathfrak n$\~module. We can choose 
$\mathbb U_{c}^+ \subset 
( D_{b_{+}} \otimes L ( \vec{\Lambda} ) )\cap M $. This implies 
that
$$\bigl(( D_{b_{+}} \otimes L ( \vec{\Lambda} )) \cap M \bigr) / 
\bigl(( D_{b_{+}} \otimes 
L ( \vec{\Lambda} )) \cap M' \bigr) \supset 
\mathbb D_{c_-}.
$$

As a consequence, if  $\mathbb W_{c}$ occurs in a filtration of 
$D_{b_{-}} \otimes L ( \vec{\Lambda} )$ (upon a suitable 
character twist) from Corollary 
\ref{MODULEMAIN-W}, then $\mathbb D_{c_-}$ occurs 
in $D_{b_{+}} \otimes L ( \vec{\Lambda} )$ for the induced 
filtration of $D_{b_{-}} \otimes L ( \vec{\Lambda} )$ (with 
the same character twist). 
This correspondence preserves the multiplicities.

For $s = 1$ case, the above count of multiplicities results in 
the inequality:

%\extra{ please split into 2-3 phrases and explain
%summand, subquotient of what: --- I completely rewrite 
%the whole paragraph.}

\begin{equation}
q^{b^2 / 2}\mathsf{gch} \, D_{b_{+}} \otimes L ( \vec{\Lambda} ) 
\ge \sum_{c \in ( - b + \omega + Q ) 
\cap P_{-}} q^{(b_- - c_-)^2 / 2} \cdot 
q^{-c^2 / 2}
\mathsf{gch} \, 
\mathbb D_{c}, \label{INEQ-DLU}
\end{equation}
where $q^{(b_- - c_-)^2 / 2}$ is
due to Theorem \ref{MODULEMAIN-W} and formula
(\ref{thexpanfin1}) for $c = c _-$.
 
The equality holds here if and only if $D_{b_+} 
\otimes L ( \Lambda )$ has a filtration by 
$\{ \mathbb D_c \}_{c \in P_-}$. 
However we know that (\ref{INEQ-DLU}) is actually 
an equality due to (\ref{thexfin0}), which gives the
required.

This argument remains essentially unchanged
for any  $s > 1$. The left-hand side
of (\ref{INEQ-DLU}) will then 
correspond to (\ref{pggbar}) with $c, b \in P_-$. 
See Theorem \ref{D-CHARS}, 
Corollary \ref{E-D-DAG}, and (\ref{pggmix}).

Thus  $D_{b_+} \otimes 
L ( \vec{\Lambda} )$ has a required filtration by 
$\{ \mathbb D_c \}_{c \in P_-}$ for every $s > 0$ 
and we can use Corollary \ref{L-EXACT} 
and Theorem \ref{CRIT-FILT} to establish that
$\mathscr D_i ( D_b \otimes L ( \vec{\Lambda} ) )$ has a 
filtration by $\{ \mathbb D_c \}_{c \in P}$.
This concludes the proof of the theorem. \sq 
\vskip 0.2cm

%\subsubsection{\sf Final decomposition}
The following application is one of the main results of
this half of the paper; it finally interprets formula
$(\ref{pggmix})$ representation-theoretically. 

\begin{thm}\label{MODULEMAIN-U}
Let $b \in P$. For each sequence 
$\mathbf v = ( \pi_1, \ldots, \pi_p )$  
of elements in $\Pi$, the tensor product 
$$D_b \otimes L ( \pi_1 \Lambda_0 ) \otimes \cdots \otimes 
L ( \pi_p \Lambda_0 )$$
admits a decreasing separable filtration by 
$\{ \mathbb D_c \}_{c \in P}$ with possible
character twists. 
The multiplicities in 
this filtration are given by Theorem \ref{MAINTHM}, namely by 
the $r = 0$ case of formula $(\ref{pggmix})$. 
\end{thm}

{\it Proof.}
In view of Theorem \ref{EXT-UL}, the existence of a filtration 
is a direct consequence of Theorem \ref{CRIT-FILT}. For its 
character, use (\ref{pggmix}) and Corollary 
\ref{E-D-DAG}. \sq

\subsubsection{\sf Some perspectives}

\begin{conj}\label{dual-orth}
For any $c \in P$, the module $\mathbb D_c$ is free 
over the following graded polynomial ring $R_c$, similar to that
in Theorem \ref{free}:
$$R_c\equal
\otimes _{i=1}^n \C [X_{i,1}, \ldots, X_{i,m_i}]
^{\mathfrak S_{m_i}}, \hbox{\,\, where\, \,} 
 m_i\!=  ( \al^{\vee}_i,  c )^{\malt},\, \ \deg X_{i, j}\!=\!1.
$$
See Section \ref{SEC:MAINTH} for the notation 
$(\cdot,\cdot)^{\malt}$.
%\extra{What is $R_c$? ---- I will update later. This is bit 
%difficult to write down.} 
Moreover, we conjecture the existence of a 
maximal (universal) 
cyclic $\widetilde{\mathfrak b}^-$-modules 
$\mathsf D_c$ such that
\begin{enumerate}
\item the module $\mathsf D_c$ maps surjectively 
onto $D_c$ 
(as $\widetilde{\mathfrak b}^-$-modules);
\item $\mathrm{end}_{\widetilde{\mathfrak b}^-} 
( \mathsf D_c ) \cong R_{-c}$ and
 $\mathsf D_c$ is free as an $R_{-c}$-module; 
\item and also, $h^0_c \cdot \mathsf{gch} \, \mathsf D_c = 
\mathsf{gch} \, D_c\for h^0_c$ from (\ref{limeval}).
\end{enumerate}
%\extra{what is
%$b$? --- Sorry. Typo}
%\extra{ : changed to  $\mathsf D$. 
%Also check next : --- I have checked.}
For such $\mathsf D_c$, we conjecture that
$\mathring{\mathbb D}_c \equal  \C_0 \otimes_{R_c} 
\mathbb D_c$ (local Demazure slices)
 satisfy the following orthogonality 
relations:
$$\dim \, \mathrm{Ext}^p _{(\widetilde{\mathfrak b}^-, 
\widetilde{\mathfrak h})} ( \mathsf D_b \otimes_{\C} 
\C_{m \delta}, 
\mathring{\mathbb D}_{c}^{\vee} ) = 
\delta_{p,0} \delta_{b,c} \delta_{m,0}.$$
\end{conj}
\vskip -.7cm \sq

\comment{
\begin{cor}\label{CRIT-FILT-D}
Let $M$ be a $( \widetilde{\mathfrak b}^-, 
\widetilde{\mathfrak h} )$\~module, whose set of 
$\widetilde{\mathfrak h}$\~weights is contained in 
$\bigcup_{j = 1}^m \left( \Lambda_j + \sum_{i = 0}^n 
\Z_{\le 0} \alpha_i\right)$ for a finite subset 
$\{\Lambda_j\}_{j = 1}^m \subset \widetilde{\mathfrak h}^*$.

Then $M$ admits a decreasing 
separable filtration whose associated graded are direct 
sums of $\{\mathsf D_b \otimes_{\C} \C_{k \Lambda_0 + 
m \delta}\}_{b \in P, k, m \in \Q}$
if and only if
$$\mathrm{Ext}^p _{(\widetilde{\mathfrak b}^-, 
\widetilde{\mathfrak h})} ( \mathsf D_b \otimes_{\C} 
\C_{m \delta}, 
\mathring{\mathbb D}_{c}^{\vee} ) = \{ 0 \} \hskip 5mm b \in P, 
k, m \in \Q, p > 0.$$
\end{cor}

{\it Proof.}
Use Conjecture \ref{dual-orth} and follow
Theorem \ref{CRIT-FILT}.\sq
}

Conjecture \ref{dual-orth} for $\mathfrak g$ of type 
$\mathsf{ADE}$ were essentially proved in \cite{FKM}. It 
yields an analogue of Theorem \ref{CRIT-FILT} for 
$\{ \mathsf D_b\}_b$ (instead of $\{ \mathbb D_b\}_b$).
Taking them into account, our ``explanation" of the 
remaining two decompositions in Theorem \ref{MAINTHM} 
(but not the exact formulas for the coefficients there)
is as follows:

\begin{conj}\label{CONJ:ANY}
Let $\Lambda$ be a level-one fundamental weight. Then
for any $b,c 
\in P$, $m \in \mathbb Q$ \,and $p > 0$, we have :
\begin{align*}
\mathrm{Ext}^p_{(\widetilde{\mathfrak b}^-, 
\widetilde{\mathfrak h})} 
( D_b \otimes L ( \Lambda ),\mathring{\mathbb D}_c^{\vee} \otimes 
\C_{\Lambda_0 + m \delta} ) & = \{ 0 \}, \\
\mathrm{Ext}^p_{(\widetilde{\mathfrak b}^-, 
\widetilde{\mathfrak h})} 
( \mathring{\mathbb D}_b \otimes L ( \Lambda ),  
D_c^{\vee} 
\otimes \C_{\Lambda_0 + m \delta} ) & = \{ 0 \}.
\end{align*}
\end{conj}
\vskip -1.1cm \sq

The equalities here
are equivalent to each other. 
The first corresponds to the decomposition in  
(\ref{pggbar}), the second interprets (\ref{pggdag}). 
We follow
the proof of Corollary \ref{MODULEMAIN-W} 
and use the relations:
$$\mathsf{gch} \, D_b ^{\vee} = q^{-b^2/2} w_0 
(\overline{E}_{b^{\iota}}),\ 
\hskip 3mm \mathsf{gch} \, \mathring{\mathbb D}_b = q^{b^2/2} 
w_0 (E^{\dag *}_{b^{\iota}}), \where b \in P.$$
See Proposition \ref{D-CHARS}, Corollary 
\ref{E-D-DAG}, and Conjecture \ref{dual-orth}.

\comment{
\extra{
YOU DID NOT USED $\mathbb D$ in the paper without mathring;
please mentioned this when $\mathring{ \mathbb D}$
occurs for the first time. Also, we never used LOCAL 
D-slices $\mathring D$  before this place, right? 
Please mention how they are defined,  when $end$ is 
discussed. --- For the former, I understand. I changed 
$\mathring{\mathbb D}$ to $\mathbb D$. For the latter, 
I do not understand what you mean... 
That is in Conjecture \ref{dual-orth}, and there seem 
nothing to be added. Also, I edited the latter half of the 
comments after Conjecture \ref{CONJ:ANY} as it looks more 
confusing than being helpful to me.}}

%\extra{REMIND HERE THE CONNECTION OF $D_b$, $\mathring{D}_b$
%with E-bar and E-dag}
%, and possibly comment greater what
%exactly follows from this Conjecture. BELOW FEEL FREE TO
%EDIT AND ADD IF ABSOLUTELY NECESSARY SOME REFS, BUT VERY FEW.
%WE ALREADY PROVIDED MANY IN THE PAPER. THIS IS JUST A REVIEW. 
%--- I think this is more than enough.}

\vskip 0.2cm
\subsubsection{\sf Conclusion}
As we discussed in the Introduction, the nonsymmetric
Rogers-Ramanujan-type identifies from our paper can be seen 
(upon some simplifications) as 
partitions of the symmetric ones, so they are formally
not ``brand new". However they are an important developments,  
since powerful DAHA tools can be used to study them,
which are non-existent in the symmetric theory. 
The DAHA technique of intertwiners is a major example, 
which is exactly the bridge between 
the $\overline{E},E^\dag$\~polynomials and Demazure modules. 

\vfil
In contrast to general  
$E$\~polynomials (for any $q,t$), the
coefficients of $\overline{E}_b,E^\dag_b$ are
{\em positive integers\,}, which alone is very remarkable.
This also holds for the coefficients of  
expansions of the products of theta-functions in terms 
of these polynomials. The 
interpretation of such integrality and positivity 
via the theory of level-one thick Demazure modules
(in the twisted setting) is what we do in this paper.
%including connections with the global and local Weyl 
%modules, level-one Demazure modules, and 
%Kirillov-Reshetikhin modules (closely related to
%Weyl and Demazure modules). Establishing the
%connection of the $E^\dag$\~polynomials with the 
%and further applications is the main purpose of the
%second half of our paper (in the twisted setting). 
 
%\extra{We didn't used them in this paper. 
%Actually, most of the (classical version of) 
%Kirillov-Reshetikhin modules are (special case of) 
%Demazure module} and 
%The theory of nonsymmetric polynomials,
%is a fundamental discovery. However these polynomials are 
%generally difficult to interpret within 
%the Lie theory (with some reservation concerning type $A$).
%Anyway the positivity of the coefficients does not hold in
%the general $q,t$\~setting. 
%\extra{Sorry. The latter statement is what I do not agree. 
%All of Vasserot's construction of polynomial representations 
%of DAHA, Haiman's realization of modified Macdonald polynomials, 
%and Braverman-Finkelberg-Shiraishi suggests that they naturally 
%appear geometrically. The problem is that they sometimes show up 
%as complexes, instead of a single sheaf.}
%From this viewpoint, the limits $t\to 0,\infty$ are really 
%exceptional. Their interpretation via the standard level-one
%Demazure modules (thin) and Demazure slices of thick
%modules appeared very clarifying.

\vfil
The restriction to the {\em level-one\,} Demazure modules
requires an explanation.
This case is exceptional in  Kac-Moody theory;
the corresponding (basic) integrable representations and
Demazure modules have many unique
features. The characters of the level-one 
Demazure modules (thin and thick) 
have important applications in the {\em classical\,} 
finite-dimensional Lie theory; in a sense,
they are ``no worse" than the characters of 
finite-dimensional simple Lie algebras.
\vfil

Our paper is mostly based on the following special feature of the 
Demazure level-one modules: they
constitute a remarkable 
families of {\em orthogonal\,} polynomials. 
More exactly, we use very much 
the orthogonality relation between the limits 
$t\to 0$ and  $t\to\infty$,
interpreted as some homological $Ext$\~duality in the
(based on prior works).  This is
the key for our
identification of the characters of Demazure slices with
$E^\dag$\~polynomials divided by their norms. 
%\vskip 0.2cm

\vfil
The connection between the limits $t\to 0$ and $ t\to\infty$,
is  quite non-trivial in the nonsymmetric
theory vs. the symmetric theory. Algebraically (combinatorially),
the $\overline{E}$\~polynomials
($t=0$) can be obtained from the 
$E^\dag$\~polynomials ($t=\infty$), but the latter are
significantly more involved than the former. The Kac-Moody
interpretation shed light on this asymmetric behavior; let
us touch it upon.
%\vskip 0.2cm
\vfil
  
The level-one Kac-Moody modules are naturally {\em unions\,}
(inductive limits) of the usual (thin) Demazure modules. 
On the other hand, the natural 
associate graded of the level-one Kac-Moody modules are direct
{\em sums\,} of Demazure slices (quotients of thick Demazure
modules).  
For instance, generally there are
no natural embeddings between neighboring (in the
Bruhat sense) Demazure slices in contrast
to usual (thin) Demazure modules.

At the level of graded characters of level-one
integrable modules, which are essentially
theta-functions, these two fundamental relations
become as follows.
The theta-function associated with a given 
root system is: $(a)$\ the limit of
$\overline{E}$-polynomials and $(b)$\ a certain 
Rogers-Ramanujan sum of $E^\dag$\~polynomials 
divided by their norms.
%\vskip 0.2cm

The nonsymmetric Rogers-Ramanujan sums are dual-purpose in
our paper. First, we use them 
to prove the coincidence of the characters of the
Demazure slices with the $E^\dag$\~polynomials divided
by their norm. This was expected in \cite{ChFB}, but
the representation-theoretical
tools at that time were insufficient to approach this problem. 
Second, we obtain that the 
multiplicities of the Demazure slices in tensor products
of level-one Kac-Moody modules can be found via the
DAHA-based ``numerical" machinery. Let us comment on the
latter.

The existence of
such a decomposition of tensor products of level-one
Kac-Moody modules in terms of Demazure
slices is quite a theorem, even without exact formulas
for the multiplicities. The latter are involved (unless in
the case of a single level-one Kac-Moody module); they are
remarkable {\em string functions\,} in some cases. 
Potentially, we can extend our approach (the expansion
via Demazure slices combined with the DAHA formulas)
to the decomposition of any integrable Kac-Moody modules,
but this is beyond this paper.
\vskip 0.2cm

To conclude, we think that this paper $(a)$\ 
significantly extends the
usage of nonsymmetric Macdonald polynomials
in the theory of (at least) integrable 
representations of affine Lie algebra, $(b)$\ clarifies 
the fundamental (but somewhat mysterious) role of the 
level-one Demazure characters in this theory and beyond,
$(c)$\ interprets the positivity phenomenon in the 
theory of non-symmetric Macdonald polynomials upon 
$t\to\infty$, which opens a road to their categorification
and that of the corresponding $q$\~Whittaker function. 
%\vfil

{\bf Acknowledgements.} The first author thanks Hiraku 
Nakajima and RIMS for the kind invitation in 2017, which resulted
in this paper, and MSRI for the invitation in 2018 (where this
work was continued).

\vskip -0.4cm 
\medskip
\bibliographystyle{unsrt}

\end{document}